\documentclass[mnsc,nonblindrev]{informs3}   
\usepackage{multirow}
\OneAndAHalfSpacedXI 

\usepackage{bbm}
\usepackage[ruled,lined,linesnumbered]{algorithm2e}
\usepackage{subcaption}
\usepackage[table,xcdraw]{xcolor}
\usepackage[flushleft]{threeparttable}
\usepackage{pdflscape}
\usepackage{natbib}
\usepackage{hhline}
 \bibpunct[, ]{(}{)}{,}{a}{}{,}%
 \def\bibfont{\small}%
\usepackage{hyperref}
\hypersetup{
    colorlinks=true,
    linkcolor=red,
    citecolor = blue,
    filecolor=magenta,      
    urlcolor=cyan,
}
\newcommand{\tabnote}[1]{%
    \vspace{0.5em} 
    \footnotesize\textit{Note:} #1\par 
}
\newcommand{\fignote}[1]{%
    \vspace{0.5em} 
    \footnotesize\textit{Note:} #1\par 
}
\usepackage{endnotes}

\TheoremsNumberedThrough     
\ECRepeatTheorems

\EquationsNumberedThrough    

\MANUSCRIPTNO{}

\begin{document}



\RUNTITLE{Multiclass Queue Scheduling Under Slowdown}

\TITLE{Multiclass Queue Scheduling Under Slowdown: An Approximate Dynamic Programming Approach }





%
\ARTICLEAUTHORS{%
\AUTHOR{Jing Dong}
\AFF{Decision, Risk, and Operations, Columbia Business School, \EMAIL{jing.dong@gsb.columbia.edu}} 
\AUTHOR{Berk G\"{o}rg\"{u}l\"{u}}
\AFF{DeGroote School of Business, McMaster University, \EMAIL{gorgulub@mcmaster.ca}}
\AUTHOR{Vahid Sarhangian}
\AFF{Department of Mechanical and Industrial Engineering, University of Toronto, \EMAIL{sarhangian@mie.utoronto.ca}}
} 

\ABSTRACT{%
In many service systems, especially those in healthcare, customer waiting times can result in increased service requirements. Such service slowdowns can significantly impact system performance. Therefore, it is important to properly account for their impact when designing scheduling policies. Scheduling under wait-dependent service times is challenging, especially when multiple customer classes are heterogeneously affected by waiting. In this work, we study scheduling policies in multiclass, multiserver queues with wait-dependent service slowdowns. We propose a simulation-based Approximate Dynamic Programming (ADP) algorithm to find close-to-optimal scheduling policies. The ADP algorithm (i) represents the policy using classifiers based on the index policy structure, (ii) leverages a coupling method to estimate the differences of the relative value functions directly, and (iii) uses adaptive sampling for efficient state-space exploration. Through extensive numerical experiments, we illustrate that the ADP algorithm generates close-to-optimal policies that outperform well-known benchmarks. 
We also provide insights into the structure of the optimal policy, which reveals an important trade-off between instantaneous cost reduction and preventing the system from reaching high-cost equilibria. Lastly, we conduct a case study on scheduling admissions into rehabilitation care to illustrate the effectiveness of the ADP algorithm in practice.}%


\KEYWORDS{multiclass queue, service slowdown, approximate dynamic programming, variance reduction}

\maketitle
\section{Introduction}\label{sec:introduction}

In the queueing literature, customers' service requirements are typically assumed to be independent of their waiting times. However, in many service systems, service times can be longer when the system is experiencing higher congestion. Such service \emph{slowdowns} are particularly common in healthcare settings. For example, \cite{huang2010impact} identify an association between long delays in admission to the emergency department (ED) and prolonged inpatient length-of-stay (LOS). \cite{chan2017impact} find that longer waiting times in admission to the Intensive Care Unit (ICU) can lead to longer ICU LOS. \cite{gorgulu2022} find that delays in transitioning from acute to rehabilitation care can result in longer rehabilitation LOS. These empirical studies also find that delays can have heterogeneous impacts on different types of patients.

Service slowdowns due to delays can significantly impact system performance and result in extended periods of high congestion and poor quality of service. Therefore, it is important to properly account for their impact when determining staffing and scheduling policies. Scheduling under wait-dependent service times is challenging, especially when there are multiple customer classes \citep{shitreatment}. Existing works mainly focus on performance analysis of single-class queues (\citealt{dong2015service,wu2019service,d2022m}), 
or scheduling in clearing systems (with no future arrivals) \citep{shitreatment, alidaee1999scheduling}. Moreover, simple scheduling policies such as the $c\mu$-rule \citep{van1995dynamic}, which perform well in the absence of slowdown, could perform very poorly when delays lead to increased service requirements.

 


In this paper, we propose a simulation-based Approximate Dynamic Programming (ADP) algorithm to find well-performing scheduling policies for multiclass multiserver queues with wait-dependent service slowdown. We start by considering a class of Markovian systems with preemption where the service rate is a function of the number of customers in the system. Our algorithmic development leverages specific structural properties of the problem through three main components: 
(i) a classifier to characterize the index policy structure; (ii) a coupling construction to estimate the value function differences directly; and (iii) an adaptive sampling algorithm for efficient state-space exploration. We conduct extensive numerical experiments and illustrate that our proposed algorithm finds near-optimal policies, outperforms other benchmarks, and scales well for large problem instances.

We further propose a heuristic to connect systems with wait-dependent service times to the proposed Markovian system with state-dependent service rates. This allows us to extend the policy derived from the Markovian system to systems with wait-dependent service times. In addition, we propose an extension of our ADP algorithm to systems where preemption is not allowed. We illustrate the effectiveness of these extensions in a case study on scheduling admissions to rehabilitation care. The case study is based on the empirical finding that admission delays have heterogeneous adverse effects on rehabilitation length-of-stay (service times) and other outcomes for different categories of patients. 

Our main contributions can be summarized as follows. 





\begin{enumerate}

    \item  We investigate scheduling policies for multiclass queues with class-dependent service slowdown. A notable phenomenon that can arise in systems experiencing service slowdown is meta-stability. We characterize the equilibrium points of the fluid model of the Markovian system with state-dependent service slowdown and show that it could have multiple stable equilibria, including some associated with long queues. Our results demonstrate that scheduling in systems prone to meta-stability requires carefully balancing the tradeoffs between the instantaneous cost reduction rate and avoiding entrapment in ``bad" equilibria.

    \item We proposed a simulation-based ADP algorithm to solve the scheduling problem. The algorithm incorporates several innovative methods to reduce computational time. These methods are more broadly applicable to other complex queueing control problem. Notably, by exploiting the specific structure of the problem, we propose a coupling method that directly estimates the value function differences, which bypasses some of the challenges in relative value function estimation. Compared to existing methods, the proposed coupling results in a significant reduction in simulation lengths as well as the variance of the estimates (up to 93\%) which collectively reduce computational cost.

    \item We provide insights into the structure of the optimal policy for scheduling queues with class-dependent service slowdowns. In contrast to simple queuing systems, where the optimal policy focuses on myopically maximizing the instantaneous rate of reducing the holding cost, in the presence of slowdown, the optimal policy aims to also prevent the system from being stuck at the high-cost equilibria, hence balancing two objectives. Therefore, the optimal policy is affected by the change in system load generated by the increased slowdown or changes in the arrival rates.
    \item  We conduct a case study on scheduling admissions to rehabilitation care using real data. In this setting, rehabilitation length-of-stays are affected heterogeneously by admission delays for Neuro/MSK and Medicine patients. We find that the ADP policy provides a $2.5$ days reduction in average waiting times compared to historical averages. This respectively translates to, on average, $2.17$ and $4.93$ increase in discharge Functional Independence Measure (FIM) (a standard metric for measuring rehabilitation outcomes taking values between 18 and 126) for Neuro/MSK and Medicine patients, respectively. Furthermore, compared to a naive first-come-first-served (FCFS) policy, using the ADP policy results in $10.76$ days reduction in average waiting times and respectively $25.40$ and $16.74$ point increases in the discharge FIM scores for Neuro/MSK and Medicine patients. 


    
\end{enumerate}

The rest of the paper is organized as follows. Section \ref{sec:lit_rev} discusses the related literature. Section \ref{sec:model_desc} describes the main model and its connection to the wait-time dependent model. Section \ref{sec:adp} presents our simulation-based ADP algorithm, and Section \ref{sec:nonpreemptive} extends it to the non-preemptive setting. In Section \ref{sec:num_experiments}, we conduct extensive numerical experiments for a two-class queue to illustrate the performance of our algorithm and generate insights into the structure of the optimal policy. Section \ref{sec:fiveclass} presents further numerical experiments with five classes of customers to demonstrate the scalability of the algorithm. Section \ref{sec:case_study} presents a case study on scheduling patient admission to rehabilitation care. Section \ref{sec:conclusion} concludes our paper with a discussion of future research directions.

\section{Literature Review}\label{sec:lit_rev}
Our work relates to two main streams of literature: (i) scheduling in queueing systems, especially those with wait-dependent service slowdown, and (ii) approximate dynamic programming, especially its application to queueing systems. We next provide a brief summary of some related works.

\textit{Scheduling in queueing systems.}
Scheduling multiclass queues is extensively studied in the Operations Research literature. In most of the previous works, customers' service requirements are assumed to be independent of the load of the system or their waiting times. When service times are exogenous, the $c\mu$ policy (or generalizations of it) which prioritizes customers with a higher cost of waiting and a shorter service requirement is known to achieve near-optimal performance \citep{van1995dynamic,mandelbaum2004scheduling}. \cite{atar2004scheduling, atar2010cmu}, and \cite{kim2018dynamic} study the scheduling problem where customers can abandon the system if waited for too long. For example, \cite{atar2010cmu} show that when the patience time is exponentially distributed and the system is overloaded, the $c \mu / \theta$ policy which prioritizes customers with a higher cost of waiting, a shorter service requirement, and a smaller abandonment rate achieves near-optimal performance. \cite{stolyar2004maxweight} study the MaxWeight policy where the size of the queue also enters the priority decision. \cite{dai2005maximum} studies the MaxPressure policy, which is the same as the MaxWeight policy in parallel server systems, but takes into account future flow dynamics when dealing with a network of queues. The MaxPressure policy is shown to be throughput optimal and achieves good performance cost-wise in many applications.

Our work considers the case where the service rate depends on the load of the system. Scheduling under load (wait)-dependent service times is challenging, especially when there are multiple customer classes that are affected differently by the system load (wait). Myopic policies where one gives priority to the class with the highest cost-reduction rate (like the c$\mu$ rule) may be highly sub-optimal. \cite{hu2022optimal} examine multiclass queues where customers can transition from one class to another if waited for too long. They show that when the system is congested, a more forward-looking modified $c\mu/\theta$ policy that takes the class transition into account should be applied. Most other works on queues with load (wait)-dependent service requirements focus on queues with a single class of customers. Thus, there is no multi-class scheduling decision involved (see, e.g., \citealt{dong2015service,chan2017impact,selen2016snowball}). \cite{wu2019service} study queueing systems where customers' patience times and service times are correlated. Note that the correlation can be exogenous, i.e., customers with longer service requirements are more patient, or endogenous, i.e., waiting can generate more service needs \citep{wu2021service}.

Scheduling problems with a fixed number of jobs (no arrivals) and wait-dependent processing times have been studied in machine scheduling (for a single machine). \cite{browne1990scheduling} consider a model where the processing times are linearly increasing in the wait. They show that sorting the jobs based on $\mathbb{E}[S_i]/\alpha_i$ minimizes the makespan, where $S_i$ is the initial service time and $\alpha_i$ is the deterioration rate. Following \cite{browne1990scheduling}, other studies have considered different slowdown functions (e.g., piecewise linear, convex) and objectives (e.g., total flowtime, total tardiness); see \cite{alidaee1999scheduling} and \cite{cheng2004concise} for comprehensive reviews of these works. Recently, motivated by treatment planning in mass casualty events, \cite{shitreatment} studies the scheduling problem where jobs can have heterogeneous slowdown functions. The goal is to minimize the makespan, among all schedules that minimize the number of overdue jobs. The key difference between our work and these scheduling works is that we consider a queueing system with customer arrivals instead of a clearing system. As we will demonstrate later, the arrival process can have a major effect on the optimal scheduling policy.

\textit{Approximate dynamic programming (ADP) for scheduling in queueing systems.} Due to the complicated system dynamics under service slowdown, analytical characterization of the optimal scheduling policy is challenging. As such, in this work, we adopt an ADP approach to find near-optimal policies. ADP algorithms have been developed for various sequential decision making problems in Operations Research, including queueing control problems. A common approach in the literature involves approximating the value function using a parametric function of the system state \citep{powell2007approximate}. 

There are two main ways to leverage the value function approximation.  
The first one utilizes the linear programming (LP) formulation of the MDP (e.g., \citealt{de2003linear,adelman2008relaxations,veatch2015approximate,bhat2023}). This approach requires a linear parameterization of the value function, which may constrain its capacity to approximate complex value functions -- a challenge that arises in our problems due to service slowdowns. In addition, the corresponding LP can be very large in scale. Further algorithmic developments are required to solve the LP efficiently \citep{topaloglu2009using}. 
The second line of literature focuses on dynamic programming-based algorithms, such as value iteration and policy iteration (see, e.g., \citealt{maxwell2013tuning,  moallemi2008approximate, dai2019inpatient}). In general, the performance of the algorithm depends on the choice of the basis functions for value function approximation. \cite{moallemi2008approximate}, \cite{chen2009approximate}, and \cite{dai2019inpatient} propose to choose the basis function based on the ``limiting" fluid queue. More general approximations, such as neural network-based approximations can be used as well, see, e.g., \cite{van1997neuro}. \cite{liu2022rl} consider a Reinforcement Learning (RL)-based approach to solve queueing network control problems. 
\cite{dai2022queueing}
investigate applying deep RL algorithms for queueing control problems with long-run average cost formulation. 
See \cite{walton2021learning} for a review of recent developments on applying RL to queueing systems.


Our proposed algorithm is based on (approximate) policy iteration. Compared to the previous literature, we approximate the policy directly rather than the value function. This is achieved by exploiting the special structure of queueing scheduling problems. As we illustrate in the paper, approximating the policy can be a much easier task than approximating the value function, which is quite complicated in our setting. Parameterizing the policy directly has been used in the policy gradient method \citep{sutton1999policy, kakade2001natural}. Compared to the general policy gradient literature, especially \cite{dai2022queueing}, 
our approach utilizes more domain knowledge of the structure of the underlying problem. In particular, we find that the optimal scheduling policy follows an index-based priority rule. Thus, the policy can be approximated by a classifier that classifies different states into different priority orders. 

A key challenge in applying policy iteration or policy gradient in queueing applications is to accurately estimate the (relative) value function or advantage function, especially under the long-run average cost formulation, due to the high variance and slow mixing \citep{asmussen1992queueing,dai2022queueing}. Various variance reduction techniques have been developed to address the challenge \citep{henderson2002approximating,henderson2020variance}. Since only the difference of the relative value function is relevant in characterizing the optimal scheduling policy in queueing scheduling problems, we propose a coupling construction to estimate the differences directly. The coupling construction leads to a significant reduction in the variance of the estimate and the simulation time compared to the regeneration method used in \cite{dai2022queueing}. To further improve the computational efficiency, we develop an adaptive sampling approach, which allows us to dynamically allocate the sampling budget to different states. This adaptive sampling approach is related to the ranking and selection literature (see \citealt{chen2015ranking,hong2021review} for reviews), but the specific application and the objective considered here is different. 

Lastly, it has been well established in the literature that a good initial policy (initialization) is important in learning the optimal control in queueing systems \citep{chen1999value, liu2022rl, dai2022queueing}. To derive a good initial policy, we develop a fluid-based policy iteration algorithm. The fluid-based approximation is similar to that of \cite{chen1999value} but with two important distinctions. First, we consider a system with service slowdown where the fluid model can have more than one equilibrium point. Second, due to the complexities of the fluid control problem, it is not possible to obtain an analytical solution.

\section{Model Description}\label{sec:model_desc}
Our main model is a multiclass, multiserver queue with state-dependent service rates. We start by allowing preemption in scheduling and later extend to the non-preemptive case (see Section \ref{sec:nonpreemptive}). We also discuss how this model serves as an approximation for queues with wait-dependent service times. The system consists of $I$ distinct classes of customers and $C$ identical servers; each server can serve only one customer at a time. Class $i$ customers arrive according to a Poisson process with rate $\lambda_i$, $i=1,\dots, I$. The service requirements are exponentially distributed with rates that are class-dependent functions of the system state. In particular, the service rate for class $i$ customers at time $t$ takes the form $\bar f_i(X_i(t))$, where $X_i(t)$ is the number of class $i$ customers in the system at time $t$ and $\bar f_i$ is a non-negative decreasing function. The system has a finite capacity $\kappa_i$ for class $i$ customers. Specifically, if a class $i$ customer arrives at the system when the system already has $\kappa_i$ class $i$ customers (including both the customers in service and waiting to be served), this arriving customer is blocked. 

Let $\{X(t)=(X_1(t), \dots, X_I(t)); t\geq 0\}$ denote the number of customers in the system from each class, and $\{Z(t)=(Z_1(t), \dots, Z_I(t)); t\geq 0\}$ denote the number of customers in service from each class. Note that $X_i(t)-Z_i(t)$ is the number of class $i$ customers waiting in the queue, and $X(t)$ takes values in $\mathcal{X}:=\{0,1,\dots, \kappa_1\}\times \{0,1,\dots, \kappa_2\} \times \dots \times \{0,1,\dots, \kappa_I\}$. The process $Z(t)$ is determined by the scheduling policy $\pi$, which is assumed to be non-anticipatory, i.e., it can only depend on the current or past states of the system. We assume preemption is allowed. Specifically, $Z(t)$ for an admissible policy $\pi$ and $t \geq 0$ satisfies,
\begin{equation}
     Z_i(t) \leq X_i(t), \forall i \in \{1, \dots, I\}\ \text{and}\ \sum\limits_{i=1}^{I} Z_i(t) \leq C,
\end{equation}

Let $h_i$ denote the cost of holding a class $i$ customer for one unit of time in the system and $b_i$ denote the cost of blocking a class $i$ customer.
Our goal is to minimize the long-run average cost of holding customers waiting and blocking incoming customers. Let $B(t)=(B_1(t), \dots, B_I(t))$ denote the number of customers blocked by time $t$ from each class. Since $X(t), Z(t), B(t)$ depend on the scheduling policy $\pi$, we write them as $X^{\pi}(t), Z^{\pi}(t), B^{\pi}(t)$ when we want to mark the dependence explicitly. (We sometimes suppress the dependence on the scheduling policy when it is understood from the context). We define the expected total cost incurred over $[0,T]$ as,
\[
\Gamma^{\pi}(T) = \mathbb{E}\left[\int_{t=0}^{T}\sum\limits_{i=1}^{I} h_i X_i^{\pi}(t)dt + \sum\limits_{i=1}^{I}  b_i B^{\pi}_i(T) \right].
\]
Our objective is to minimize the long-run average cost, which can be expressed as
\begin{align} \label{eq:main}
    \min\limits_{\pi \in \Pi}\ \limsup_{T \rightarrow \infty}\  \frac{1}{T} \Gamma^{\pi}(T),
\end{align}
where $\Pi$ is the set of all admissible policies. 

We use uniformization to formulate the problem as a discrete-time Markov decision process \citep{puterman2014markov}. As the state space is finite and the action space at each state is compact, it is without loss of generality to only consider deterministic Markovian policies \citep{puterman2014markov}. In particular, at decision epoch $t$, the scheduling policy can be viewed as a mapping from the state of the Markov chain, $X(t)$, to the allocation of the servers, $Z(t)$. Note that the preemption assumption implies
that an optimal policy does not allow deliberate idleness.

\subsection{Connection to Wait-Dependent Service Times}
For our primary motivating applications, the scheduling problem may deviate from the proposed model in two important ways. First, the service time depends on the waiting time rather than the state of the system. Second, preemption may not be feasible. However, incorporating these features substantially increases the complexity of the problem. We start with the simpler model to gain insights into the algorithmic development. We discuss the extension of our algorithms to the non-preemptive case in Section \ref{sec:nonpreemptive}. In what follows, we propose an approach to connect the wait-dependent service time model to the simpler state-dependent service time model.


For the model with wait-dependent service times, let $W_{ij}$ denote the waiting time of the $j$-th arriving customer from class $i$. The service of time for the customer $(i,j)$ follows a general distribution with mean $g_i(W_{ij})$. In this case, to have a Markovian system descriptor, we need to keep track of the waiting times of customers in the queue and the amount of time each job has spent in service. This renders the corresponding MDP too complicated to solve even when there are only two classes of customers and a single server.

The simpler model with a state-dependent service time can be viewed as an approximation to the model with wait-dependent service times. 
Let $Q_i(t)$ denote the queue length for class $i$ at time $t$. We also denote $A_i(s,t)$ as the number of class $i$ arrivals on $(s,t]$ and $R_i(t)$ as the total residual waiting time of all class $i$ customers in the system at time $t$. 
Based on the sample-path construction of Little's law (Theorem 2 in \citealt{kim2013statistical}), we have 
\[
\frac{1}{A_i(t,t+\delta)}\sum_{j=Q_i(t)+1}^{Q_i(t)+A_i(t,t+\delta)}W_{ij}
=\frac{\delta}{A_i(t,t+\delta)}\frac{1}{\delta}\int_{t}^{t+\delta}Q_i(s)ds
+\frac{R_i(t)-R_i(t+\delta)}{A_i(t,t+\delta)}.
\]
Note that if $R_i(t)\approx R_i(t+\delta)$, then we have 
\[
W_i(t) \approx \frac{Q_i(t)}{\lambda_i},
\]
where $W_i(t)$ denotes the waiting time of a newly arrived class $i$ customer at time $t$. 
Because we are interested in the number in system rather than the queue length, we utilize the analogous relations between (waiting times, queue length) and (time in system, number in system). Thus, for a given number in system $x_i$, and the wait-dependent service time function $g_i(w_i)$, we propose the following heuristic to find the waiting time $w_i$ such that the time in system ($w_i + g_i(w_i)$) is equal to $x_i /\lambda_i$,
\begin{equation} \label{eq:optimize}
    \tilde{w_i}(x_i) := \argmin\limits_{w_i\geq 0} |w_i + g_i(w_i) - x_i / \lambda_i|.
\end{equation}
Hence, we set the state-dependent service rate as $\bar f_i(x_i)= 1/\tilde{w_i}(x_i)$. We illustrate the performance of this heuristic approach in the case study in Section \ref{sec:case_study}.

\subsection{Fluid Approximation and Equilibrium Analysis}\label{sec:fluid_approx}
In this section, we study a deterministic fluid approximation of the main model. The fluid model provides insights into the long-run behavior of the system under load-dependent service slowdown.
Let $\Bar{X}(t)=(\Bar{X}_1(t), \dots, \Bar{X}_I(t))$, where $\Bar{X}_i(t)$ denote the amount of class $i$ fluid in the system at time $t$. We also write $\Bar{Z}(t)=(\Bar{Z}_i(t), \dots, \Bar{Z}_I(t))$, where $\Bar{Z}_i(t)$ denote the amount of capacity allocated to serve class $i$ fluid at time $t$. The fluid model satisfies the following Ordinary Differential Equation (ODE) with right-hand-side discontinuity. For $i \in \{1,\dots, I\}$,
\begin{align}\label{eq:fluid_dynamics}
\dot{\Bar{X}}_i(t)=\lambda_i\mathbbm{1}\{\Bar{X}_i(t)<\kappa_i\} - \Bar{f}_i(\Bar{X}_i(t))\Bar{Z}_i(t),
\end{align}
where $\Bar{Z}_i(t)$ is determined by the scheduling policy satisfying,
\[
\Bar{Z}_i(t)\leq \Bar{X}_i(t),\ \forall i \in \{1, \dots, I\}\ \text{and}\ \sum\limits_{i=1}^{I}\Bar{Z}_i(t) \leq C.
\]

Note that the ODE \eqref{eq:fluid_dynamics} that describes the dynamics of our fluid model is discontinuous at $\Bar{X}_i(t) = \kappa_i$, but continuous in time. Thus, we first establish that the ODE is well-defined and a solution in the sense of Filippov \citep{filippov2013differential} exists. 

\begin{proposition}\label{prop:existance_fluid}
    For any initial condition $\Bar{X}(0)$, there exists a solution to the problem defined by \eqref{eq:fluid_dynamics}.
\end{proposition}



Next, we examine the limiting behavior of the fluid model as $t \rightarrow \infty$. Let $\dot{\bar{X}}(t)=\varphi^{\bar{\pi}}(\bar{X}(t))$, $\bar{X}(t)\in \mathbbm{R}^{I}_{+}$ and 
$\Phi^{\bar{\pi}}(\bar{x}_0, t)$ denote the value of the fluid trajectory at time $t$ starting from $\Bar{X}(0) =\bar{x}_0$ under policy $\bar{\pi}$. We use the following definitions of equilibria.

\begin{definition}{(Equilibrium).} A state $\Bar{x}_e$ is defined as an equilibrium point of the system \eqref{eq:fluid_dynamics} under policy $\bar{\pi}$ if
$\Phi^{\bar{\pi}}(\Bar{x}_e, t) = \Bar{x}_e$, $\forall t\geq0$.
\end{definition}

We distinguish between two types of equilibrium points: those that satisfy $\varphi(\Bar{x}_e)=0$, representing stationary points of the ODE, and those that emerge within regions of discontinuity, which we refer to as \emph{pseudo equilibria} following \cite{bernardo2008piecewise}.

\begin{definition}{(Locally 
Asymptotically Stable Equilibrium).} An equilibrium point $\Bar{x}_e$ is locally asymptotically stable 
if there exists $\epsilon>0$ such that for any initial point $\bar x_0$ satisfying $||\Bar{x}_0 - \Bar{x}_e|| < \epsilon$, $\lim\limits_{t \rightarrow \infty} \Phi^{\bar{\pi}}(\Bar{x}_0, t) = \Bar{x}_e$.
\end{definition}

We refer to a locally asymptotically stable stationary point as a stable equilibrium. Similarly, we refer to a locally asymptotically stable pseudo equilibrium as a stable pseudo equilibrium.

\begin{definition}{(Meta-stability).} A stochastic system is meta-stable if its corresponding fluid limits have multiple stable (pseudo) equilibria.  
\end{definition}

When there is load-dependent service slowdown, the fluid model can exhibit meta-stability \citep{dong2022metastability}. To illustrate this, we consider a simple two-class example with linear slowdown functions: \begin{equation} \label{eq:slowdown1}
    \bar{f}_i(x_i) = \mu_i - a_i x_i, 
\end{equation}
where $a_i$ is the slowdown rate and $\mu_i$ is the maximum service rate for class $i$ customers. We make the following assumptions on the system parameters: (i) $\lambda_i / (\mu_i - a_i C) \leq C$, (ii) $\mu_i > a_i \kappa_i$, (iii) $\kappa_i \geq C$. The first assumption guarantees that when there are fewer class $i$ customers than the number of servers, the effective service rate is enough to keep up with class $i$ arrivals. The second assumption states that the effective service rate does not go to zero. This can also be guaranteed by truncating the effective service rate at some positive lower bound, but for simplicity, we impose this assumption. The last assumption states that the system size for class $i$ customers should be at least as large as the number of servers. Furthermore, for illustration purposes, we assume a simple strict priority policy that prioritizes class 1 customers. We next characterize the corresponding equilibrium points.


\begin{proposition}\label{prop:fluid_characterize}
For the two-class system with service rate function \eqref{eq:slowdown1}, assuming $\lambda_i / (\mu_i - a_i C) \leq C$, $\mu_i \geq a_i \kappa_i$, $\kappa_i \geq C$,  and under strict priority to class 1, the ODE \eqref{eq:fluid_dynamics} has one stable equilibrium and two stable pseudo equilibria:
\begin{align*}
    \Bar{x}_e \in 
    \Bigg\{ \left(\frac{\mu_1 - \sqrt{\mu_1^2 - 4 a_1 \lambda_1}}{2a_1},\frac{\mu_2 - \sqrt{\mu_2^2 - 4 a_2 \lambda_2}}{2a_2}\right),\ \left(\frac{\mu_1 - \sqrt{\mu_1^2 - 4 a_1 \lambda_1}}{2a_1},\kappa_2\right), \left(\kappa_1,\kappa_2\right) \Bigg\}.    
\end{align*}

\end{proposition}

\begin{figure}[!h]
     \centering
     \begin{subfigure}[b]{0.47\textwidth}
         \centering
         \includegraphics[width=0.67\textwidth]{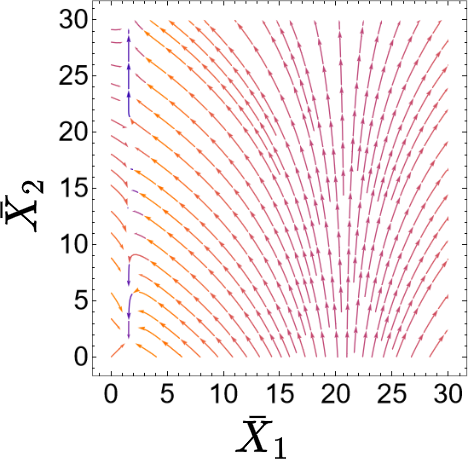}
         \caption{Vector field under meta-stability.}
         \label{fig:vector_field}
     \end{subfigure}\hfill
    \begin{subfigure}[b]{0.47\textwidth}
         \centering
         \includegraphics[width=.96\textwidth]{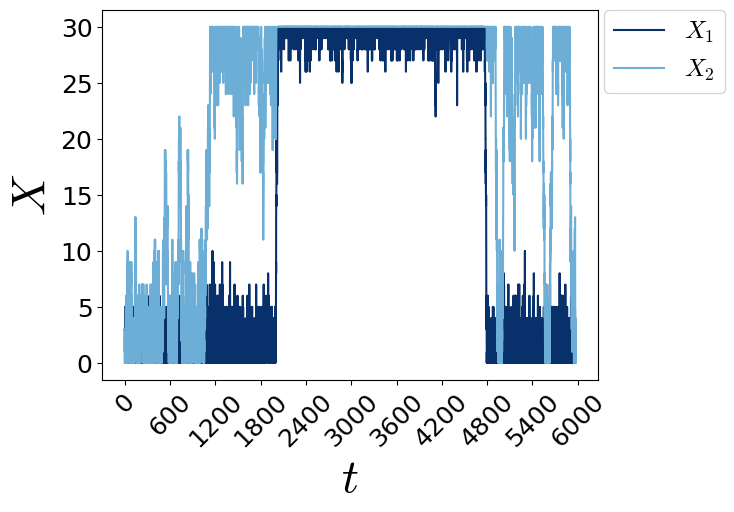}
         \caption{Sample paths under meta-stability.}
         \label{fig:stochastic_trajectory_fluid}
     \end{subfigure}
     \caption{Vector field of the fluid model and an example of sample paths for the stochastic system under meta-stability with parameters $\lambda_1 = \lambda_2=1.5$, $\mu_1=\mu_2=1$, $a_1=0.03$, $a_2=0.02$, $\kappa_1=\kappa_2= 30$, $C=4$.}
     \label{fig:fluid_flow_stochastic}
     \end{figure}

Proposition \ref{prop:fluid_characterize} states that under a simple strict priority policy, the fluid model \eqref{eq:fluid_dynamics} exhibits meta-stability (a more complicated case with four equilibrium points is provided in Appendix \ref{sec:additional_vector_field_example}). Among the three equilibrium points characterized in Proposition \ref{prop:fluid_characterize}, we refer to the first equilibrium as ``good," because the equilibrium is associated with zero fluid queue, and the second and third equilibria as ``bad," because they are associated with a positive queue. Figure \ref{fig:vector_field} illustrates the vector field for a system with meta-stability. Meta-stability in the fluid model implies that the initial state can determine which equilibrium point the fluid trajectory converges to. With random fluctuations, the stochastic system can evolve to either equilibrium as illustrated in Figure \ref{fig:stochastic_trajectory_fluid}. Therefore, developing appropriate control policies that prevent the system from getting stuck in the bad equilibrium regions is critical.



\section{Approximate Dynamic Programming}\label{sec:adp}
In this section, we develop an ADP approach to solve the scheduling problem \eqref{eq:main}. We start by providing an overview of the key ideas. 

Let $\Lambda=\sum_{i=1}^{I}\lambda_i+ C\max\limits_{1\leq i\leq I}\bar f_i(0)$ denote the maximum possible transition rate of the system, i.e., the uniformization factor. We also define $\mathcal{Z}(x):= \{z: z_i \leq x_i\ \forall i, \sum_{i=1}^I z_i\leq C \}$ and write $\mathbf{e}_i$ as the unit vector with the $i$-th element equal to 1.
For the long-run average cost problem, we need to solve the Bellman equation,
\begin{align} \label{eq:bellman}
     v^*(x) &= \min\limits_{z \in \mathcal{Z}(x)} \frac{1}{\Lambda}\left\{\sum\limits_{i=1}^{I} h_i x_i + \lambda_i b_i \mathbbm{1}_{\{x_i = \kappa_i\}} - \gamma^* + \sum\limits_{i=1}^{I} \lambda_i\mathbbm{1}_{\{x_i < \kappa_i\}} v^{*}(x + \mathbf{e}_i)
    + \sum\limits_{i=1}^{I} z_i\Bar{f}_i(x_i) v^{*}(x-\mathbf{e}_i) 
    \nonumber \right.\\ 
    &\left. + \Bigl(\Lambda - \sum\limits_{i=1}^{I} \lambda_i\mathbbm{1}_{\{x_i < \kappa_i\}} - \Bar{f}_i(x_i) z_i \Bigr) v^*(x) \right\}, \quad \forall x\in \mathcal{X}, &
\end{align}
where $\gamma^*$ is the optimal long-run average cost and $v^*(\cdot)$ is the optimal relative value function. 
Note that if we view the right-hand side of \eqref{eq:bellman} as an optimization problem with respect to $z$, then the optimal action $z^*$ solves the following linear program:
\begin{align}\label{eq:LP}
    \max_{z} \sum\limits_{i=1}^{I} \Bar{f}_i(x_i)D^{*}_i(x) z_i ~~
    \mbox{ s.t. } z\in \mathcal{Z}(x),
\end{align}
where $D^{*}_i(x) := v^{*}(x)-v^{*}(x-\mathbf{e}_i)$. We make a few observations from the linear program \eqref{eq:LP}. 

First, the optimal solution follows a strict priority rule where the class with a higher $\Bar{f}_i(x_i)D^{*}_i(x)$-value receives a higher priority. This suggests that to characterize the scheduling policy, we can divide the state space into different regions, each corresponding to a different prioritization order, i.e., different order of the values of $\Bar{f}_i(x_i)D^{*}_i(x)$'s. This motivates us to use a classifier to approximate the scheduling policy. We discuss this in more detail in Section \ref{sec:adp_simplified}.

Second, the policy depends on $D^{*}_i(\cdot)$, the differences of the relative value functions, rather than the relative value function itself. Given a policy $\pi$, the corresponding relative value function can be characterized as 
\begin{equation}\label{eq:value_pi}
v^{\pi}(x)=\mathbb{E}\left[\sum_{k=0}^{\infty} \left(c(X^{\pi}(t_k))-\gamma^*\right)| X^{\pi}(t_0)=x \right],
\end{equation}
where $c(x)=\sum_{i=1}^{I}\frac{h_i}{\Lambda}x_i + \frac{\lambda_i}{\Lambda} b_i\mathbbm{1}_{\{x_i = \kappa_i\}}$.
Estimating $v^{\pi}(x)$ can be challenging due to its large variance \citep{dai2022queueing}. The value function difference $D_i^{\pi}(x)$ can in contrast be estimated more efficiently. Note that
\[
D_i^{\pi}(x)=\mathbb{E}\left[\sum_{k=0}^{\infty} c(X^{\pi}(t_k))-c(\tilde X^{\pi}(t_k))  |X^{\pi}(t_0)=x, \tilde X^{\pi}(t_0)=x-\mathbf{e}_i\right]
\]
where $X^{\pi}$ and $\tilde X^{\pi}$ are two systems under policy $\pi$ starting from $x$ and $x-\mathbf{e}_i$ respectively. Under appropriate coupling construction between $X^{\pi}$ and $\tilde X^{\pi}$, we only need to simulate the two systems until they coincide, after which the cost difference between the two will be zero, leading to substantial simulation time and variance reduction; see Section \ref{sec:coupling} for more details.

\begin{figure}[!h]
     \centering
     \begin{subfigure}[b]{0.47\textwidth}
         \centering
         \includegraphics[width=.75\textwidth]{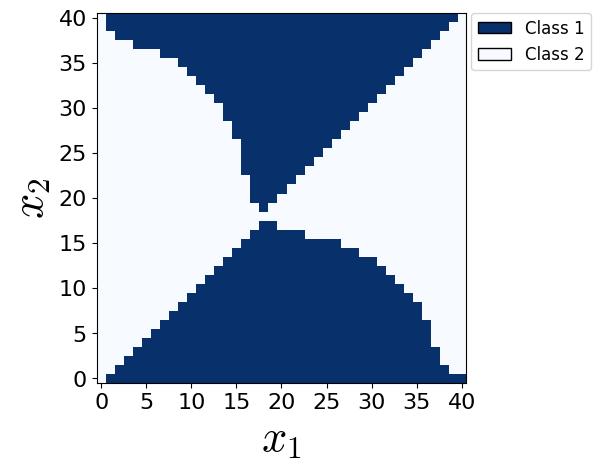}
         \caption{Optimal Policy}
     \end{subfigure}
    \begin{subfigure}[b]{0.47\textwidth}
         \centering
         \includegraphics[width=.75\textwidth]{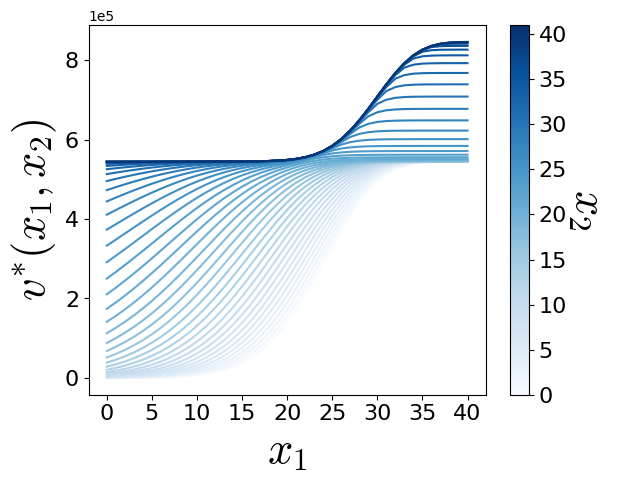}
         \caption{$v^*(x_1,x_2)$}
     \end{subfigure}
     \caption{A visual illustration of optimal policy and relative value function under the parameters $\lambda = (0.3,0.3), \mu = (0.9,0.9), a = (0.02,0.02), \kappa = (40,40), h = (1,1), b = (0,0), \bar{f}_i(x)= \mu_i - a_i x_i$.}
     \label{fig:new_case}
     \end{figure}

To demonstrate the observations discussed above, we consider a simple two-class example with linear slowdown, i.e., $ \bar{f}_i(x)= \mu_i - a_i x_i$ with $\mu_i>a_i\kappa_i$.
Figure \ref{fig:new_case} plots (a) the optimal policy and (b) some optimal relative value function values for a symmetric system with  $\lambda = (0.3,0.3), \mu = (0.9,0.9), a = (0.02,0.02), \kappa = (40,40), h = (1,1), b = (0,0)$.
Note that in this two-class example, there are only two possible priorities: prioritizing class 1 over class 2 and vice versa. Figure \ref{fig:new_case} (a) shows how the state space can be divided into different regions according to the optimal policy. The darker region is where class 1 is prioritized and the lighter region is where class 2 is prioritized. We also note that the relative value function as depicted in Figure \ref{fig:new_case} (b) takes various shapes and can be hard to approximate using a parameterized family of functions, in addition to the high variance incurred in simulation-based estimation.



Our proposed ADP algorithm is based on a simulation-based policy iteration approach. In what follows, we first introduce the policy iteration algorithm and then provide more details about the key features of our algorithm: (i) representing the policy using a classifier; (ii) estimating the value function differences utilizing a coupling construction; and (iii) using adaptive sampling to learn the policy at different states, which is then used to train the policy classifier. We also propose using the optimal fluid policy as a good initial policy. Section \ref{sec:fluid_ADP} provides more details as to how to find the optimal fluid policy.

\subsection{Sampling-Based ADP: Approximate Policy Iteration}\label{sec:adp_simplified}
A classic approach to solving \eqref{eq:bellman} is policy iteration. Given a policy $\pi_n$ at the $n$-th iteration, it involves two key steps:
\begin{enumerate}
\item Policy evaluation: obtain the relative value function under the current policy, i.e., $v^{\pi_n}(x)$ as defined in \eqref{eq:value_pi} for all $x\in \mathcal{X}$.
\item Policy improvement: obtain the updated policy $\pi_{n+1}$ by solving
\[
\min\limits_{z \in \mathcal{Z}(x)} \sum\limits_{i=1}^{I} z_i\Bar{f}_i(x_i) v^{\pi_n}(x-\mathbf{e}_i) 
\nonumber + \Bigl(\Lambda - \sum\limits_{i=1}^{I} \lambda_i\mathbbm{1}_{\{x_i < \kappa_i\}} - \Bar{f}_i(x_i) z_i \Bigr) v^{\pi_n}(x), \quad \forall x\in \mathcal{X},
\]
which can equivalently be written as
$\max\limits_{z \in \mathcal{Z}(x)} \Bar{f}_i(x_i)D^{\pi_n}(x)z_i$.
\end{enumerate}
Steps 1 and 2 can be very computationally intensive or even prohibitive in practice when $|\mathcal{X}|$ is large. To address the curse of dimensionality, we propose the following simulation-based approximate policy iteration. Given a policy $\hat\pi_n$ at the $n$-th iteration:
\begin{enumerate}
\item Use simulation to estimate $D_i^{\hat\pi_n}(x)$, $i=1,\dots, I$, at a selected set of states, i.e., for $x\in \mathcal{X}_0$ where $\mathcal{X}_0$ is a small subset of $\mathcal{X}$. We simulate $D_i^{\hat\pi_n}(x)$, $x\in \mathcal{X}_0$, under a coupling construction to increase the estimation efficiency (see Algorithm \ref{algo:coupling}).
\item Update the policy by training a classifier that classifies each state $x$ to an order of the indexes $\bar f_1(x_1)D_1^{\hat\pi_n}(x), \bar f_2(x_2)D_2^{\hat\pi_n}(x), ..., \bar f_I(x_I)D_I^{\hat\pi_n}(x)$. For example, when $I=2$, we will classify each state as either $\bar f_1(x_1)D_1^{\hat\pi_n}(x)\geq \bar f_2(x_2)D_2^{\hat\pi_n}(x)$, in which case class 1 is prioritized over class 2, or $\bar f_1(x_1)D_1^{\hat\pi_n}(x)<\bar f_2(x_2)D_2^{\hat\pi_n}(x)$, in which case class 2 is prioritized over class 1.
\end{enumerate}

\begin{figure}[]
     \centering
     \begin{subfigure}[b]{0.285\textwidth}
         \centering
         \includegraphics[width=1.03\textwidth]{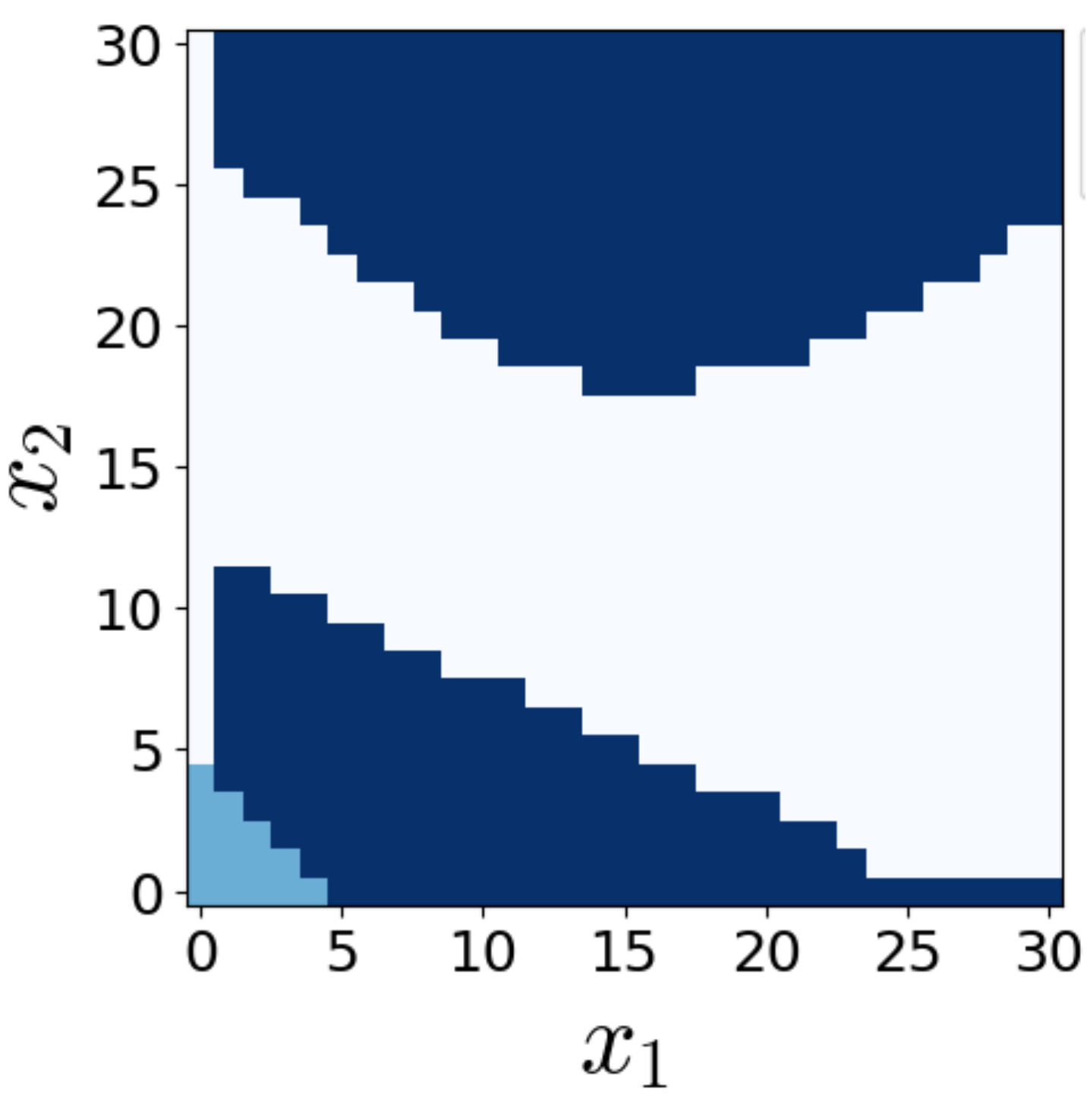}
         \caption{Optimal Policy}
     \end{subfigure}
     \begin{subfigure}[b]{0.285\textwidth}
         \centering
         \includegraphics[width=1.03\textwidth]{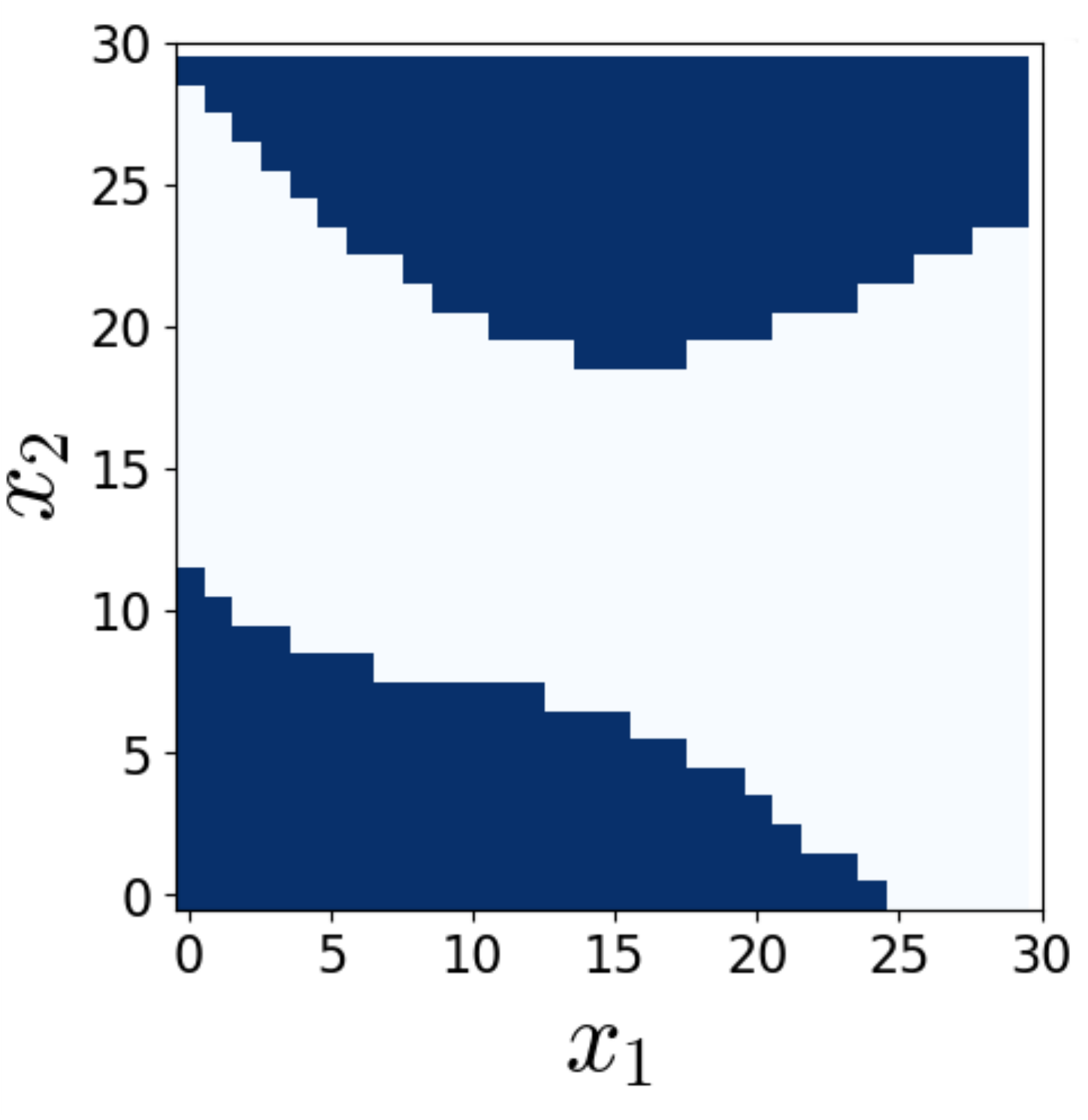}
         \caption{Logistic Regression}
     \end{subfigure}
    \begin{subfigure}[b]{0.41\textwidth}
         \centering
         \includegraphics[width=1.03\textwidth]{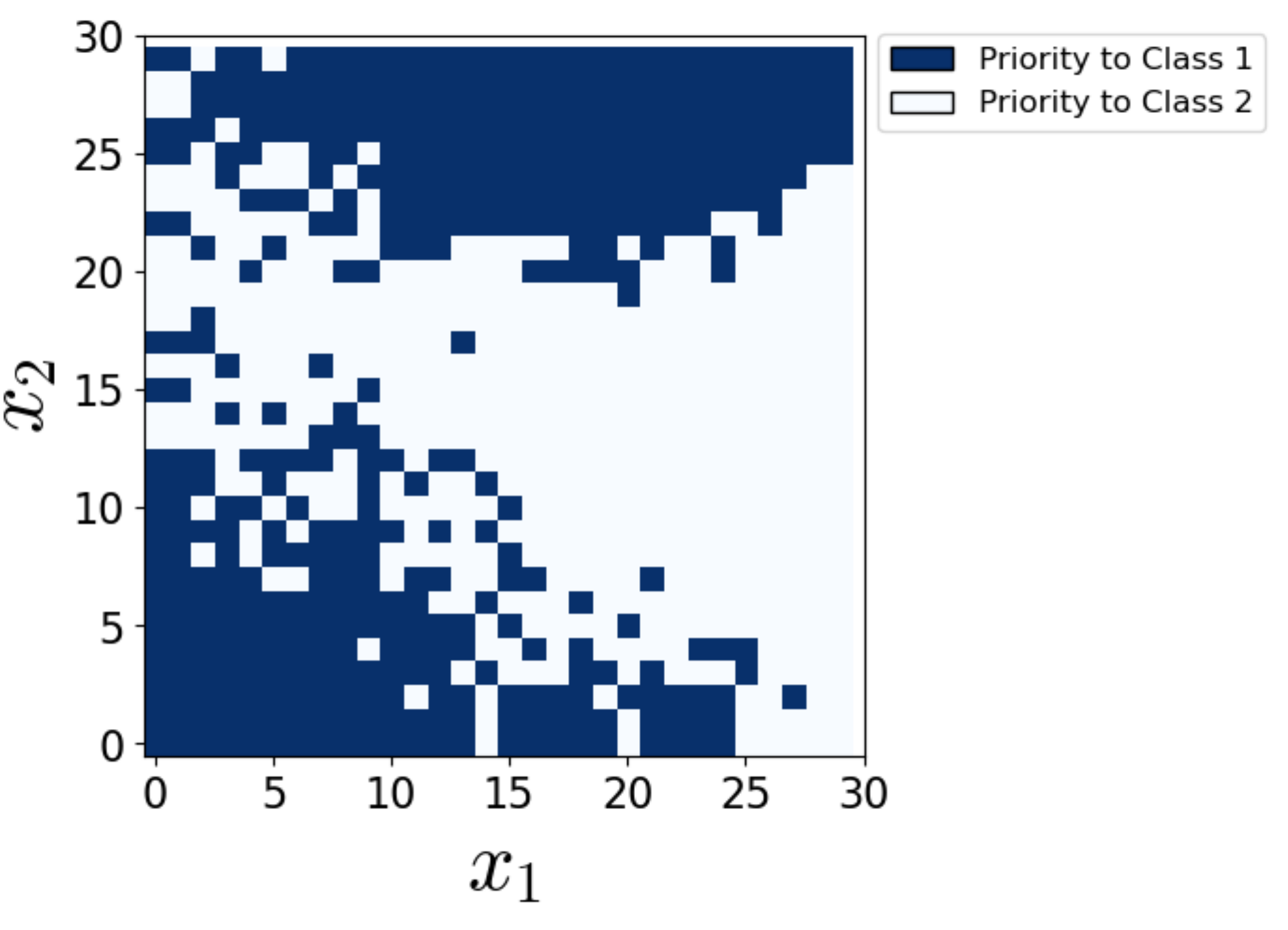}
         \caption{Direct Estimation}
     \end{subfigure}
     \caption{Optimal policy, ADP policy generated when Logistic Regression is used as a generalization method and ADP policy directly obtained for all states with $\lambda = (1.5,1.5)$, $\mu = (1,1)$, $a = (0.01,0.02)$, $\kappa = (30,30)$, $h = (3.5,1)$, $b = (0,0)$, $C=4$, $\bar{f}_i(x)= \mu_i - a_i x_i$.}
     \label{fig:general_example}
     \end{figure}

Note that based on the structure of the Bellman equation, the optimal policy is an index rule where the index for class $i$ is $\bar f_i(x_i)D_i^{\hat\pi_n}(x)$. Thus, we can characterize the policy by dividing the state space into different regions, each corresponding to a priority order. We utilize classifiers to learn this division. We observe in our numerical experiments that the classifiers achieve good generalization capability even though the simulation-based estimation of $D_i^{\hat\pi_n}(x)$ may suffer from a large variance.
Figure \ref{fig:general_example} provides a two-class example to demonstrate the good performance of the classifier. Figure \ref{fig:general_example} (a), (b) and (c) respectively depict the optimal policy, the one learned by a Logistic Regression, and the one learned using direct estimation at each state. In particular, for Figure \ref{fig:general_example} (c), we first estimate the value function differences at each state and then solve the corresponding linear program \eqref{eq:LP}. For Figure \ref{fig:general_example} (b), we use the estimated value function differences at only 25 different states to train a logistic regression, where the dependent variable is $\mathbbm{1}\{\bar f_1(x_1)\hat D_1(x)>f_2(x_2)\hat D_2(x)\}$. The features are third-order polynomials of $x$. The corresponding scheduling policy gives priority to class 1 if the predicted probability is greater than $0.5$, and gives priority to class 2 otherwise.
We observe that the classifier learns a policy that is very similar to the optimal policy, while the direct estimation suffers from the estimation noise.

\begin{remark}
In a two-class example, we can use a standard binary classifier.
When there are more than two customer classes, there are $I!$ possible priority rankings. In this case, we can either train an ML model to learn the ranking directly (e.g., leaning-to-rank method, e.g., \citealt{joachims2002optimizing}) or train a multiclass classifier (e.g., the multinomial logit model) and rank the classes based on the predicted probabilities. 
\end{remark}

Since $D_i^{\hat\pi_n}(x)$'s are only used to determine the order of $\bar f_i(x_i)D_i^{\hat\pi_n}(x)$'s, the estimation accuracy of $D_i^{\hat\pi_n}(x)$ can be adapted to the state $x$. In particular, when the $\bar f_i(x_i)D_i^{\hat\pi_n}(x)$-values are far from each other, a rather crude estimate of $D_i^{\hat\pi_n}(x)$'s are sufficient.
In Algorithm \ref{algo:adaptive_sampling}, we propose an adaptive sampling scheme to achieve a more efficient allocation of the sampling budget to different states. 



Algorithm \ref{algo:main} summarizes the main ADP algorithm. 

\begin{algorithm}
   Input an initial policy $\pi^{0}(\cdot)$ (We suggest using the fluid optimal policy as the initial policy, which can be obtained using Algorithm \ref{algo:initial}).
  Set $n = 0$ and input $N$ and $n_{max}$\;
    Uniformly sample $N$ points from the state-space $\mathcal{X}$ and denote the sampled set of points as $\mathcal{X}_0$ \; 
    \tcc{Policy evaluation:}
  Call Algorithm \ref{algo:adaptive_sampling} to estimate $D_i^{\pi_n}(x)$, $i=1,\dots, I$, for $x\in \mathcal{X}_0$ and denote the output as $\mathcal{T}_n=\{\hat D_i^{\pi_n}(x); i=1,\dots, I, x\in \mathcal{X}_0\}$\;
  \tcc{Policy improvement:}
  Using $\{\bar f_i(x_i)\hat D_i^{\pi_n}(x); i=1,\dots, I, x\in \mathcal{X}_0\}$ as the input data, train a classifier that classifies each state into a priority order based on the descending order of $\bar f_i(x_i)\hat D_i^{\pi_n}(x)$'s. Let $\pi^{n+1}$ denote the resulting policy\;
     Set $n = n+1$. If $n< n_{max}$, go to Step 3; otherwise, terminate and output $\pi^n$.
  \caption{Sampling-Based Approximate Policy Iteration}
  \label{algo:main}
\end{algorithm}

\subsection{Estimation of Value Function Differences via Coupling} \label{sec:coupling}
Estimating the relative value function for MDPs under the long-run average cost formulation is a challenging task since it involves an infinite time-horizon. 
\cite{cooper2003convergence} proposes a regenerative method to estimate the relative value functions. The idea involves selecting a regeneration state and estimating the relative value function at any given state based on the expected cost incurred until the system reaches that regeneration state. The regenerative estimation method can still be computationally challenging. First, the estimate tends to have a high variance. Second, the selection of the regeneration state affects both the simulation time and the variance of the estimate. This is particularly difficult for systems that exhibit meta-stability where the stationary distribution is multi-modal.
In this work, we focus on estimating the value function differences, which is sufficient to characterize the optimal policy. We propose a coupling strategy that does not require the selection of a regeneration state and provides estimates with lower variance.

Given a state $x$, we simulate $(I+1)$ coupled systems with initial point values $x$, $x-\mathbf{e}_1$, ..., $x-\mathbf{e}_I$ respectively. To make the dependence on the initial state explicit, we use $X(t;x)$ to denote a system starting from $x$. Define,
\[
\tau(x) = \inf\{k \geq 0; X(t_k;x) = X(t_k;x-\mathbf{e}_i)\ \forall i \in \{1,\dots,I\}\}.
\]
We construct the coupling in a way that $\tau(x)$ is small with a high probability, and for any $k\geq \tau(x)$, $X(t_k;x) = X(t_k;x-\mathbf{e}_i)$ for $i=1,2,\dots, I$. Under such coupling construction, we have,
\[
D_i(x)=\mathbb{E}\left[\sum_{k=0}^{\tau(x)}c(X(t_k;x))-c(X(t_k;x-\mathbf{e}_i))\right].
\]
Note that the coupling also creates positive correction between $c(X(t_k;x))$ and $c(X(t_k;x-\mathbf{e}_i))$, which leads to a variance reduction.

We next provide a construction of the coupling. We consider a discrete-event simulation construction with $I+1$ systems in parallel. At an even time $t$, two independent uniform random variables, $U_1$ and $U_2$, are generated. Define,
\[
\mu_{max} = \max \left\{
\sum\limits_{i=1}^{I} \bar{f}_i(X_i(t;x))Z_i(t;x),\sum\limits_{i=1}^{I} \bar{f}_i(X_i(t;x-\mathbf{e}_k))Z_i(t;x-\mathbf{e}_k),k=1,2,\dots,I
\right\},
\]
i.e., $\mu_{max}$ is the maximum possible departure rate across all systems.
Let $\Lambda$ = $\sum_{k=1}^{I}\lambda_k+\mu_{max}$.
We generate the event types in the $I+1$ systems as follows: 
\begin{enumerate}
\item If, 
\[
\sum_{k=1}^{i-1}\lambda_k/\Lambda<U_1\leq \sum_{k=1}^{i}\lambda_k/\Lambda,
\] 
for some $i=1,\dots, I$, there is a class $i$ arrival in all systems. 
\item If $U_1>\sum_{k=1}^{I}\lambda_k/\Lambda$, there is potentially a departure. We inspect each individual system.\\
For the system starting from $x$, if $U_1\leq \sum_{k=1}^{I}\lambda_k/\Lambda
+ \sum\limits_{i=1}^{I}\bar{f}_i(X_i(t;x))Z_i(t;x)$, there is a departure; otherwise, there is a phantom event. If there is a departure, we further check whether,
\[
\frac{ \sum_{k=1}^{i-1} \bar{f}_k(X_k(t;x))Z_k(t;x)}{\sum_{k=1}^{I} \bar{f}_k(X_k(t;x))Z_k(t;x)} < U_2\leq \frac{ \sum_{k=1}^{i} \bar{f}_k(X_k(t;x))Z_k(t;x)}{\sum_{k=1}^{I} \bar{f}_k(X_k(t;x))Z_k(t;x)}
\]
for some $i=1,\dots, I$, which indicates that there is a class $i$ departure from this system.\\
Similarly, for the system that starts from $x-\mathbf{e}_j$, $j=1,\dots, I$, there is a departure if $U\leq \sum_{k=1}^{I}\lambda_k/\Lambda + \sum\limits_{i=1}^{I}f(X_i(t;x-\mathbf{e}_j))Z_i(t;x-\mathbf{e}_j)$; otherwise, there is a phantom event. If there is a departure, we further check if,
\[
\frac{ \sum_{k=1}^{i-1} \bar{f}_k(X_k(t;x-\mathbf{e}_j))Z_k(t;x-\mathbf{e}_j)}{\sum_{k=1}^{I} \bar{f}_k(X_k(t;x-\mathbf{e}_j))Z_k(t;x-\mathbf{e}_j)} < U_2\leq \frac{ \sum_{k=1}^{i} \bar{f}_k(X_k(t;x-\mathbf{e}_j))Z_k(t;x-\mathbf{e}_j)}{\sum_{k=1}^{I} \bar{f}_k(X_k(t;x-\mathbf{e}_j))Z_k(t;x-\mathbf{e}_j)}
\]
for some $i=1,\dots, I$, which indicates that there is a class $i$ departure from this system.
\end{enumerate}

Note using the same $U_1$ and $U_2$ for the $I+1$ systems leads to a coupling. It is straightforward to see that once two systems coincide under this coupling, they will have identical sample paths afterwards.
We also remark that the coupling construction is not unique, i.e., other constructions can also lead to similar performance. The algorithm to estimate the value function differences under the coupling construction is summarized in Algorithm \ref{algo:coupling}. 

\begin{algorithm}
  \caption{Value Function Difference Estimation.}
  \label{algo:coupling}
  Input $T$ and a policy $\pi$;\\
  \For{$l\in\{1,\dots, n\}$}{
  Simulate $I+1$ systems, $X(t_k;x)$ and $X(t_k;x-\mathbf{e}_i)$ $\forall i \in \{1,\dots,I\}$, under policy $\pi$ and the coupling construction described above until $\min\{\tau(x),T\}$\; 
  Calculate $\mathcal{D}_{i,l}=\sum_{k=0}^{\min\{\tau(x),T\}}c(X(t_k;x))-c(X(t_k;x-\mathbf{e}_i))$ $\forall i \in \{1,\dots,I\}$\; 
  }
  Output $\mathcal{D}_{i,1}, \dots, \mathcal{D}_{I,n}$.
\end{algorithm}

\begin{figure}
    \centering
    \includegraphics[scale=0.35]{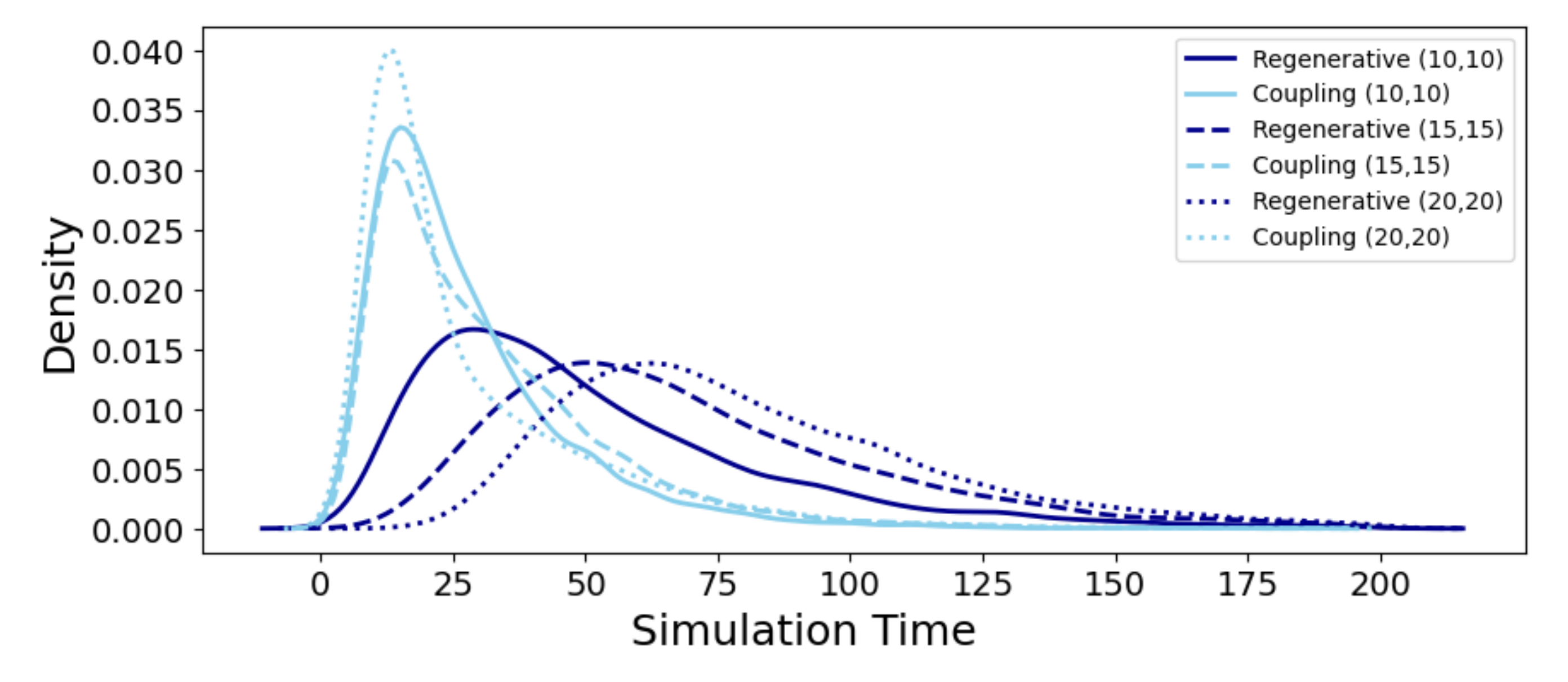}
    \caption{Distribution of the simulation times using the coupling method versus the regeneration method. (Arrival rates: $\lambda = (1.5,1.5)$; base service rates $\mu = (1,1)$; slowdown rates: $a = (0.0103, 0.0203)$; blocking thresholds: $\kappa = (30,30)$; blocking costs: $b=(0,0)$; holding costs: $h = (3,1)$. Initial states: $(10,10), (15,15), (20,20)$.)}
    \label{fig:sim_time_dist}
\end{figure}

We next illustrate the advantages of the coupling method compared to the standard regenerative method through some numerical experiments. 
We consider a system with two customer classes and linear slowdown, i.e., $\bar f(x_i)=\mu_i-a_ix_i$. The detailed system parameters are listed in the caption of Figure \ref{fig:sim_time_dist}. 
Figure \ref{fig:sim_time_dist} plots the distribution of the {\em simulation time} using the coupling method versus the regenerative method starting from $x=(10,10),(15,15),(20,20)$. The simulation time for the coupling method is the time for the $I+1$ systems to coincide, i.e., $\tau(x)$. The simulation time for the regenerative method is the maximum regeneration time for the $I+1$ systems. In particular, for a fixed regeneration state $x_0$, let 
$\varsigma(x)=\inf\{k\geq 0: X(t_k;x)=x_0\}$,
then, the simulation time is $\max\{\varsigma(x),\varsigma(x-\mathbf{e}_i), i=1,\dots, I\}$. Note that since our goal is to estimate $D_i(x), i\in\{1,\ldots, I\}$, we simulate $I+1$ systems for each $x$. We use common random numbers across the $I+1$ systems to increase the correlation between the estimates of $v(x)$ and $v(x-\mathbf{e}_i)$.
We set $(1,1)$ as the regeneration state because it is the closest state to $(\lambda_1/\mu_1,\lambda_2/\mu_2)$. We make two observations from the figure. First, the coupling times are in general shorter than the regeneration times.
This is because the coupling method does not require hitting a single regeneration state. The system can be fully synchronized in multiple ways.
Second, as the starting state gets further away from the regenerative state, the regeneration time tends to be longer. However, the coupling method does not suffer from this since the systems can also coincide at the truncation/blocking boundary.

Table \ref{tab:estimation_acc} compares the standard deviation of the estimate of $\hat D(x)$'s using the coupling method versus the regenerative method. The standard deviation is estimated using $10^3$ replications. We observe that the coupling method achieves a much smaller standard deviation. In particular, the coupling method provides a variance reduction in the range of $78\%$ to $93\%$ compared to the regenerative method. 
This is 
likely because the coupling method creates a stronger correlation between $X(t_k;x)$ and $X(t_k;x-\mathbf{e}_i)$.  
Smaller estimation variance together with smaller simulation time translates to more efficient estimation. Table \ref{tab:estimation_acc} also shows the number of replications needed and the corresponding total running time to achieve a 99\% probability of correct priority ranking. We observe the coupling method requires fewer replications and shorter running times than the regeneration method.

\begin{table}[!h]
\centering
\caption{Standard deviation (std) of the value function difference estimates, the number of replications needed and the corresponding running time required to achieve 99\% probability of correct priority ranking using coupling method versus regeneration method. (The std is estimated based on $10^3$ replications.)} \label{tab:estimation_acc}
\resizebox{\textwidth}{!}{\begin{tabular}{c|cc|cc|cc}
\hline
\textbf{}        & \multicolumn{2}{c|}{\textbf{Std of $\hat{D}_1(x_1,x_2),\hat{D}_2(x_1,x_2)$}} & \multicolumn{2}{c|}{\textbf{Average Running Time (sec)}} & \multicolumn{2}{c}{\textbf{Average \# of Replications}} \\ \hline
$(x_1,x_2)$      & \textbf{Regenerative}                   & \textbf{Coupling}                  & \textbf{Regenerative}         & \textbf{Coupling}        & \textbf{Regenerative}        & \textbf{Coupling}        \\ \hline
\textbf{(10,10)} & 261, 283                                & 67, 93                             & 134.19                        & 12.66                    & 2426.60                      & 582.00                   \\
\textbf{(10,15)} & 225, 300                                & 62, 112                            & 58.08                         & 6.86                     & 860.90                       & 222.95                   \\
\textbf{(10,20)} & 216, 366                                & 70, 172                            & 28.63                         & 4.83                     & 355.60                       & 103.15                   \\
\textbf{(15,10)} & 244, 273                                & 70, 87                             & 162.96                        & 10.45                    & 2557.00                      & 428.40                   \\
\textbf{(15,15)} & 199, 276                                & 65, 109                            & 52.56                         & 5.18                     & 724.30                       & 131.40                   \\
\textbf{(15,20)} & 131, 329                                & 56, 153                            & 18.65                         & 3.83                     & 211.20                       & 60.35                    \\
\textbf{(20,10)} & 211, 243                                & 65, 78                             & 109.97                        & 8.79                     & 1559.70                      & 353.00                   \\
\textbf{(20,15)} & 153, 247                                & 52, 93                             & 32.81                         & 3.76                     & 413.50                       & 67.30                    \\
\textbf{(20,20)} & 81, 284                                 & 35, 140                            & 12.34                         & 3.39                     & 124.95                       & 38.70                    \\ \hline
\end{tabular}}

\label{tab:regen_std}
\end{table}

\subsection{Adaptive Sampling Approach}\label{sec:adaptive_sampling}
To deal with the curse of dimensionality, we only sample from a subset of states to estimate the value function differences. We further develop an adaptive sampling scheme to dynamically allocate the sampling budget depending on how difficult it is to determine the ranking among $\Bar{f}_i(x_i)D_i^{\pi}(x)$'s. 
The idea is that if $|\Bar{f}_i(x_i)D_i^{\pi}(x)-\Bar{f}_j(x_j)D_j^{\pi}(x)|$ is large for $i\neq j$, it is easier to determine the order. Thus, a smaller number of samples would be sufficient to estimate $D_i^{\pi}(x), i\in\{1,\ldots,I\}$. Let $R_i(x)=\Bar{f}_i(x_i)D_i^{\pi}(x)$. 
We denote $\hat R_i(x)$ as the estimate for $R_i(x)$ based on sample averages, i.e.,
\[
\hat R_i(x)=\Bar{f}_i(x_i)\frac{1}{n}\sum_{k=1}^{n} \mathcal{D}_{i,k}(x),
\]
where $\mathcal{D}_{i,k}(x),i\in\{1,\ldots,I\}$ are iid samples of $\sum_{k=0}^{\tau(x)}c(X(t_k;x))-c(X(t_k;x-\mathbf{e}_i))$.
We also denote $\hat V_{ij}(x)$ as the estimate for $Var(\hat R_i(x)-\hat R_j(x))$ based on sample variance, i.e.,
\[
\hat V_{ij}(x)=\frac{1}{(n-1)^2}\sum_{k=1}^{n}[(\Bar{f}_i(x_i)\mathcal{D}_{i,k}(x)-\Bar{f}_j(x_j)\mathcal{D}_{j,k}(x)) - (\hat R_i(x)-\hat R_j(x))]^2.
\]
For a target confidence level $\alpha$, let $\phi_{\alpha}$ denote the $(1-\alpha)/2$-th upper percentile of the standard Normal distribution.
We divide the simulation into several iterations. In each iteration, we continue generating more samples associated with $x$ if,
\[
\min_{1\leq i,j\leq I, i\neq j}\frac{|\hat R_i(x)-\hat R_j(x)|}{\sqrt{\hat V_{ij}(x)}}\leq \phi_{\alpha}.
\]
The algorithm is summarized in Algorithm \ref{algo:adaptive_sampling}.

\begin{algorithm}
  Input $n_{step}$, $\tilde{n}_{max}$, $\alpha$,
  a policy $\pi$, and a subset of states $\mathcal{X}_0$. Set $\mathcal{M}^{0}=\mathcal{X}_0$, $k=0$, and $n=0$\;
  For each point $x \in \mathcal{M}^{k}$, call Algorithm \ref{algo:coupling} to generate $n_{step}$ samples of $\mathcal{D}_i^{\pi}(x)$, $i=1,\dots, I$. \;
  Set $n=n+n_{step}$. Use the $n$ samples generated so far to calculate $\hat{R}_i(x)$ and $\hat V_{ij}(x)$ for $x\in \mathcal{M}^{k}$, $i,j=1\dots I$, $i\neq j$.\;
  For each $x\in \mathcal{M}^{k}$, if
  \[
    \min_{1\leq i,j\leq I, i\neq j}\frac{|\hat R_i(x)-\hat R_j(x)|}{\sqrt{\hat V_{ij}(x)}}\leq \phi_{\alpha},
 \]
 add $x$ to $\mathcal{M}^{k+1}$. Set $k=k+1$\;
 If $\mathcal{M}^{k}$ empty or $n\geq \tilde n_{max}$, terminate and output $\hat{D}_i(x)$ as the sample average of the generated $\mathcal{D}_i^{\pi}(x)$'s for all $x\in\mathcal{X}_0$, $i=1,\dots, I$; otherwise, go to Step 2.
  \caption{Adaptive Sampling} \label{algo:adaptive_sampling}
\end{algorithm}

We next demonstrate the benefits of applying the adaptive sampling approach through some numerical examples. Consider a two-class system linear service slowdown. We apply the proposed ADP algorithm with and without adaptive sampling respectively. We set $\tilde n_{max}=2000$. In the ADP algorithm without adaptive sampling, each $D_i(x)$ is estimated based on the average of $2000$ samples. For the ADP algorithm with adaptive sampling, we set $\alpha=0.95$. In this case, the number of samples generated can be much smaller than $2000$ for some $x$. 

We observe that the policies are almost identical, and the policies learned with and without the adaptive sampling lead to a long-run average cost of $14.72 \pm 0.20$ and $14.71 \pm 0.21$ respectively. Meanwhile, the number of samples generated using the adaptive sampling algorithm is only $47\%$ of that without adaptive sampling. In other words, the adaptive sampling achieves approximately the same performance as the ADP algorithm without adaptive sampling but reduces the sampling cost by $53\%$.

\begin{figure}[]
 \centering
 \begin{subfigure}[b]{0.465\textwidth}
     \centering
     \includegraphics[width=\textwidth]{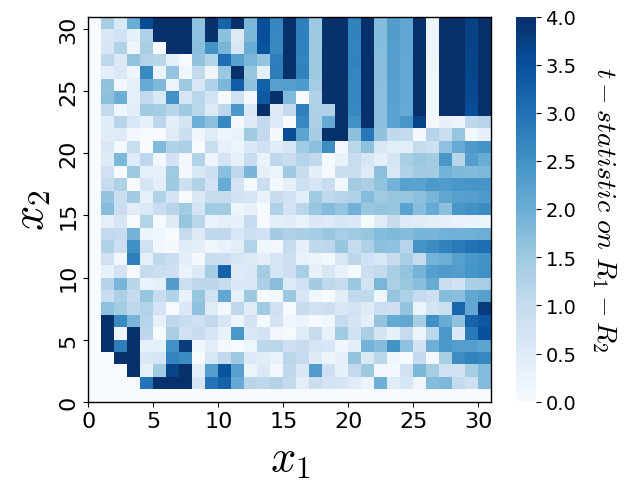}
     \caption{t-statistics on $R_1 - R_2$}
     \label{fig:adap_tstat}
 \end{subfigure}
 \begin{subfigure}[b]{0.48\textwidth}
     \centering
     \includegraphics[width=\textwidth]{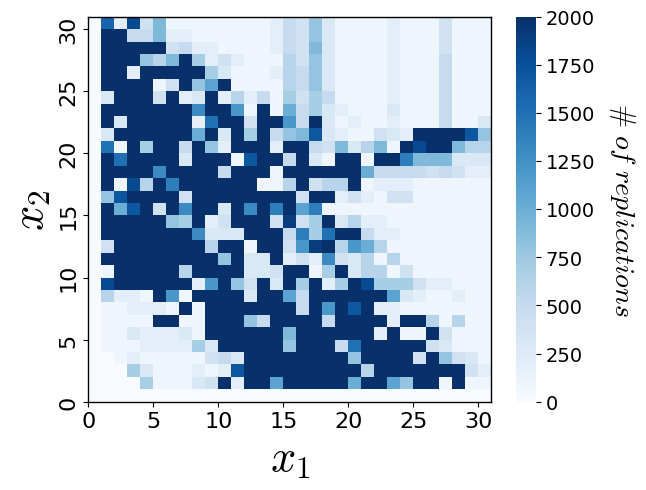}
     \caption{Number of replications}
     \label{fig:adap_nofRef}
 \end{subfigure}
 \caption{t-statistics on  $R_1-R_2$ calculated over 100 replications and the number of replications taken in different states of the sample space for identifying the correct action with a high probability guarantee ($\alpha = 0.95$).}
\label{fig:adaptive_sampling_nofRep}
 \end{figure}

Figure \ref{fig:adaptive_sampling_nofRep}  illustrates how the sample budget is allocated among different states. Figure \ref{fig:adaptive_sampling_nofRep} (a) shows the value of $t$-statistic defined as $|\hat R_1(x)-\hat R_2(x)|/\sqrt{\hat V_{12}(x)}$, 
whereas Figure \ref{fig:adaptive_sampling_nofRep} (b) shows the number of samples generated for different $x$. We observe that when it is hard to tell $R_1(x)$ and $R_2(x)$ apart (in the middle region when the t-statistic is small), the algorithm generates more samples. On the other hand, when it is easier to decide the order between $R_1(x)$ and $R_2(x)$ (when $x$ is close to zero or $[30,30]$), much fewer samples are generated. 

More details regarding the implementation of the algorithm including guidelines on setting its inputs (e.g., $N$ in Algorithm \ref{algo:main}, $T$ in Algorithm \ref{algo:coupling}) are provided in Section \ref{sec:implementation_details}. 

\section{Extension to the Non-preemptive Setting}\label{sec:nonpreemptive}
In this section, we extend our ADP method to the non-preemptive setting. 
To keep the system state Markovian, we need to keep track of the number of customers in service $Z$ in addition to the number of customers in the system $X$. With a slight abuse of notation, we continue to denote $\mathcal{X}$ as the set of all states. We define $\mathcal{X}(A)\subset \mathcal{X}$ as the set of actionable states, i.e., the states in which an action needs to be taken. We focus on \emph{non-idling} policies due to their practicality, although this assumption can be relaxed at the expense of a larger set of actionable states and action sets. The actionable states are those with an available server and more than one customer waiting to be served. Because ongoing services cannot be interrupted, an action determines which class of customers in the queue is assigned the available server. 



 Let $\mathcal{I}(x,z):= \{i\in \{1,\ldots,I\}: x_i>z_i\}$ denote the set of admissible actions in a given state $(x,z)$. For the long-run average cost problem, we solve the following Bellman equation,
\begin{align} \label{eq:bellman_nonpreemptive}
     v^*(x,z) &= \frac{1}{\Lambda}\left\{\sum\limits_{i=1}^{I} h_i x_i + \lambda_i b_i \mathbbm{1}_{\{x_i = \kappa_i\}} - \gamma^* + \sum\limits_{i=1}^{I} \lambda_i\mathbbm{1}_{\{x_i < \kappa_i, \sum_{i=1}^{I}z_i =C\}} v^{*}(x + \mathbf{e}_i,z) \right. \nonumber \\
     &+\sum\limits_{i=1}^{I} \lambda_i\mathbbm{1}_{\{x_i < \kappa_i,\sum_{i=1}^{I}z_i <C\}} v^{*}(x + \mathbf{e}_i,z+\mathbf{e}_i)+ \sum\limits_{i=1}^{I} z_i\Bar{f}_i(x_i) v^{*}(x-\mathbf{e}_i,z-\mathbf{e}_i) 
    \nonumber \\ 
    &\left. + \Bigl(\Lambda - \sum\limits_{i=1}^{I} \lambda_i\mathbbm{1}_{\{x_i < \kappa_i\}} - \Bar{f}_i(x_i) z_i \Bigr) v^*(x,z) \right\}, \quad \forall x,z\in \mathcal{X}-\mathcal{X}(A) \\
    v^*(x,z) &= \min\limits_{i \in \mathcal{I}(x,z)} \left\{
    v^{*}(x,z+\mathbf{e}_i)
    \right\}, \forall x,z\in \mathcal{X}(A), \label{eq:bellman_nonpreemptive_actionable}
\end{align}
where $\gamma^*$ is the optimal long-run average cost and $v^*(\cdot)$ is the optimal relative value function. 

Let $i' = \min\{i:i\in \mathcal{I}(x,z)\}$ be the smallest indexed action in $\mathcal{I}$. We define $\tilde D^{*}_i(x,z) := v^{*}(x,z+\mathbf{e}_{i'}-v^{*}(x,z+\mathbf{e}_i)), \forall i \in \mathcal{I}(x,z)$. The optimal action is then to serve the class with the largest $\tilde D_{i}^{*}(x,z)$.
Note that similar to the preemptive system, the optimal action only depends on the value function differences.
For the non-preemptive system, the value function differences are defined slightly differently and thus require a different estimation algorithm.
Let $S$ denote the enlarged state space that includes both $X$ and $Z$. We define,
\[
\tilde \tau(x) = \inf\{k \geq 0; S(t_k;x,z+\mathbf{e}_{i'}) = S(t_k;x,z+\mathbf{e}_i),\ \forall i \in \mathcal{I}\}.
\]
The updated value function difference estimation algorithm for the non-preemptive system is provided in Algorithm \ref{algo:coupling_nonpreemptive}. 
\begin{algorithm}
  \caption{Value Function Difference Estimation for Nonpreemptive Systems.}
  \label{algo:coupling_nonpreemptive}
  Input $T$ and a policy $\pi$;\\
  \For{$l\in\{1,\dots, n\}$}{
  Simulate $|\mathcal{I}|$ systems, $S(t_k;x,z+\mathbf{e}_{i'})$ and $S(t_k;x,z+\mathbf{e}_i)$ $\forall i \in \mathcal{I}$, under policy $\pi$ and the coupling construction described in Section \ref{sec:coupling} until $\min\{\tilde \tau(x),T\}$\; 
  Calculate 
  $\mathcal{D}_{i,l}=\sum_{k=0}^{\min\{\tilde \tau(x),T\}}c(S(t_k;x,z+\mathbf{e}_{i'}))-c(S(t_k;x,z+\mathbf{e}_i))$ $\forall i \in \mathcal{I}$\; 
  }
  Output $\mathcal{D}_{i,1}, \dots, \mathcal{D}_{|\mathcal{I}|,n}$.
\end{algorithm}

\section{ADP Performance and Structure of the Optimal Policy}\label{sec:num_experiments}

In this section, we first show that the proposed ADP learns the near-optimal policy. Next, we provide some insights into the structure of the optimal policy.
To be able to calculate the optimal policy, i.e., solve the MDP exactly, we focus on a two-class system with $C=4$ and $\kappa=(30,30)$. We demonstrate the performance of our algorithm on larger-scale problems in Section \ref{sec:fiveclass}. We consider a linearly decreasing service rate function in our experiments, i.e.,
$
\bar{f}_i (x_i) =\mu_i - a_i x_i,
$
where $a_i$ is the slowdown rate and $\mu_i$ is the maximum service rate for class $i$ customers.

\subsection{The performance of the ADP Algorithm}\label{sec:compare_benchmark}

In this subsection, we investigate the performance of our ADP algorithm. We refer to the scheduling policy learned by our ADP algorithm as the ADP policy. We also investigate the performance of benchmark scheduling policies that have been proposed in the literature for parallel server queues without service slowdown. We adapt these benchmarks to heuristically capture service slowdown. 
We consider five benchmark policies: 
\begin{enumerate}
    \item The $c\mu$ policy ($h\bar f(0)$): it gives priority to the class with a higher $h_i \bar{f}_i(0)$ value. 
    \item The adapted state-dependent $c\mu$ policy ($h\bar f(x)$): it gives priority to the class with a higher $h_i \bar{f}_i(X_i(t))$ value, where $X_i(t)$ is the state of the queue at the decision epoch. 
    \item The adapted state-dependent maximum pressure (Max Pressure) policy: it gives priority to the class with a larger $h_i X_i(t) \bar{f}_i(X_i(t))$ value. 
    \item The shortest-queue-first policy (SQF): it gives priority to the class with fewer customers in the system
    \item The longest-queue-first policy (LQF): it gives priority to the class with more customers in the system. 
\end{enumerate}
We conduct extensive sensitivity analysis on the performance of these policies by considering different worst-case system loads, i.e., $\lambda_1/(C\bar f_1(\kappa_1))+\lambda_2/(C\bar f_2(\kappa_2))$; different values of $\mu_1-\mu_2$; different values of $\lambda_1-\lambda_2$; and different blocking costs, $b=(b_1,b_2)$. 
We also vary the holding cost of class 1, $h_1$ (the holding cost of class 2 is fixed at $1$), and the service slowdown rates, $a=(a_1,a_2)$. Our main observation is that the ADP policy achieves robust and near-optimal performance. While the benchmarks can perform well in some cases, they could have very large optimality gaps in other cases. This highlights the importance of carefully designing the scheduling policy when facing service slowdown. We present the main insights and some supporting numerical evidence here. The extensive numerical results are provided in Section \ref{sec:experiment_blocking}).


\subsubsection{Worst-case system load.}\label{sec:worst-case_load} We examine the effect of the worst-case system load, $\lambda_1 / (\mu_1 - a_1 \kappa_1) + \lambda_2 / (\mu_2 - a_2 \kappa_2)$, on the performance of the policies. We fix the service rates $\mu = (1,1)$, arrival rates $\lambda= (1.5,1.5)$, and blocking cost $b=(0,0)$, and vary the worse-case system load by changing the slowdown rates. Table \ref{tab:slowdown_effect_experiments_summary} summarizes the optimal long-run average costs and the optimality gaps for different policies (the reported numbers are averages over systems with different holding costs). We observe the following. (1) The ADP policy yields a very low optimality gap, on average 0.5\%. (2) The optimality gap of the ADP policy tends to increase with the worst-case system load, from $0.21\%$ to $1.41\%$. This is likely because as the system load increases, the estimates of relative value function differences suffer from a larger variance, which leads to more estimation errors. 
(3) As expected, the average cost under the optimal policy increases with the worst-case system load. (4) Policies $h\bar f(x)$, Max Pressure, SQF, and LQF in general do not perform well. (5) The policy $h\bar{f}(0)$ performs well when the worst-case system load is relatively low, i.e., $1$ or $1.25$. However, its performance can degrade significantly and become highly suboptimal when the worst-case system load is high. The average optimality gap is $24.05\%$ when the worst-case system load is $1.5$.


\begin{table}[]
\centering
\caption{Average optimal cost and average \% optimality gaps over different holding costs under ADP policy and five benchmark policies and varying service slowdown rates.}
\label{tab:slowdown_effect_experiments_summary}
\begin{tabular}{l|c|cccccc}
\hline
\textbf{\begin{tabular}[c]{@{}l@{}}Worst-case\\ Load\end{tabular}} & \textbf{Optimal}                                           & \textbf{ADP}                                              & \textbf{$h \bar f(0)$}                                        & \textbf{$h \bar f(x)$}                                        & \textbf{\begin{tabular}[c]{@{}c@{}}Max\\ Pressure\end{tabular}} & \textbf{SQF}                                               & \textbf{LQF}                                                \\ \hline
\textbf{1}                                                         & \begin{tabular}[c]{@{}c@{}}8.14\\ $\pm$ 0.04\end{tabular}  & \begin{tabular}[c]{@{}c@{}}0.21\\ $\pm$ 0.03\end{tabular} & \begin{tabular}[c]{@{}c@{}}0.13\\ $\pm$ 0.03\end{tabular}  & \begin{tabular}[c]{@{}c@{}}1.91\\ $\pm$ 0.11\end{tabular}  & \begin{tabular}[c]{@{}c@{}}6.17\\ $\pm$ 0.09\end{tabular}       & \begin{tabular}[c]{@{}c@{}}13.00\\ $\pm$ 0.20\end{tabular}  & \begin{tabular}[c]{@{}c@{}}15.13\\ $\pm$ 0.24\end{tabular}  \\ \hline
\textbf{1.25}                                                      & \begin{tabular}[c]{@{}c@{}}9.73\\ $\pm$ 0.06\end{tabular}  & \begin{tabular}[c]{@{}c@{}}0.29\\ $\pm$ 0.10\end{tabular} & \begin{tabular}[c]{@{}c@{}}3.59\\ $\pm$ 0.19\end{tabular}  & \begin{tabular}[c]{@{}c@{}}9.37\\ $\pm$ 0.36\end{tabular}  & \begin{tabular}[c]{@{}c@{}}18.78\\ $\pm$ 0.42\end{tabular}      & \begin{tabular}[c]{@{}c@{}}12.29\\ $\pm$ 0.19\end{tabular} & \begin{tabular}[c]{@{}c@{}}70.08\\ $\pm$ 1.57\end{tabular}  \\ \hline
\textbf{1.5}                                                       & \begin{tabular}[c]{@{}c@{}}11.68\\ $\pm$ 0.09\end{tabular} & \begin{tabular}[c]{@{}c@{}}1.41\\ $\pm$ 0.13\end{tabular} & \begin{tabular}[c]{@{}c@{}}24.05\\ $\pm$ 0.59\end{tabular} & \begin{tabular}[c]{@{}c@{}}39.79\\ $\pm$ 0.86\end{tabular} & \begin{tabular}[c]{@{}c@{}}64.23\\ $\pm$ 1.24\end{tabular}      & \begin{tabular}[c]{@{}c@{}}24.17\\ $\pm$ 0.49\end{tabular} & \begin{tabular}[c]{@{}c@{}}479.55\\ $\pm$ 4.29\end{tabular} \\ \hline
\end{tabular}

\raggedright
\tabnote{Other system parameters are fixed at $\lambda = (1.5,1.5)$, $\mu=(1,1)$, $h \in \{(1,1),(1.5,1),(2,1),(2.5,1),(3,1),(3.5,1)\}$, $b=(0,0)$, $C=4$, $\kappa=(30,30)$. Detailed long-run average costs for each holding costs under optimal policy, ADP policy and five benchmark policies can be found in Table \ref{tab:slowdown_effect_experiments} in the Online Supplementary.}
\end{table}

\subsubsection{Service rate heterogeneity.}\label{sec:service_rate_heterogeneity} We next examine the effect of service rate heterogeneity on the performance of different policies. We fix the arrival rates $\lambda= (1.5,1.5)$ and blocking cost $b=(0,0)$, and vary the difference between $\mu_1$ and $\mu_2$ while keeping the worse-case load at $1.5$ by changing the slowdown rates $a=(a_1,a_2)$. We keep the difference between $a_1$ and $a_2$ at $0.01$. Table \ref{tab:service_effect_experiments_summary} summarizes the long-run average cost and the average optimality gaps for different policies under different system parameters. We observe that (1) the ADP policy yields a very small optimality gap in all cases, less than $1.08\%$.
(2) Even though the worst-case system load is kept the same, as the difference between the service rates increases, the optimal cost tends to slightly increase. 
(3) Benchmark policies perform quite poorly in most cases.

\begin{table}[]
\centering
\caption{Average optimal cost and average \% optimality gaps over different holding costs under ADP policy and five benchmark policies and varying service and slowdown rates.}
\label{tab:service_effect_experiments_summary}
\begin{tabular}{cc|c|cccccc}
\hline
\textbf{$\mu$}          & \textbf{$a$}              & \textbf{Optimal}                                          & \textbf{ADP}                                              & \textbf{$h \bar{f}(0)$}                                    & \textbf{$h \bar{f}(x)$}                                    & \textbf{\begin{tabular}[c]{@{}c@{}}Max\\ Pressure\end{tabular}} & \textbf{SQF}                                               & \textbf{LQF}                                                \\ \hline
\textbf{(0.975, 1.025)} & \textbf{(0.0107, 0.0207)} & \begin{tabular}[c]{@{}c@{}}8.57\\ $\pm$ 0.04\end{tabular} & \begin{tabular}[c]{@{}c@{}}0.96\\ $\pm$ 0.09\end{tabular} & \begin{tabular}[c]{@{}c@{}}35.00\\ $\pm$ 0.64\end{tabular} & \begin{tabular}[c]{@{}c@{}}44.05\\ $\pm$ 0.56\end{tabular} & \begin{tabular}[c]{@{}c@{}}83.62\\ $\pm$ 1.27\end{tabular}      & \begin{tabular}[c]{@{}c@{}}24.90\\ $\pm$ 0.42\end{tabular} & \begin{tabular}[c]{@{}c@{}}496.53\\ $\pm$ 4.28\end{tabular} \\
\textbf{(0.952, 1.052)} & \textbf{(0.0111, 0.0211)} & \begin{tabular}[c]{@{}c@{}}8.73\\ $\pm$ 0.05\end{tabular} & \begin{tabular}[c]{@{}c@{}}0.71\\ $\pm$ 0.10\end{tabular} & \begin{tabular}[c]{@{}c@{}}28.22\\ $\pm$ 0.56\end{tabular} & \begin{tabular}[c]{@{}c@{}}30.00\\ $\pm$ 0.52\end{tabular} & \begin{tabular}[c]{@{}c@{}}82.34\\ $\pm$ 1.38\end{tabular}      & \begin{tabular}[c]{@{}c@{}}19.19\\ $\pm$ 0.37\end{tabular} & \begin{tabular}[c]{@{}c@{}}504.25\\ $\pm$ 4.19\end{tabular} \\
\textbf{(0.930, 1.080)} & \textbf{(0.0115, 0.0215)} & \begin{tabular}[c]{@{}c@{}}8.99\\ $\pm$ 0.05\end{tabular} & \begin{tabular}[c]{@{}c@{}}0.71\\ $\pm$ 0.11\end{tabular} & \begin{tabular}[c]{@{}c@{}}22.45\\ $\pm$ 0.50\end{tabular} & \begin{tabular}[c]{@{}c@{}}22.62\\ $\pm$ 0.49\end{tabular} & \begin{tabular}[c]{@{}c@{}}85.03\\ $\pm$ 1.71\end{tabular}      & \begin{tabular}[c]{@{}c@{}}14.05\\ $\pm$ 0.31\end{tabular} & \begin{tabular}[c]{@{}c@{}}506.20\\ $\pm$ 4.04\end{tabular} \\
\textbf{(0.910, 1.110)} & \textbf{(0.0119, 0.0219)} & \begin{tabular}[c]{@{}c@{}}9.23\\ $\pm$ 0.05\end{tabular} & \begin{tabular}[c]{@{}c@{}}1.08\\ $\pm$ 0.15\end{tabular} & \begin{tabular}[c]{@{}c@{}}16.46\\ $\pm$ 0.42\end{tabular} & \begin{tabular}[c]{@{}c@{}}16.07\\ $\pm$ 0.41\end{tabular} & \begin{tabular}[c]{@{}c@{}}91.83\\ $\pm$ 2.19\end{tabular}      & \begin{tabular}[c]{@{}c@{}}9.97\\ $\pm$ 0.27\end{tabular}  & \begin{tabular}[c]{@{}c@{}}505.54\\ $\pm$ 3.77\end{tabular} \\ \hline
\end{tabular}

\raggedright
\tabnote{Other system parameters are fixed at $\lambda = (1.5,1.5)$, $h \in \{(1,1),(1.25,1),(1.5,1),(1.75,1),(2,1),(2.25,1)\}$, $b=(0,0)$, $C=4$, $\kappa=(30,30)$. Detailed long-run average costs for each holding costs under optimal policy, ADP policy and five benchmark policies can be found in Table \ref{tab:service_effect_experiments} on the Online Supplementary.}
\end{table}
\subsubsection{Arrival rate heterogeneity.} We next examine the effect of arrival rate heterogeneity on the performance of different policies. We fix $\mu = (1,1)$, $a=(0.0103,0.0203)$, $b=(0,0)$, and $h=(1.5,1)$, and vary $\lambda$ while making sure that the best- and worst-case system loads are 0.75 and 1.5 respectively. 
Table \ref{tab:arrival_effect_experiments} summarizes the long-run average costs and the average optimality gaps for different policies under different arrival rates. We observe the following. (1) The ADP policy yields a very small optimality gap, $0.72\%$ on average. 
(2) The benchmark policies, which are agnostic to the arrival rates, can perform quite poorly, especially when $\lambda_1$ is very different from $\lambda_2$. When $\lambda_1$ and $\lambda_2$ are close to each other, $h\bar{f}(0)$ and $h\bar{f}(x)$ perform quite well.

\begin{table}[]
\caption{Average optimal cost and average \% optimality gaps under ADP policy and five benchmark policies and varying arrival rates.}
\centering
\label{tab:arrival_effect_experiments}
\scalebox{1}{
\begin{tabular}{l|c|cccccc}
\hline
$\lambda$                   & \textbf{Optimal}                                          & \textbf{ADP}                                               & \textbf{$h \bar{f}(0)$}                                    & \textbf{$h \bar{f}(x)$}                                    & \textbf{\begin{tabular}[c]{@{}c@{}}Max\\ Pressure\end{tabular}} & \textbf{SQF}                                               & \textbf{LQF}                                                \\ \hline
\textbf{0.5,2.5}            & \begin{tabular}[c]{@{}c@{}}9.88\\ $\pm$ 0.06\end{tabular} & \begin{tabular}[c]{@{}c@{}}0.09\\ $\pm$ 0.26\end{tabular}  & \begin{tabular}[c]{@{}c@{}}65.98\\ $\pm$ 1.04\end{tabular} & \begin{tabular}[c]{@{}c@{}}65.98\\ $\pm$ 1.04\end{tabular} & \begin{tabular}[c]{@{}c@{}}154.65\\ $\pm$ 2.30\end{tabular}     & \begin{tabular}[c]{@{}c@{}}65.81\\ $\pm$ 1.04\end{tabular} & \begin{tabular}[c]{@{}c@{}}508.26\\ $\pm$ 5.42\end{tabular} \\ \hline
\textbf{0.7,2.3}            & \begin{tabular}[c]{@{}c@{}}8.89\\ $\pm$ 0.05\end{tabular} & \begin{tabular}[c]{@{}c@{}}0.19\\ $\pm$ 0.23\end{tabular}  & \begin{tabular}[c]{@{}c@{}}60.80\\ $\pm$ 1.01\end{tabular} & \begin{tabular}[c]{@{}c@{}}60.80\\ $\pm$ 1.01\end{tabular} & \begin{tabular}[c]{@{}c@{}}143.52\\ $\pm$ 2.47\end{tabular}     & \begin{tabular}[c]{@{}c@{}}59.95\\ $\pm$ 0.99\end{tabular} & \begin{tabular}[c]{@{}c@{}}555.71\\ $\pm$ 6.80\end{tabular} \\ \hline
\textbf{0.9,2.1}            & \begin{tabular}[c]{@{}c@{}}8.84\\ $\pm$ 0.05\end{tabular} & \begin{tabular}[c]{@{}c@{}}0.28 \\ $\pm$ 0.20\end{tabular} & \begin{tabular}[c]{@{}c@{}}41.33\\ $\pm$ 0.85\end{tabular} & \begin{tabular}[c]{@{}c@{}}41.33\\ $\pm$ 0.85\end{tabular} & \begin{tabular}[c]{@{}c@{}}116.86\\ $\pm$ 2.35\end{tabular}     & \begin{tabular}[c]{@{}c@{}}39.06\\ $\pm$ 0.83\end{tabular} & \begin{tabular}[c]{@{}c@{}}544.53\\ $\pm$ 7.41\end{tabular} \\ \hline
\textbf{1.1,1.9}            & \begin{tabular}[c]{@{}c@{}}9.18\\ $\pm$ 0.06\end{tabular} & \begin{tabular}[c]{@{}c@{}}1.16\\ $\pm$ 0.38\end{tabular}  & \begin{tabular}[c]{@{}c@{}}19.25\\ $\pm$ 0.70\end{tabular} & \begin{tabular}[c]{@{}c@{}}19.25\\ $\pm$ 0.70\end{tabular} & \begin{tabular}[c]{@{}c@{}}85.08\\ $\pm$ 2.01\end{tabular}      & \begin{tabular}[c]{@{}c@{}}15.89\\ $\pm$ 0.67\end{tabular} & \begin{tabular}[c]{@{}c@{}}514.42\\ $\pm$ 7.79\end{tabular} \\ \hline
\textbf{1.3,1.7}            & \begin{tabular}[c]{@{}c@{}}9.30\\ $\pm$ 0.06\end{tabular} & \begin{tabular}[c]{@{}c@{}}3.85\\ $\pm$ 0.52\end{tabular}  & \begin{tabular}[c]{@{}c@{}}4.34\\ $\pm$ 0.53\end{tabular}  & \begin{tabular}[c]{@{}c@{}}4.34\\ $\pm$ 0.53\end{tabular}  & \begin{tabular}[c]{@{}c@{}}64.23\\ $\pm$ 1.75\end{tabular}      & \begin{tabular}[c]{@{}c@{}}4.96\\ $\pm$ 0.52\end{tabular}  & \begin{tabular}[c]{@{}c@{}}503.62\\ $\pm$ 7.76\end{tabular} \\ \hline
\textbf{1.5,1.5}            & \begin{tabular}[c]{@{}c@{}}8.77\\ $\pm$ 0.05\end{tabular} & \begin{tabular}[c]{@{}c@{}}0.04\\ $\pm$ 0.33\end{tabular}  & \begin{tabular}[c]{@{}c@{}}0.04\\ $\pm$ 0.33\end{tabular}  & \begin{tabular}[c]{@{}c@{}}0.04\\ $\pm$ 0.33\end{tabular}  & \begin{tabular}[c]{@{}c@{}}58.41\\ $\pm$ 1.72\end{tabular}      & \begin{tabular}[c]{@{}c@{}}16.41\\ $\pm$ 0.64\end{tabular} & \begin{tabular}[c]{@{}c@{}}537.89\\ $\pm$ 8.36\end{tabular} \\ 
\hhline{=|=======}
\textbf{Avg. Opt. Gap (\%)} & -                                                         & \begin{tabular}[c]{@{}c@{}}0.72\\ $\pm$ 0.10\end{tabular}  & \begin{tabular}[c]{@{}c@{}}31.85\\ $\pm$ 0.53\end{tabular} & \begin{tabular}[c]{@{}c@{}}31.85\\ $\pm$ 0.53\end{tabular} & \begin{tabular}[c]{@{}c@{}}103.54\\ $\pm$ 0.94\end{tabular}     & \begin{tabular}[c]{@{}c@{}}32.84\\ $\pm$ 0.49\end{tabular} & \begin{tabular}[c]{@{}c@{}}527.90\\ $\pm$ 2.24\end{tabular} \\ \hline 
\end{tabular}}

\raggedright
\tabnote{Other system parameters are fixed at $\mu = (1,1)$, $a = (0.0167,0.0167)$, $h=(1.5,1)$, $b=(0,0)$, $C=4$, $\kappa=(30,30)$.}
\end{table}


\subsection{Structure of the optimal policy}\label{sec:insights}
In this section, we examine the structure of the optimal policy. With service slowdown, especially when the system exhibits meta-stability, a good policy balances two objectives: (1) reducing holding cost by giving priority to the class with a higher cost of waiting and a larger service rate, and (2) preventing the system from reaching the bad equilibrium region, e.g., by giving priority to the class that is subject to higher slowdown. 
For standard queueing models without service slowdown, myopically maximizing the instantaneous cost reduction rate often achieves good performance, as demonstrated through the well-known $c\mu$ rule. 
However, in the presence of service slowdown, myopically maximizing the instantaneous cost reduction rate may drive the system to a bad equilibrium regime. 
Thus, an optimal policy tends to give more weight to maximize the myopic cost reduction rate when either the system is not at much risk of converging to a bad equilibrium or is already trapped in a bad equilibrium region. On the other hand, the optimal policy tends to put more weight on preventing the system from moving to a bad equilibrium region when that risk is high but still controllable. In the following, we illustrate these points through numerical examples. 



\begin{figure}[]
     \centering
     \begin{subfigure}[b]{0.285\textwidth}
         \centering
         \includegraphics[width=\textwidth]{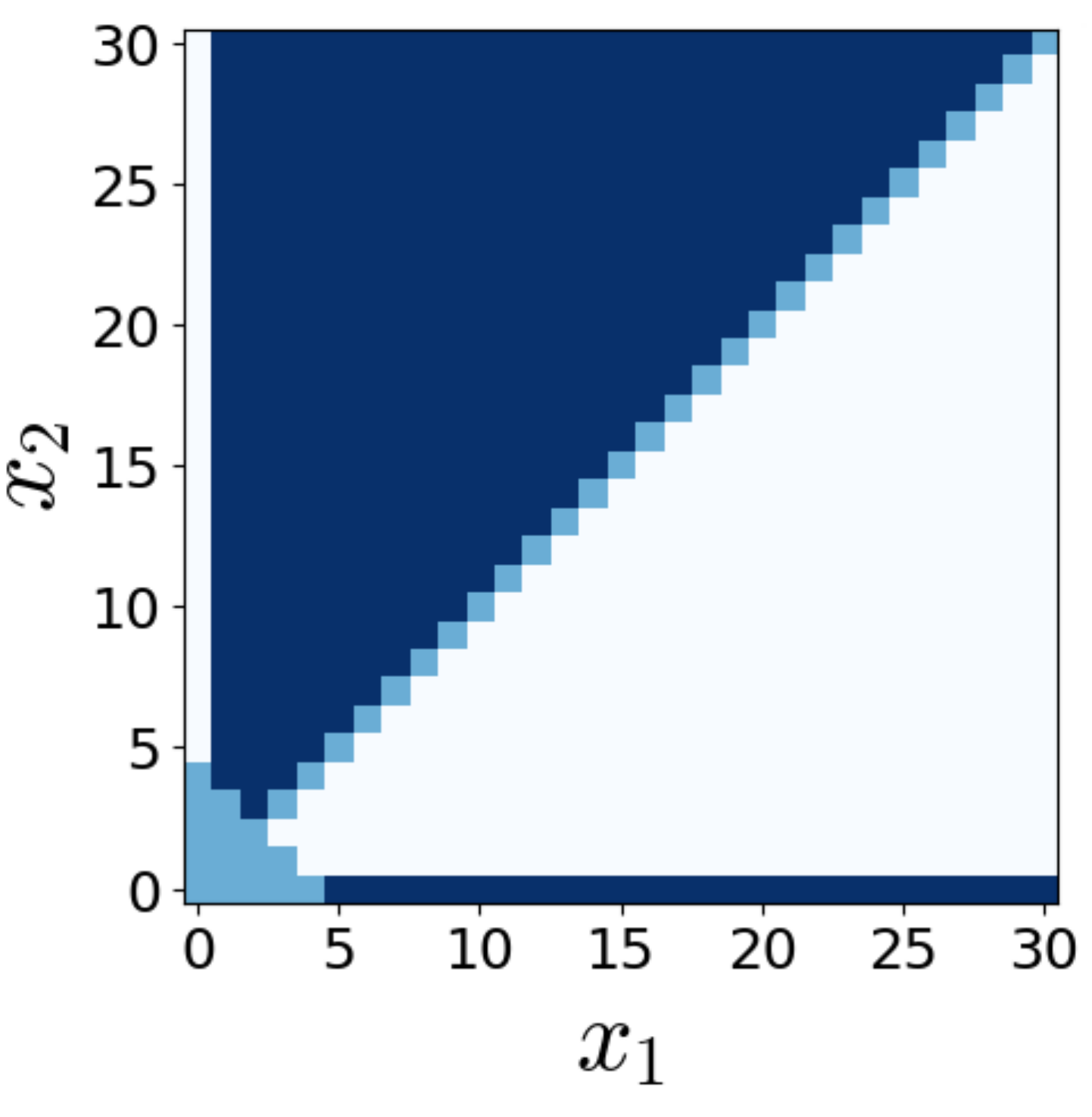}
         \caption{Worst-case load = 1}
         \label{fig:observation_1_1}
     \end{subfigure}
    \begin{subfigure}[b]{0.41\textwidth}
         \centering
         \includegraphics[width=\textwidth]{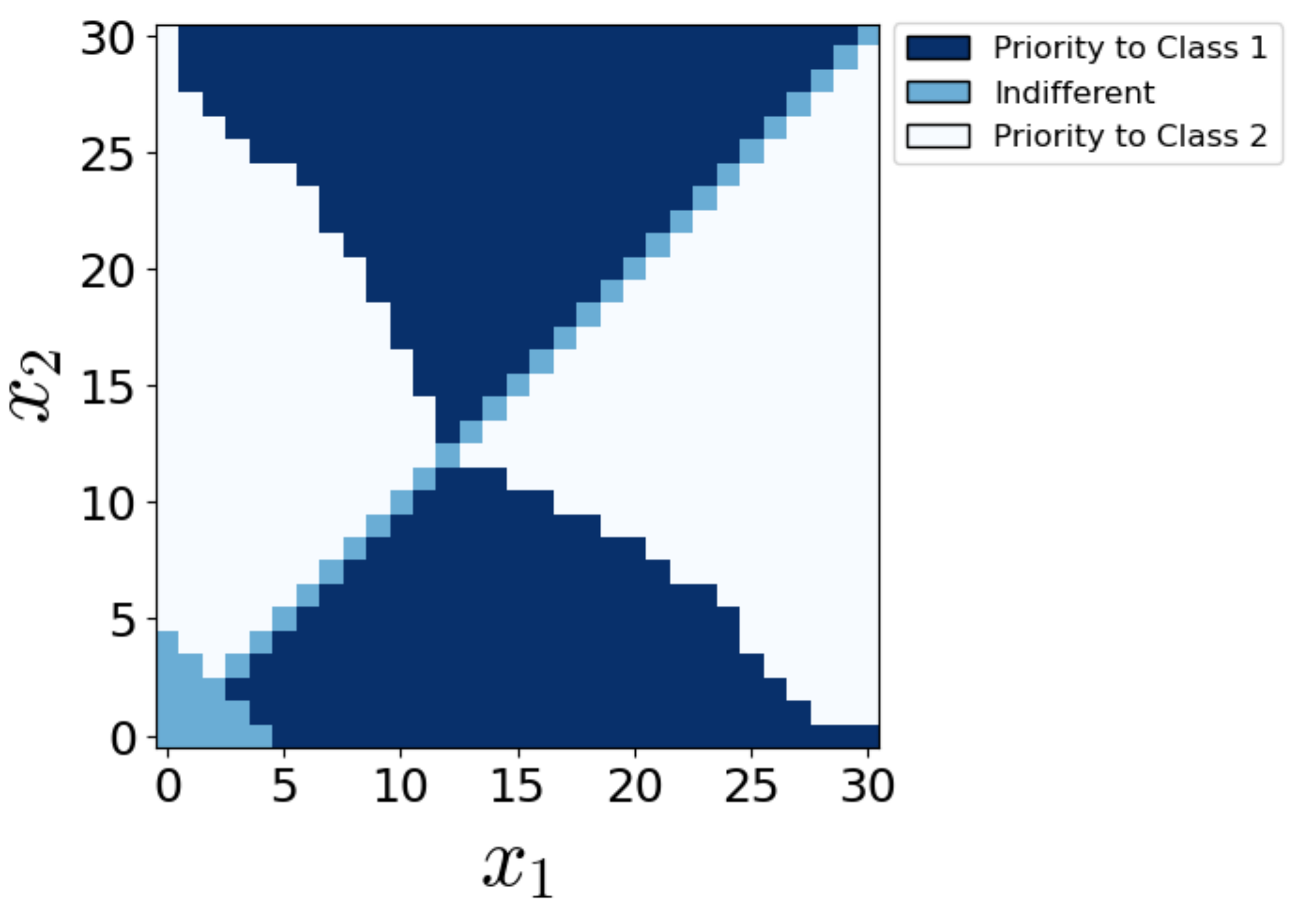}
         \caption{Worst-case load = 2}
     \end{subfigure}
     \caption{Optimal policies under different worst-case loads. Other parameters are fixed at $\lambda = (1.5,1.5), \mu = (1,1), \kappa = (30,30), h = (1,1), b = (0,0), C=4, \bar{f}_i(x)= \mu_i - a_i x_i$.}
     \label{fig:observation_1}
     \end{figure}

Figure \ref{fig:observation_1} illustrates the optimal policy for a symmetric system with two different worst-case system loads of $1$ and $2$. When the worst-case load is 1 (Figure \ref{fig:observation_1} (a)), due to blocking, there is no risk of converging to a bad equilibrium. Thus, the optimal policy reduces the holding cost at the maximum rate by prioritizing the class with a larger service rate, i.e., the class with fewer customers. When the worst-case load is $2$ (Figure \ref{fig:observation_1} (b)), the system can converge to a bad equilibrium. In this case, when the system is near the origin, the optimal policy minimizes the risk of converging to a bad equilibrium by prioritizing the class at a higher risk of deterioration, i.e., the class with more customers. However, when the system is already around a bad equilibrium, the optimal policy myopically maximizes the cost reduction rate by prioritizing the class with larger service rate.

\begin{figure}[]
     \centering
     \begin{subfigure}[b]{0.285\textwidth}
         \centering
         \includegraphics[width=1.03\textwidth]{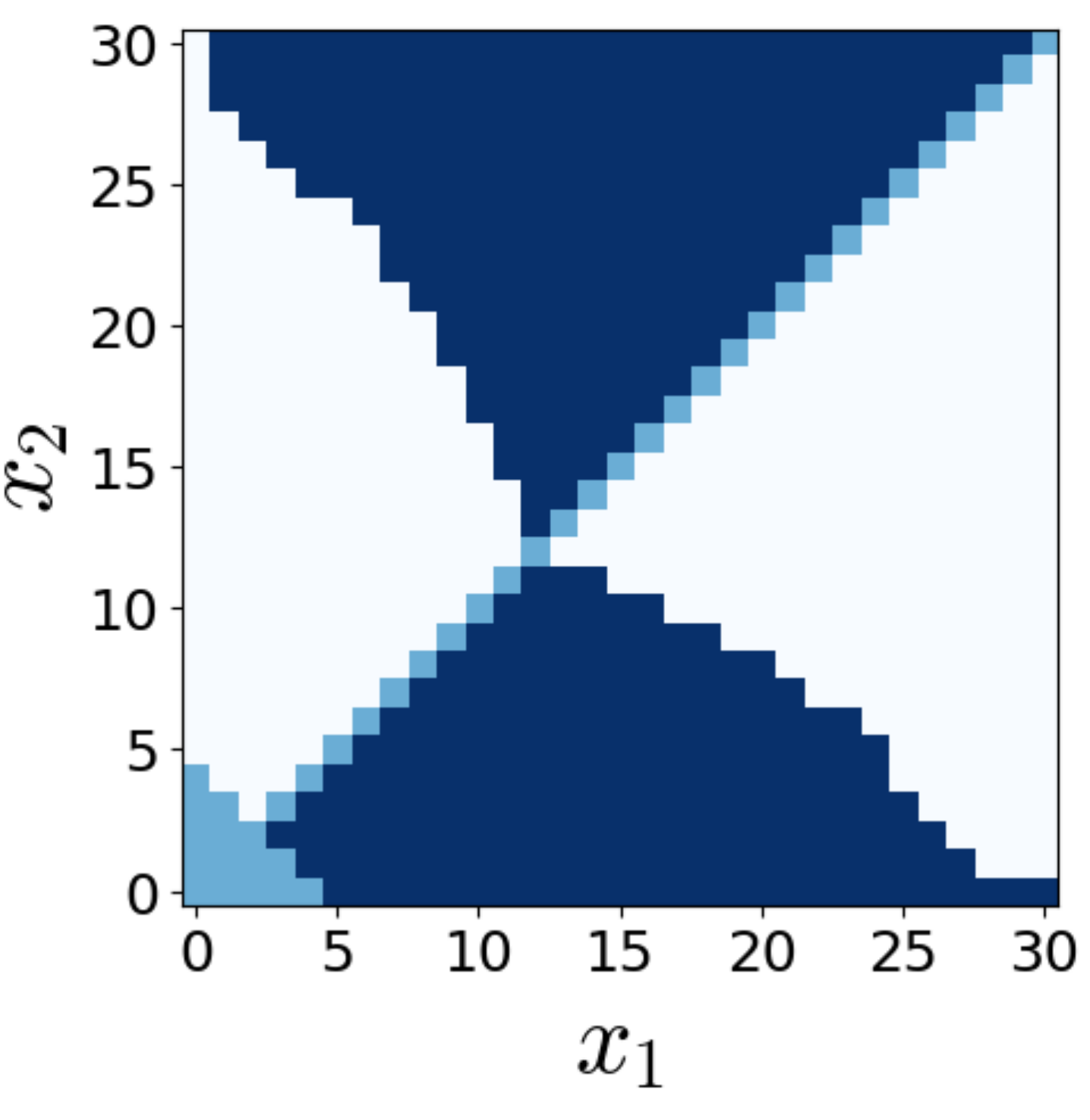}
         \caption{$a=( 0.02083, 0.02083)$ \\ ($a_2-a_1 = 0$)}
     \end{subfigure}
    \begin{subfigure}[b]{0.285\textwidth}
         \centering
         \includegraphics[width=1.03\textwidth]{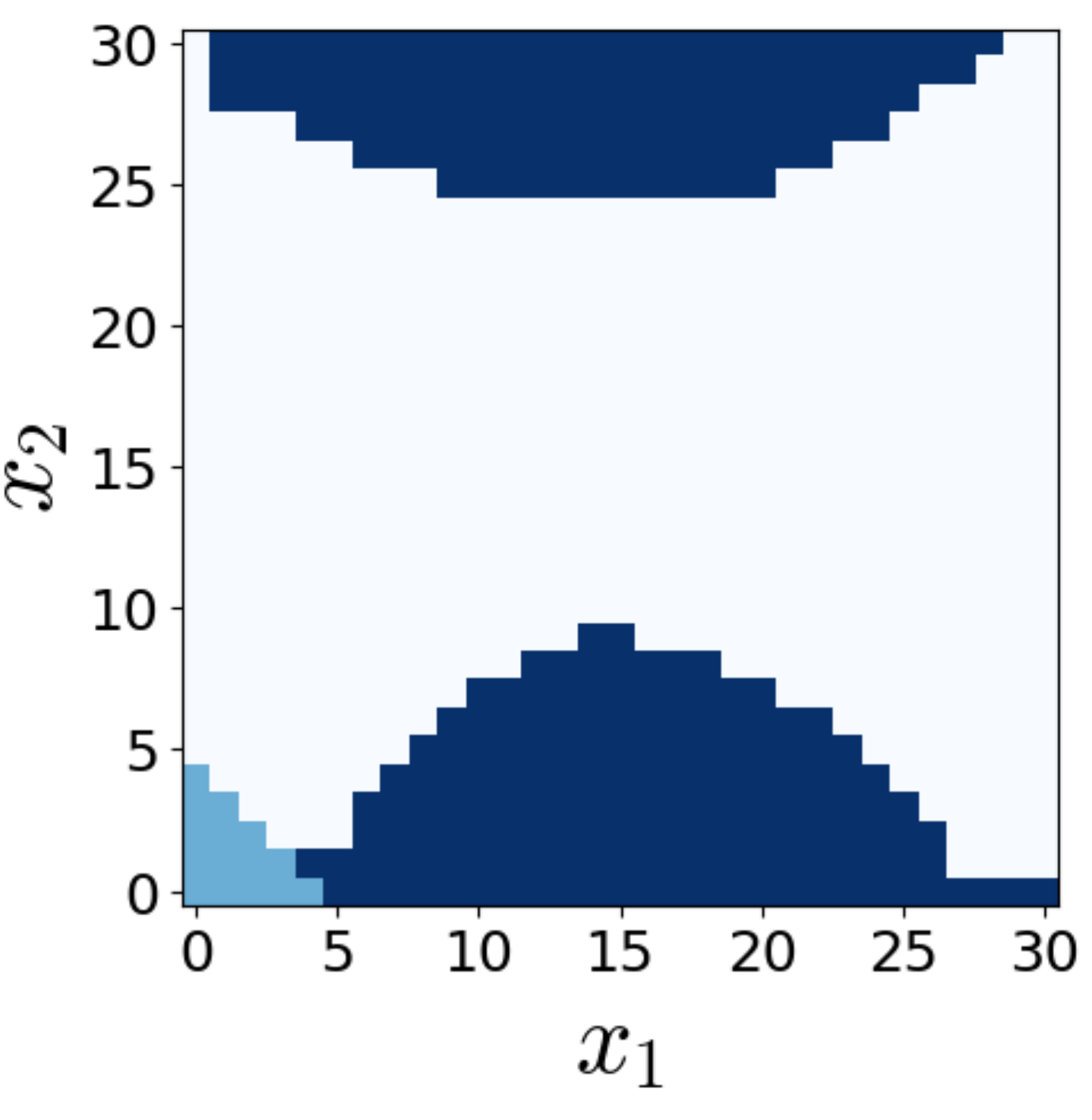}
         \caption{$a =( 0.02031,  0.02131)$ \\ ($a_2-a_1 = 0.001$)}
     \end{subfigure}
    \begin{subfigure}[b]{0.41\textwidth}
         \centering
         \includegraphics[width=1.03\textwidth]{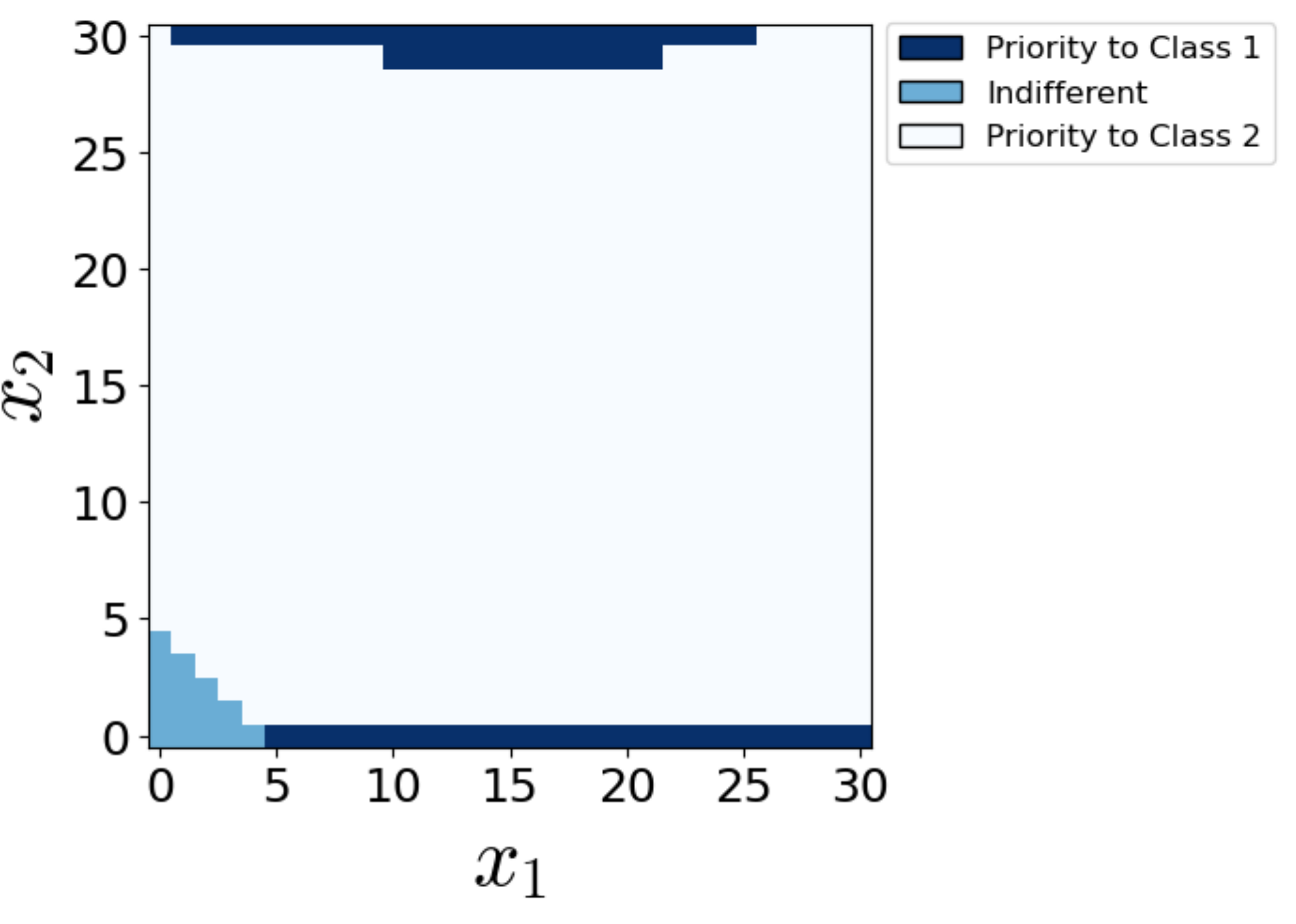}
         \caption{$a =( 0.01408,  0.02408)$ \\ ($a_2-a_1 = 0.01$)}
     \end{subfigure}
     \caption{Optimal policies under different slowdown rates in an overloaded system (Worst-case load is 2). Other parameters are fixed at $\lambda = (1.5,1.5), \mu = (1,1), \kappa = (30,30), h = (1,1), b = (0,0), C=4, \bar{f}_i(x)= \mu_i - a_i x_i$.}
     \label{fig:observation_2_2}
     \end{figure}

Next, we fix the worst-case system load equal to 2 and increase the difference between the slowdown rates of the two classes, i.e., $a_2-a_1$. Figure \ref{fig:observation_2_2} illustrates the optimal policy for three different slow-down rate values. We observe that as the difference between $a_2$ and $a_1$ increases, class 2 is at a higher risk of deterioration than class 1 and thus receives priority in a larger region. 

\begin{figure}[]
     \centering
     \begin{subfigure}[b]{0.285\textwidth}
         \centering
         \includegraphics[width=1.03\textwidth]{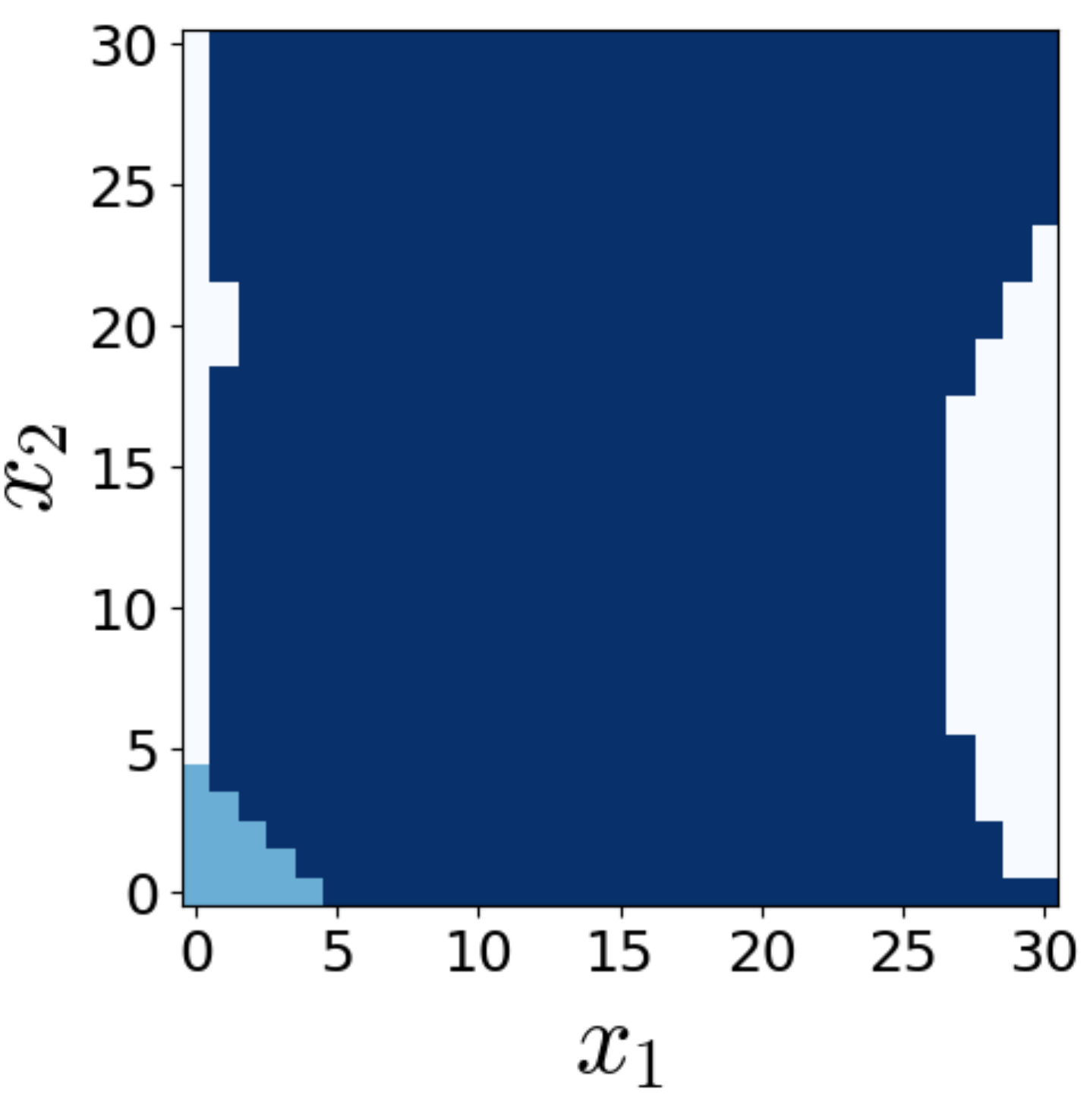}
         \caption{$\lambda = (1.5,1.5)$}
         \label{fig:observation_3_1}
     \end{subfigure}
    \begin{subfigure}[b]{0.285\textwidth}
         \centering
         \includegraphics[width=1.03\textwidth]{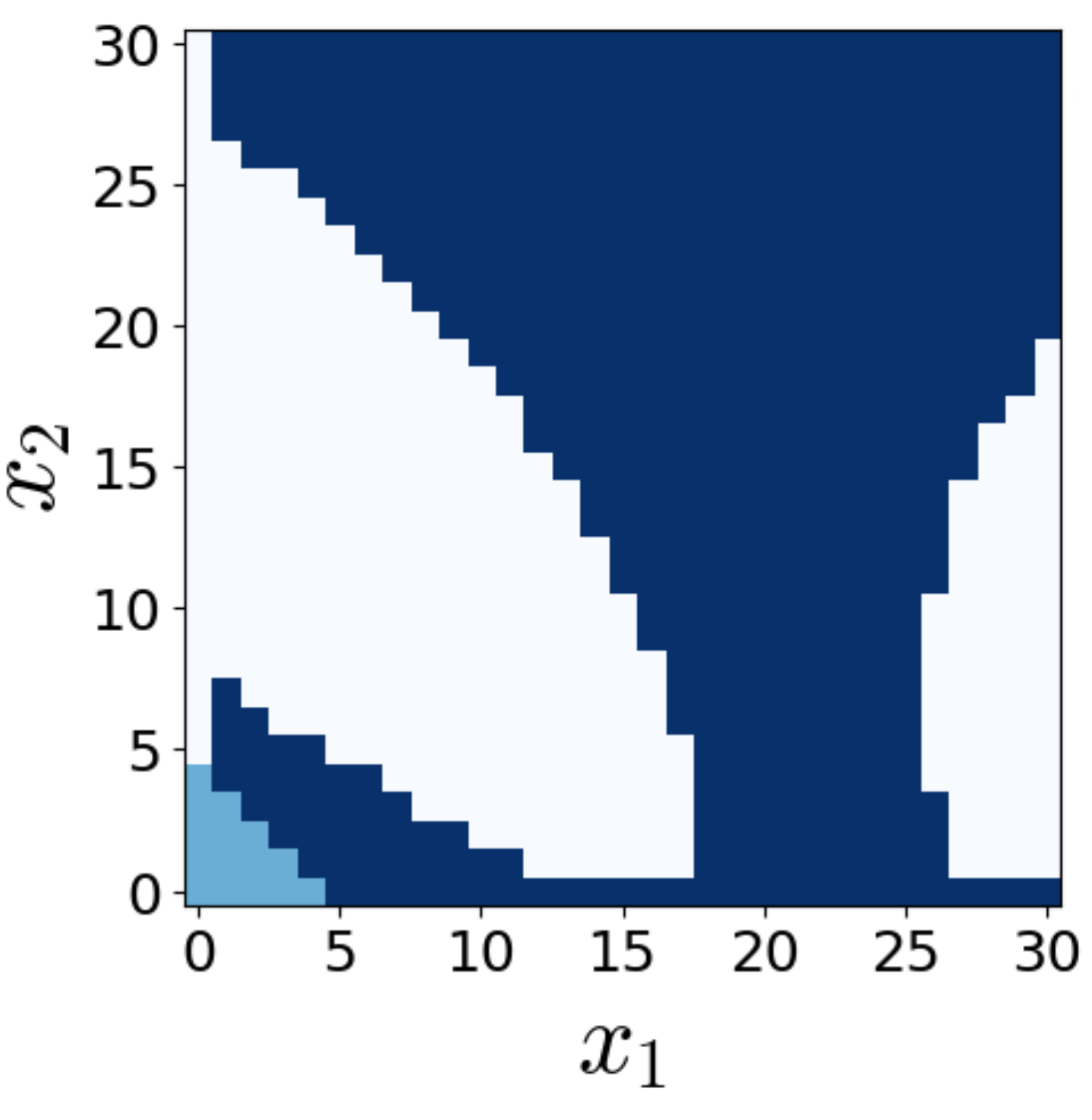}
         \caption{$\lambda= (1.2, 1.8)$}
     \end{subfigure}
    \begin{subfigure}[b]{0.41\textwidth}
         \centering
         \includegraphics[width=1.03\textwidth]{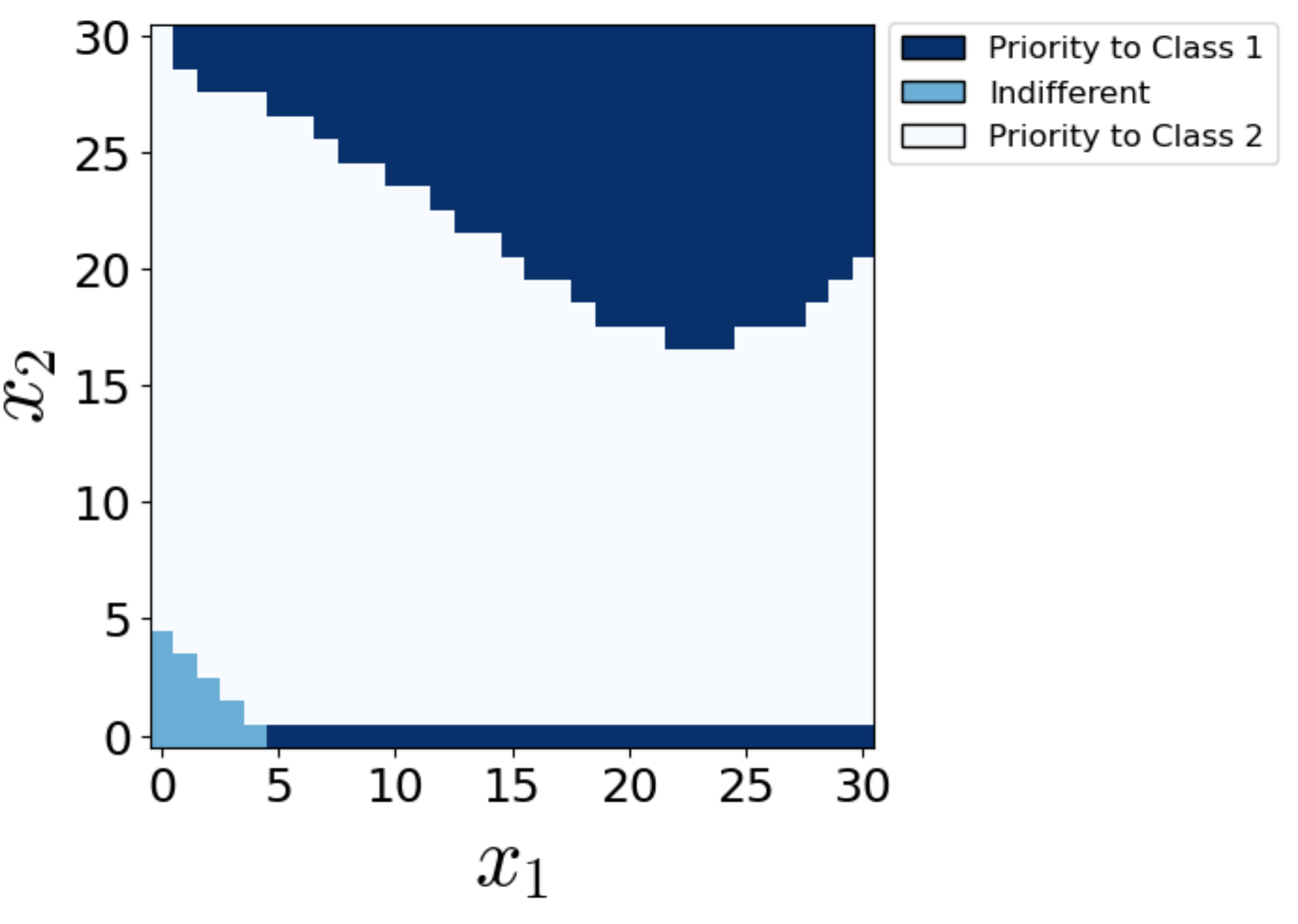}
         \caption{$\lambda=(0.8, 2.2)$}
     \end{subfigure}
     \caption{Optimal policies under different arrival rates. Other parameters are fixed at $\mu=(1,1), a = (0.0167,0167), \kappa = (30,30), h = (1.5,1), b = (0,0), C=4, \bar{f}_i(x)= \mu_i - a_i x_i$.}
     \label{fig:observation_3}
     \end{figure}

We observe that optimal scheduling under service slowdown differs from scheduling without slowdown with respect to teo important aspects. First, the optimal policy for systems with service slowdown is affected by the arrival rates. 
Figure \ref{fig:observation_3} illustrates the optimal policies under different arrival rates. Note that the alterations in arrival rates are made such that the best- and worst-case system loads remain the same, at 0.75 and 2, respectively. When the arrival rates are equal (Figure \ref{fig:observation_3_1}), the optimal policy tends to prioritize class 1 throughout since class 1 has a higher holding cost. However, as the arrival rate of class 2 increases, the optimal policy allocates more priority to class 2. 
Intuitively, when class 2 experiences service slowdown, it will have a larger impact on the system load due to the larger arrival rate and thus increase the risk of reaching a bad equilibrium.
This demonstrates that even with the same worst-case system load, arrival rates could have a significant effect on the optimal policy. 

\begin{figure}[]
     \centering
     \begin{subfigure}[b]{0.285\textwidth}
         \centering
         \includegraphics[width=1.03\textwidth]{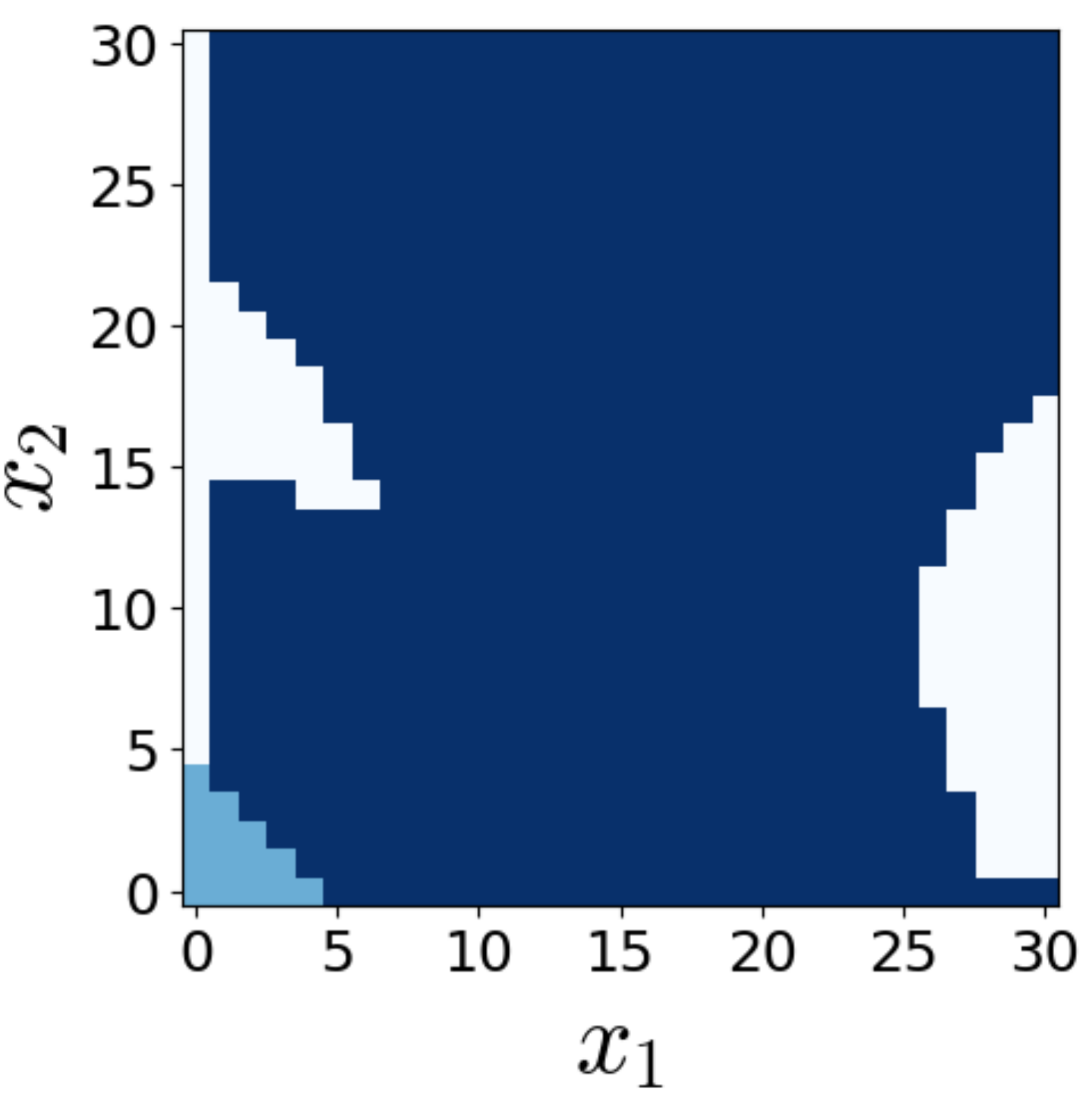}
         \caption{$b= (0,0)$}
         \label{fig:observation_4_1}
     \end{subfigure}
    \begin{subfigure}[b]{0.285\textwidth}
         \centering
         \includegraphics[width=1.03\textwidth]{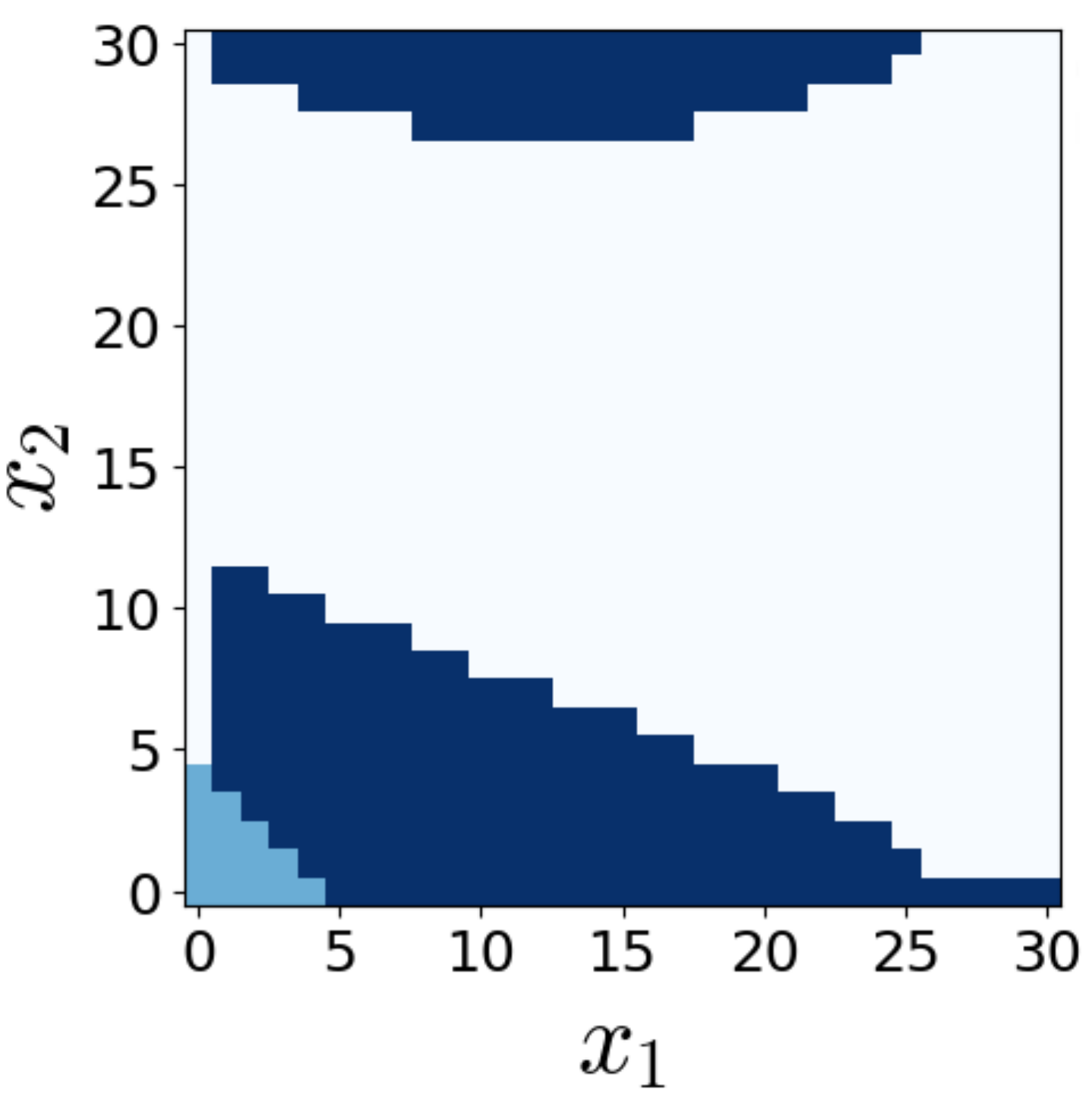}
         \caption{$b= (0,10)$}
     \end{subfigure}
    \begin{subfigure}[b]{0.41\textwidth}
         \centering
         \includegraphics[width=1.03\textwidth]{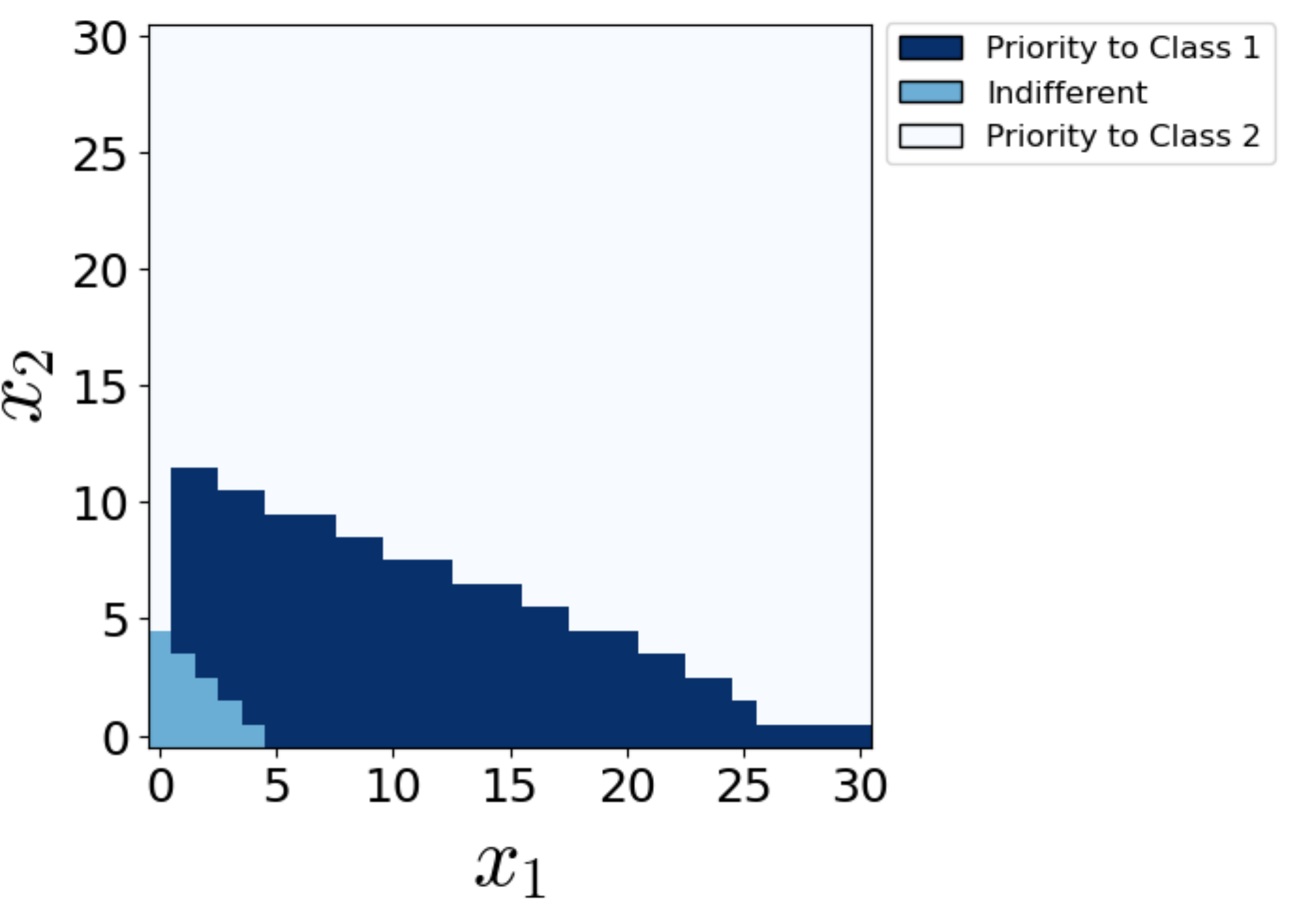}
         \caption{$b= (0,100)$}
     \end{subfigure}
     \caption{An example with $\lambda = (1.5,1.5), \mu = (1,1), a=(0.0103,0.0203), h = (5,1), \kappa = (30,30), \bar{f}_i(x)= \mu_i - a_i x_i$.}
     \label{fig:observation_4}
     \end{figure}

Second, as the blocking costs increase, preventing the system from converging and getting stuck at a bad equilibrium region becomes more evident in the optimal policy.
Figure \ref{fig:observation_4} illustrates optimal policies under varying blocking costs for class 2 customers. In this example, class 1 has a much higher holding cost than class 2. Thus, when the blocking cost for class 2 is low, we tend to prioritize class 1 in most regions. However, as the blocking cost of class 2 customers increases, class 2 is prioritized, especially closer to the blocking boundary.

\section{Experiments with more than two classes}\label{sec:fiveclass}
In this section, we illustrate the performance of the ADP algorithm for more than two classes. We consider a five-class system with $\kappa = (30,30,30,30,30)$, which increases the state-space in the previous experiments $27,000$-folds (from $31^2$ to $31^5$). In our experiments, we use a multinomial logistic regression to determine the priority order, based on the ranking of predicted class probabilities.

Table \ref{tab:five_class_experiment} provides the long-run average cost under the ADP and benchmark policies for six example cases with different system parameters. Due to the large state space, obtaining the optimal policy is no longer feasible. 
The table illustrates that the ADP policy yields significantly lower long-run average costs compared to the other benchmarks. Specifically, the average performance gap of benchmark policies with respect to ADP policy ranges from $31\%$ to $42\%$. As we highlight in Section \ref{sec:compare_benchmark}, the ADP algorithm yields larger improvements compared to the benchmarks when (1) the difference between slowdown rates for different classes increases; (2) the difference between arrival rates for different classes increases; and (3) one class has a significantly larger blocking cost than the others. This illustrates that our insights extend to systems with more than two classes. We also highlight that our ADP algorithm can learn well-performing policies using limited samples. Due to the large state space, in each policy iteration, we only sample $0.004\%$ of the states.

\begin{table}[]
\caption{Long-run average costs under the ADP policy and five benchmark policies for six parameter regimes.}
\label{tab:five_class_experiment}
\centering
\begin{tabular}{ccccccc}
\hline
\textbf{Cases}  & \textbf{ADP}                                               & \textbf{$h \mu(0)$}                                        & \textbf{$h \mu(x)$}                                        & \textbf{\begin{tabular}[c]{@{}c@{}}Max\\ Pressure\end{tabular}} & \textbf{SQF} & \textbf{LQF} \\ \hline
\textbf{Case 1} & \begin{tabular}[c]{@{}c@{}}37.93\\ $\pm$0.03\end{tabular}  & \begin{tabular}[c]{@{}c@{}}58.95\\ $\pm$0.05\end{tabular}  & \begin{tabular}[c]{@{}c@{}}59.15\\ $\pm$0.08\end{tabular}  & \begin{tabular}[c]{@{}c@{}}87.76\\ $\pm$0.05\end{tabular}       & \begin{tabular}[c]{@{}c@{}}58.43\\ $\pm$0.10\end{tabular}        & \begin{tabular}[c]{@{}c@{}}65.23\\ $\pm$0.01\end{tabular}       \\ \hline
\textbf{Case 2} & \begin{tabular}[c]{@{}c@{}}70.80\\ $\pm$1.32\end{tabular}  & \begin{tabular}[c]{@{}c@{}}112.98\\ $\pm$0.38\end{tabular} & \begin{tabular}[c]{@{}c@{}}86.47\\ $\pm$0.91\end{tabular}  & \begin{tabular}[c]{@{}c@{}}120.03\\ $\pm$0.63\end{tabular}      & \begin{tabular}[c]{@{}c@{}}78.62\\ $\pm$0.45\end{tabular}        & \begin{tabular}[c]{@{}c@{}}85.12\\ $\pm$0.22\end{tabular}       \\ \hline
\textbf{Case 3} & \begin{tabular}[c]{@{}c@{}}46.89\\ $\pm$1.85\end{tabular}  & \begin{tabular}[c]{@{}c@{}}51.77\\ $\pm$0.39\end{tabular}  & \begin{tabular}[c]{@{}c@{}}51.79\\ $\pm$0.39\end{tabular}  & \begin{tabular}[c]{@{}c@{}}50.43\\ $\pm$0.22\end{tabular}       & \begin{tabular}[c]{@{}c@{}}63.56\\ $\pm$0.45\end{tabular}        & \begin{tabular}[c]{@{}c@{}}67.89\\ $\pm$0.36\end{tabular}       \\ \hline
\textbf{Case 4} & \begin{tabular}[c]{@{}c@{}}11.78\\ $\pm$0.12\end{tabular}  & \begin{tabular}[c]{@{}c@{}}34.15\\ $\pm$0.13\end{tabular}  & \begin{tabular}[c]{@{}c@{}}33.86\\ $\pm$0.15\end{tabular}  & \begin{tabular}[c]{@{}c@{}}24.76\\ $\pm$0.44\end{tabular}       & \begin{tabular}[c]{@{}c@{}}29.79\\ $\pm$0.34\end{tabular}        & \begin{tabular}[c]{@{}c@{}}32.54\\ $\pm$0.22\end{tabular}       \\ \hline
\textbf{Case 5} & \begin{tabular}[c]{@{}c@{}}10.73\\ $\pm$0.07\end{tabular}  & \begin{tabular}[c]{@{}c@{}}10.02\\ $\pm$0.02\end{tabular}  & \begin{tabular}[c]{@{}c@{}}11.79\\ $\pm$0.35\end{tabular}  & \begin{tabular}[c]{@{}c@{}}20.67\\ $\pm$0.56\end{tabular}       & \begin{tabular}[c]{@{}c@{}}38.16\\ $\pm$0.48\end{tabular}        & \begin{tabular}[c]{@{}c@{}}42.22\\ $\pm$0.32\end{tabular}       \\ \hline
\textbf{Case 6} & \begin{tabular}[c]{@{}c@{}}200.68\\ $\pm$1.17\end{tabular} & \begin{tabular}[c]{@{}c@{}}394.92\\ $\pm$2.14\end{tabular} & \begin{tabular}[c]{@{}c@{}}384.90\\ $\pm$2.37\end{tabular} & \begin{tabular}[c]{@{}c@{}}387.98\\ $\pm$2.20\end{tabular}      & \begin{tabular}[c]{@{}c@{}}250.03\\ $\pm$3.86\end{tabular}       & \begin{tabular}[c]{@{}c@{}}236.27\\ $\pm$3.95\end{tabular}      \\ \hline
\end{tabular}

\raggedright
\tabnote{
  \textbf{Case 1:} $\lambda=(0.5,1,1.5,2,2.5),$ $\mu=(1,1,1,1,1),$ $a = (0.033,0.0134,0.0134,0.028,0.028),$ $h=(1,1.01,1.02,1.03,1.04),$ $b=(0,0,0,0,0),$ $\kappa = (30,30,30,30,30),$ $C=8$. 
  \textbf{Case 2:} $\lambda=(0.5,1,1.5,2,2.5),$ $\mu=(1,1,1,1,1),$ $a = (0.03,0.03,0.03,0.03,0.03),$ $h=(1,1.1,1.2,1.3,1.4),$ $b=(0,0,0,0,0),$ $\kappa = (30,30,30,30,30),$ $C=8$. 
  \textbf{Case 3:} $\lambda=(0.5,1,1.5,2,2.5),$ $\mu=(1,1,1,1,1),$ $a = (0.03,0.03,0.03,0.03,0.03),$ $h=(1.4,1.3,1.2,1.1,1),$ $b=(0,0,0,0,0),$ $\kappa = (30,30,30,30,30),$ $C=8$. 
  \textbf{Case 4:} $\lambda=(1.2,1.2,1.2,1.2,1.2),$ $\mu=(1,1,1,1,1),$ $a = (0.032,0.03,0.028,0.026,0.024),$ $h=(1,1.1,1.2,1.3,1.4),$ $b=(0,0,0,0,0),$ $\kappa = (30,30,30,30,30),$ $C=8$. 
  \textbf{Case 5:} $\lambda=(1.2,1.2,1.2,1.2,1.2),$ $\mu=(1,1,1,1,1),$ $a = (0.032,0.03,0.028,0.026,0.024),$ $h=(1.4,1.3,1.2,1.1,1),$ $b=(0,0,0,0,0),$ $\kappa = (30,30,30,30,30),$ $C=8$. 
  \textbf{Case 6:} $\lambda=(1.6,1.6,1.6,1.6,1.6),$ $\mu=(1,1,1,1,1),$ $a = (0.032,0.032,0.032,0.032,0.032),$ $h=(1.4,1.3,1.2,1.1,1),$ $b=(0,25,50,75,100),$ $\kappa = (30,30,30,30,30),$ $C=8$. }
\end{table}

\section{Case Study}\label{sec:case_study}
In this section, we present a case study on scheduling admissions to rehabilitation care (rehab for short). The case study is based on analyses of data from a Canadian hospital studied in \cite{gorgulu2022}. The hospital provides both acute and rehab care at the same site. Patients often experience long admission delays for rehab due to the limited rehab capacity. In particular, patients wait an average of $2.65$ days before admission to rehab. \cite{gorgulu2022} focus on the two largest patient categories in acute care: Neurology/Musculoskeletal (Neuro/MSK) and Medicine. Using an Instrumental Variable (IV) approach, they find that delays in admission to rehabilitation increase the rehab length-of-stay (service time) for both patient categories, and negatively impact the patients' functional independence measures (FIM) at both admissions to and discharge from rehab. (FIM is a standard measure for evaluating the functional capacities of rehab patients which takes values between $18$ and $126$, see, e.g., \citealt{linacre1994structure}.) Further, these impacts are heterogeneous for the two patient categories. In particular, the functionality of the Neuro/MSK patients is more sensitive to delays than Medicine patients.



\subsection{Model Calibration}
There are two patient classes $I = \{N$(Neuro/MSK), $M$(Medicine)$\}$. Arrivals to the queue correspond to bed requests from acute care to rehab and are modeled as a Poisson Process with a class-dependent rate $\lambda = (1.04,2.29)$. There are 79 identical servers (beds) and the system has a capacity limit of $\kappa_N= \kappa_M = 110$ for each patient category. Service times correspond to LOSs in rehab. We assume that service times follow a Log-Normal distribution where the mean is a class-dependent function of the patient's wait time (admission delay). 
In particular, the average service time for class $i$ patients is assumed to take the form $f_i(w) = \alpha_i + \beta_i (w \wedge \tilde{w})$ where $w$ is the wait time of the patient at the time of admission and $\tilde{w}$ is the wait threshold for the service time increase, $\alpha_i$ is the base service time and $\beta_i$ measure the sensitivity to delay. Since $f_i(w)$ is a convex increasing function, the optimal solution to \eqref{eq:optimize} can be easily computed. The structure of $f_i(w)$ leads to the following service time distribution for class $i$ patients:
\begin{align*}
    \text{Log-normal}\left(\mu = \log(\alpha_i + \beta_i (w \wedge \tilde{w}))-\frac{1}{2}\sigma^2,\ 
    \sigma^2=  \log\left(1+\left(\frac{\gamma_i}{\alpha_i+\beta_i (w \wedge \tilde{w})}\right)^2\right)\right).
\end{align*}

For Neuro/MSK and Medicine patients, the estimated $\alpha_i$, $\beta_i$, and $\gamma_i$ parameters can be found in Table \ref{tab:estimated_parameters}. We also set $\tilde{w}=9$. The details of the estimation can be found in Section \ref{sec:estimation_details}. 

The status-quo admission policy is an ad-hoc policy. Therefore, we assume a probabilistic routing policy is employed and calibrate the routing probabilities to match the empirical average waiting times in the data. Figure \ref{fig:probability_case_study} illustrates the average waiting times for Neuro/MSK and Medicine patients for different values of probability of choosing a Neuro/MSK patient next when there are patients from both categories waiting to be admitted (denoted by $p$). We observe that the empirical average waiting times correspond to $p=0.46$
\endnote{Note that the arrival rate of Medicine patients is significantly larger than Neuro/MSK patients, i.e., $\frac{\lambda_{M}}{\lambda_{N}}>2$. Therefore, waiting times of Neuro/MSK and Medicine patients can be balanced by setting $p =  \frac{\lambda_{N}}{\lambda_{N}+\lambda_{M}} \approx 0.312$. Since $p=0.46>0.312$, it indicates that Neuro/MSK patients receive a higher priority.}. Note that due to the service slowdown, the average waiting time is not monotone in $p$ for both categories. We set $h = (1.5,1)$ to reflect that the discharge FIM score for Neuro/MSK patients is impacted 1.5 times more than that of Medicine patients. We apply the non-preemptive version of our algorithm since preemption is not allowed in our setting. We discuss the performance of the policy learned by the preemptive algorithm in Section \ref{ap:case_pre}.
\begin{table}[]
\centering
\caption{Estimated $\alpha$, $\beta$ and $\gamma$ parameters for Medicine and Neuro/MSK patients.}
\label{tab:estimated_parameters}
\begin{tabular}{lccc}
\hline
 \textbf{Category}                  & $\alpha_i$ & $\beta_i$ & $\gamma_i$ \\ \hline
\textbf{Neuro/MSK} & $18.301$   & $2.2143$  & $12.020$   \\
\textbf{Medicine}  & $18.190$   & $1.0233$  & $11.371$   \\ \hline
\end{tabular}
\end{table}

\subsection{Results}

\begin{figure}[]
    \begin{subfigure}[b]{.49\textwidth}
        \centering
    \includegraphics[scale=0.55]{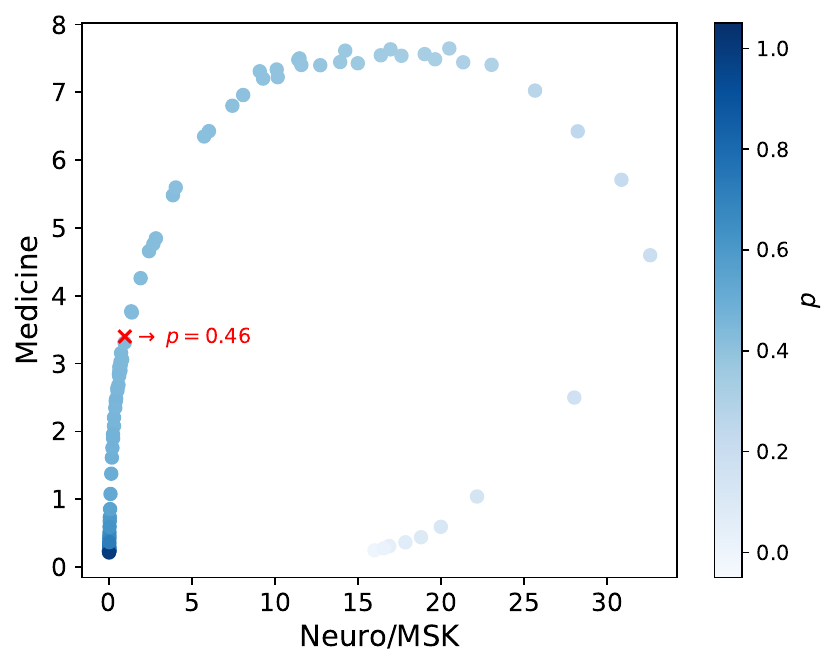}
    \caption{$p$ represents the probability of selecting Neuro/MSK over Medicine patients. }
    \label{fig:probability_case_study}
     \end{subfigure}\hfill
    \begin{subfigure}[b]{0.49\textwidth}
        \centering
    \includegraphics[scale=0.55]{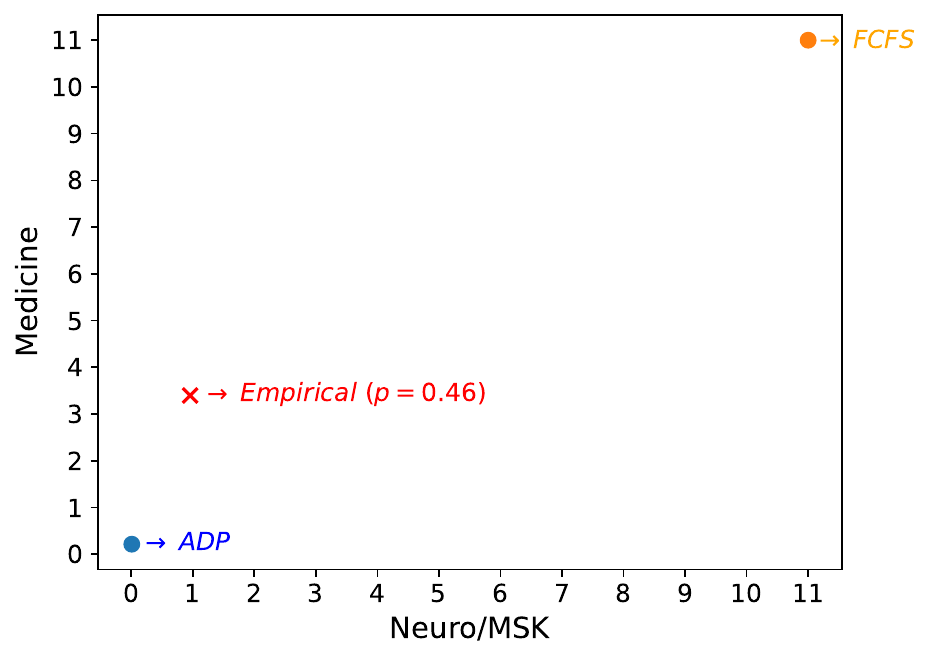}
    \caption{Comparison between  empirical, ADP and first-come-first-served (FCFS) policies.}
    \label{fig:case_study_base}
     \end{subfigure}
    \caption{Average waiting times of Neuro/MSK and Medicine patients. The point represented by x ($p=0.46$) denotes the historical waiting times.}
    \label{fig:case_study_base_all}
\end{figure}

We first focus on the estimated setting with $\beta = (2.2143,1.023)$.  Figure \ref{fig:case_study_base} illustrates the average waiting times for Neuro/MSK and Medicine patients under the ADP and FCFS policies as well as the historical values (i.e., the hospital's current policy). The results show that our ADP policy almost eliminates the waiting times of both classes. This translates to a $2.5$-day reduction in average waiting times compared to the current practice. In terms of individual categories, ADP policy respectively creates $0.95$ and $3.18$-days reductions for Neuro/MSK and Medicine patients which translates to $2.171$ $(2.3\%)$ and $4.934$ $(5.5\%)$ improvement in the discharge FIM score. We also compare the performance of the ADP algorithm to the FCFS policy. Scheduling patients according to the FCFS policy yields on average $11.03$ days of wait which is $10.756$-days larger than the average waiting times under the ADP policy. Thus, compared to FCFS, the ADP policy respectively achieves a $25.40$ $(26.14\%)$ and $16.74$ $(18.71\%)$ increase in the discharge FIM scores of Neuro/MSK and Medicine patients.

We note that the ADP policy is approximately identical to a strict priority in favor of the Neuro/MSK category. While this may seem intuitive given the higher holding cost and sensitivity to delays for Neuro/MSK patients, we find that the strict priority policy continues to be optimal even when the holding cost for Medicine patients is larger. Furthermore, as evident from the comparison with the FCFS policy, we observe that assigning higher priority to Neuro/MSK can reduce waiting times for Medicine patients as well. Finally, we note that the observed improvement relies on the feasibility of assigning strict priority to Neuro/MSK patients. In reality, there may be factors that impact the prioritization decision but are not captured in our model. Nevertheless, the results suggest the potential for a significant reduction in waiting time by accounting for the effects of slowdown when scheduling patients.

\section{Conclusion}\label{sec:conclusion}
In this work, we study scheduling in multiclass queues with service slowdown. We propose a simulation-based ADP algorithm that leverages the specific structure of the problem to generate high-quality scheduling policies. Our ADP algorithm has three key features: (i) it learns the policy directly by approximating it using a classifier; (ii) it uses a coupling method to obtain low-variance estimates of the relative value function differences directly; and (iii) it utilizes adaptive sampling to allocate the sampling budget, which leads to high computational cost reduction.

Using extensive numerical experiments, we demonstrate that our proposed ADP algorithm generates near-optimal policies and outperforms well-known benchmarks, while also scaling well to handle large-scale problems. Additionally, we provide insights into the structure of the optimal policy across various parameter regimes. Our results demonstrate that an optimal policy must strike a delicate balance between immediate cost reduction and preventing the system from moving to an undesirable equilibrium, which can occur due to congestion-induced service slowdowns. Finally, we showcase the value of the proposed algorithm through a case study on scheduling patient admissions to rehabilitation care. The case study underscores both the practical value of our algorithm and the crucial role of proper scheduling when waiting times have heterogeneous impacts on service requirements.


Our work focuses on developing an algorithmic approach to scheduling under service slowdown. Our numerical analysis offers insights into the structure of the optimal policy, which can be explored analytically in future research. Moreover, our proposed ADP algorithm, especially the coupling construction, has the potential to be applied to other complex queueing control problems. Investigating its adaptability and performance in different settings remains an avenue for future work.




\begingroup \parindent 0pt \parskip 0.0ex \def\enotesize{\normalsize} \theendnotes \endgroup

\bibliographystyle{informs2014} 
\bibliography{manuscript}

\ECSwitch

\ECHead{E-Companion}

\section{Proofs}
\subsection{Proof of Proposition \ref{prop:existance_fluid}}
We establish the existence of a solution to the problem defined by \eqref{eq:fluid_dynamics} for any initial condition $\Bar{X}(0)$. Note that $\Bar{Z}(t)$ is determined by the scheduling policy which is a function of the state $\Bar{X}(t)$, therefore, $\Bar{Z}(t)$ can be equivalently represented as a function of the state $\Bar{Z}(\Bar{X}(t))$. Note that the problem consists of a system of ODEs such that they are connected through a common capacity constraint on $\Bar{Z}(\Bar{X}(t))$. Therefore, if we show that there exists a solution for a single class system for any $\bar{Z}(\bar{X})$, it implies that there exists a solution to the problem defined by \eqref{eq:fluid_dynamics}. We first drop the subscripts for the convenience of the proof and partition the state-space $[0, \kappa]$ into two distinct regions:
\begin{align}
    \mathcal{R}_1 &= \{\bar{X}:\bar{X}<\kappa\},\\
    \mathcal{R}_2 &= \{\bar{X}:\bar{X}=\kappa\}.
\end{align}
Define,
\begin{align}
    F_1(\bar{X}) &= \lambda - \bar{f}(\bar{X})\bar{Z}(\bar{X}),\\
    F_2(\bar{X}) &= -\bar{f}(\bar{X})\bar{Z}(\bar{X}), \\
    F_3(\bar{X}) &= \psi \lambda - \psi\bar{f}(\bar{X})\bar{Z}(\bar{X}) - (1-\psi)\bar{f}(\bar{X})\bar{Z}(\bar{X}) = \psi \lambda - \bar{f}(\bar{X})\bar{Z}(\bar{X})
\end{align}

where $\psi \in [0,1]$. Consider the differential inclusion, $\dot{\bar{X}} \in \mathcal{F}(\bar{X})$ defined as:
\begin{align}
    \dot{\bar{X}} \in \mathcal{F}(\bar{X}) = 
    \begin{cases}
        F_1(\bar{X}), &\text{if }\bar{X}\in \mathcal{R}_1, \\
        F_3(\bar{X}) , &\text{if }\bar{X}\in \mathcal{R}_2.
    \end{cases}
\end{align}
By \cite{filippov2013differential}, a solution to this differential inclusion is also a solution to the problem defined by \eqref{eq:fluid_dynamics}. We use Proposition 3 of \cite{cortes2008discontinuous} to show that there exists a solution to our differential inclusion $\mathcal{F}(\bar{X})$. We need to show that $\mathcal{F}(\bar{X})$ is measurable and locally essentially bounded, i.e., bounded on a bounded neighborhood of every point, excluding sets of measure zero. Note that both $F_1(\bar{X})$ and $F_3(\bar{X})$ are continuous functions and hence measurable. Moreover, $F_1(\bar{X}) \in [\lambda - \bar{f}(0)C ,\lambda]$ and $F_3(x) \in [- \bar{f}(0)C ,0]$ $\forall \bar{X} \in [0,\kappa]$ and $\forall Z\in [0,C]$ which implies that they are bounded and close. Therefore, there exists a solution for any initial condition $\bar{X}(0)$. \halmos \endproof

\subsection{Proof of Proposition \ref{prop:fluid_characterize}}
Before providing the proof of Proposition \ref{prop:fluid_characterize}, we present the following lemma.

\begin{lemma}\label{lemma:kappa_stab}
    If $\frac{\lambda}{(\mu - a C)}\leq C < \frac{\lambda}{(\mu - a \kappa)}$, then $\kappa$ is a pseudo-equilibrium of the system \eqref{eq:fluid_dynamics} and it is asymptotically stable.
\end{lemma}
\proof{Proof of Lemma \ref{lemma:kappa_stab}.}
Let $\Bar{X}(0) = \kappa$. Assume for a contradiction that $\exists t>0$ and $\kappa - \frac{\mu C -\lambda}{a C}\geq \epsilon>0$ s.t. $\Phi^{\bar{\pi}}(\kappa, t) = \kappa-\epsilon$. Since $\Bar{X}(0) = \kappa$ and $\dot{\Bar{X}}(0)<0$, then $\exists \Tilde{t} < t$ s.t.  $\Phi^{\bar{\pi}}(\kappa,\Tilde{t}) = \kappa - \Tilde{\epsilon}>\kappa-\epsilon$. Let $\dot{\Bar{X}}=\varphi(\Bar{X})$. Note that at  $\kappa - \Tilde{\epsilon}$ we have  $\varphi(\kappa - \Tilde{\epsilon}) >0$. Moreover, for any $\kappa>y>\kappa - \Tilde{\epsilon}$ we also have $\varphi(y)>0$. This means that when the system is at $\kappa - \Tilde{\epsilon}$, the flow increases until the system reaches $\kappa$. This implies that $\forall \hat{t}>\Tilde{t}$, $\Phi^{\bar{\pi}}(\kappa, \hat{t}) > \kappa - \Tilde{\epsilon} >\kappa - \epsilon$ proving that no such $\epsilon$ can exists. This leads to a contradiction validating that $\bar{x}_e=\kappa$ is an equilibrium point. Moreover, since $\varphi(\kappa)<0$, $\bar{x}_e=\kappa$ is not a stationary point, making it a pseudo equilibrium.

We next verify stability. Let $\kappa - \frac{\mu C -\lambda}{a C}\geq \epsilon>0$, then 
\begin{align*}
    \varphi(\bar{x}_e - \epsilon) &= \lambda - \left(\mu - a \kappa + a\epsilon \right)C = \lambda - \mu C + aC \kappa - C \epsilon \\
    &\geq  \lambda - \left(\mu - a \kappa + a\epsilon \right)C \geq \lambda - \mu C + aC \kappa - aC \left(\kappa - \frac{\mu C -\lambda}{a C}\right) = 0.
\end{align*}
This implies that if $||\Bar{X}(0) - \bar{x}_e|| < \epsilon$ then $||\Phi^{\bar{\pi}}(\Bar{X}(0), t) - \bar{x}_e|| < ||\Bar{X}(0) - \bar{x}_e|| $ and hence $\bar{x}_e = \kappa$ is a stable equilibrium point. \halmos \endproof

We now turn to the proof of Proposition \ref{prop:fluid_characterize}. In order to determine potential equilibrium points, we first consider the system without truncation. Under strict priority to class 1 customers, the ODE takes the following form:
\begin{align}\label{equation:two_class_fluid_1}
    \dot{\Bar{X}}_1 &= \lambda_1 - (\mu_1 - a_1 \Bar{X}_1) (\Bar{X}_1 \wedge C) \\
    \dot{\Bar{X}}_2 &= \lambda_2 - (\mu_2 - a_2 \Bar{X}_2) ((C - \Bar{X}_1 \wedge C)\wedge \Bar{X}_2). \label{equation:two_class_fluid_2}
\end{align}
Note that $\dot{\Bar{X}}_1$ does not depend on $\Bar{X}_2$, therefore, we can independently analyze $\Bar{X}_1$. Then, by definition, an equilibrium point of the system without discontinuity (denoted by $\bar{x}_{e,1}$) is also a stationary point satisfying 
\[
  \lambda_1 - (\mu_1 - a_1 \bar{x}_{e,1}) (\bar{x}_{e,1} \wedge C) = 0.
\]

We consider two cases:
\begin{enumerate}
    \item $\bar{x}_{e,1} > C$: We can obtain the following equilibrium point.
    \begin{align*}
         \lambda_1 - (\mu_1-a_1 \bar{x}_{e,1}) C = 0 \implies \bar{x}_{e,1} = \frac{\mu_1 C - \lambda_1}{a_1 C}. 
    \end{align*}
    In order for this equilibrium point to be valid, we check whether it satisfies the case condition $(\bar{x}_{e,1} > C)$ or not: 
    \begin{align*}
        \bar{x}_{e,1} = \frac{\mu_1 C - \lambda_1}{a_1 C} > \frac{\mu_1 C - \mu_1 C + a_1C}{a_1C} = C.
    \end{align*}
    The inequality follows from the base assumption $\frac{\lambda_1}{\mu_1-a_1 C}<C$. We can conclude that $\bar{x}_{e,1} = \frac{\mu_1 C - \lambda_1}{a_1 C}$ is a valid equilibrium point.

    \item $\bar{x}_{e,1} \leq C$: We can obtain the following equilibrium points.
    \begin{align*}
        \lambda_1 - (\mu_1-a_1 \bar{x}_{e,1}) \bar{x}_{e,1} = 0 =  a_1 \bar{x}_{e,1}^2 - \mu_1 \bar{x}_{e,1}  + \lambda_1 \implies \bar{x}_{e,1} =  \frac{\mu_1 \pm \sqrt{\mu_1^2 - 4 a_1 \lambda_1}}{2 a_1}. 
    \end{align*}
     We check whether the equilibrium points satisfy the case condition $(\bar{x}_{e,1} \leq C)$ or not. First, consider $\frac{\mu_1 + \sqrt{\mu_1^2 - 4 a_1 \lambda_1}}{2 a_1}$. Based on the assumption $\frac{\lambda_1}{\mu_1-a_1 C}<C$, we can write the following:
     \begin{align*}
         \frac{\mu_1 + \sqrt{\mu_1^2 - 4 a_1 \lambda_1}}{2 a_1} &\geq \frac{\mu_1 + \sqrt{\mu_1^2 - 4 a_1 C (\mu_1 - a_1 C)}}{2 a_1} = \frac{\mu_1 + \sqrt{\mu_1^2 - 4 a_1 C \mu_1 +  4 a_1^2 C^2}}{2 a_1} = \frac{\mu_1 + \sqrt{(\mu_1-2a_1 C)^2}}{2 a_1}\\
         & = \frac{\mu_1 + |\mu_1-2a_1 C|}{2 a_1}.
     \end{align*}
     \begin{enumerate}
         \item If $\mu_1 \geq 2 a_1 C$:
         \begin{align*}
             \frac{\mu + \sqrt{\mu_1^2 - 4 a_1 \lambda_1}}{2 a_1} \geq \frac{2\mu_1-2a_1C}{2 a_1} \geq  \frac{2a_1 C-a_1 C}{ a_1} = C.  
         \end{align*}
         \item If $\mu_1 <2 a_1 C$:
         \begin{align*}
             \frac{\mu_1 + \sqrt{\mu_1^2 - 4 a_1 \lambda_1}}{2 a_1} \geq \frac{2a_1 C}{2 a_1} =  C.  
         \end{align*}
     \end{enumerate}
     The inequalities follow from the base assumption $\frac{\lambda_1}{\mu_1-a_1 C}<C$. We can conclude that $\bar{x}_{e,1} = \frac{\mu_1 + \sqrt{\mu_1^2 - 4 a_1 \lambda_1}}{2 a_1}$ is not a valid equilibrium point. Next, we consider $\frac{\mu_1 - \sqrt{\mu_1^2 - 4 a_1 \lambda_1}}{2 a_1}$. Based on the assumption $\frac{\lambda_1}{\mu_1-a_1 C}<C$, we can write the following:
     \begin{align*}
         \frac{\mu_1 - \sqrt{\mu_1^2 - 4 a_1 \lambda_1}}{2 a_1} &\leq \frac{\mu_1 - \sqrt{\mu_1^2 - 4 a_1 C (\mu_1 - a_1 C)}}{2 a_1} = \frac{\mu_1 - \sqrt{\mu_1^2 - 4 a_1 C \mu_1 +  4 a_1^2 C^2}}{2 a_1} = \frac{\mu_1 - \sqrt{(\mu_1-2a_1C)^2}}{2 a_1}\\
         & = \frac{\mu_1 - |\mu_1-2a_1 C|}{2 a_1}.
     \end{align*}
     \begin{enumerate}
         \item If $\mu_1 \geq 2 a_1 C$:
         \begin{align*}
             \frac{\mu_1 - \sqrt{\mu_1^2 - 4 a_1 \lambda_1}}{2 a_1} \leq \frac{2a_1 C}{2 a_1} =  C.  
         \end{align*}
         \item If $\mu_1 <2 a_1 C$:
         \begin{align*}
             \frac{\mu_1 - \sqrt{\mu_1^2 - 4 a_1 \lambda_1}}{2 a_1} \leq \frac{2\mu_1-2a_1 C}{2 a_1} <  \frac{2a_1 C-a_1 C}{ a_1} = C.  
         \end{align*}
     The inequalities follow from the base assumption $\frac{\lambda}{\mu_1-a_1 C}<C$. This proves that $\bar{x}_{e,1} = \frac{\mu_1 - \sqrt{\mu_1^2 - 4 a_1 \lambda_1}}{2 a_1}$ is a valid equilibrium point. 
     \end{enumerate}
\end{enumerate}

We can conclude that the system without truncation has two equilibrium points:
\begin{align*}
    \bar{x}_{e,1} \in \left\{\frac{\mu_1 - \sqrt{\mu_1^2 - 4 a_1 \lambda_1}}{2 a_1},\frac{\mu_1 C - \lambda_1}{a_1 C}\right\}.
\end{align*}

Recall that we denote $\dot{\Bar{X}}_1 = \varphi_1(\Bar{X}_1)$. Now, we examine their stability one by one. First consider $\bar{x}_{e,1}^{(1)} = \frac{\mu - \sqrt{\mu_1^2 - 4 a_1 \lambda_1}}{2 a_1}$.  Let $\sqrt{\mu_1^2 - 4 a_1 \lambda_1}/a_1\geq \epsilon>0$, then
\begin{align*}
    \varphi_1(\bar{x}_{e,1}^{(1)} + \epsilon)&= \lambda_1 - \left(\mu_1 - a_1 \left(\frac{\mu_1 - \sqrt{\mu^2_1 - 4 a_1 \lambda_1}}{2 a_1} + \epsilon \right)\right)\left(\frac{\mu_1 - \sqrt{\mu_1^2 - 4 a_1 \lambda}}{2 a_1} + \epsilon\right) \\
    &= \lambda_1 - \mu_1 \left(\frac{\mu_1 - \sqrt{\mu^2_1 - 4 a_1 \lambda_1}}{2 a_1}  + \epsilon\right) + a_1 \left(\frac{\mu_1 - \sqrt{\mu^2_1 - 4 a_1 \lambda_1}}{2 a_1}  +\epsilon\right)^2 \\
    &= \lambda_1 - \mu_1 \left(\frac{\mu_1 - \sqrt{\mu^2_1 - 4 a_1 \lambda_1}}{2 a_1} \right) - \mu_1 \epsilon + a_1 \left(\frac{\mu_1 - \sqrt{\mu^2_1 - 4 a_1 \lambda_1}}{2 a_1} \right)^2 + 2a_1 \epsilon \left(\frac{\mu_1 - \sqrt{\mu^2_1 - 4 a_1 \lambda_1}}{2 a_1} \right) + a_1 \epsilon^2\\
    &= -\mu_1\epsilon + 2a_1 \epsilon \left(\frac{\mu_1 - \sqrt{\mu^2_1 - 4 a_1 \lambda_1}}{2 a_1} \right) + a_1 \epsilon^2 = - \epsilon \sqrt{\mu^2_1 - 4 a_1 \lambda_1} + a_1 \epsilon^2 = \epsilon(a_1 \epsilon - \sqrt{\mu^2_1 - 4 a_1 \lambda_1}) \\
    &\leq \sqrt{\mu^2_1 - 4 a_1 \lambda}/a_1 (\sqrt{\mu^2_1 - 4 a_1 \lambda_1} -\sqrt{\mu^2_1 - 4 a_1 \lambda_1})=0
\end{align*}
and, 
\begin{align*}
    \varphi_1(\bar{x}_{e,1}^{(1)} - \epsilon) &= \lambda_1 - \left(\mu_1 - a_1 \left(\frac{\mu_1 - \sqrt{\mu^2_1 - 4 a_1 \lambda_1}}{2 a_1} - \epsilon \right)\right)\left(\frac{\mu_1 - \sqrt{\mu^2_1 - 4 a_1 \lambda_1}}{2 a_1} - \epsilon\right) \\
    &= \lambda_1 - \mu_1 \left(\frac{\mu_1 - \sqrt{\mu^2_1 - 4 a_1 \lambda_1}}{2 a_1}  - \epsilon\right) + a \left(\frac{\mu_1 - \sqrt{\mu^2_1 - 4 a_1 \lambda_1}}{2 a_1}  -\epsilon\right)^2 \\
    &= \lambda_1 - \mu_1 \left(\frac{\mu_1 - \sqrt{\mu^2_1 - 4 a \lambda_1}}{2 a} \right) + \mu_1 \epsilon + a_1 \left(\frac{\mu_1 - \sqrt{\mu^2_1 - 4 a_1 \lambda_1}}{2 a_1} \right)^2 - 2a_1 \epsilon \left(\frac{\mu_1 - \sqrt{\mu^2_1 - 4 a_1 \lambda_1}}{2 a_1} \right) + a_1 \epsilon^2\\
    &=\mu_1 \epsilon - 2a_1 \epsilon \left(\frac{\mu_1 - \sqrt{\mu^2_1 - 4 a_1 \lambda_1}}{2 a_1} \right) + a_1 \epsilon^2 = \epsilon \sqrt{\mu^2_1 - 4 a_1 \lambda_1} + a_1 \epsilon^2 >0.
\end{align*}

This implies that if $||\Bar{X}_1(0) - \bar{x}_{e,1}^{(1)}|| < \epsilon$ then $||\Phi^{\bar{\pi}}(\Bar{X}_1(0), t) - \bar{x}_{e,1}^{(1)}|| < ||\Bar{X}_1(0) -\bar{x}_{e,1}^{(1)}|| $ and hence $\bar{x}_{e,1}^{(1)} = \frac{\mu_1 - \sqrt{\mu^2_1 - 4 a_1 \lambda_1}}{2 a_1}$ is a stable equilibrium point.

Secondly consider $\bar{x}_{e,1}^{(2)} =\frac{\mu_1 C - \lambda_1}{a_1C}$.  Let $\epsilon>0$, then
\begin{align*}
    \varphi(\bar{x}_{e,1}^{(2)} + \epsilon) &= \lambda_1 - \left(\mu_1 - a_1 \left(\frac{\mu_1 C - \lambda_1}{a_1C} + \epsilon \right)\right)C = \lambda_1 - \mu_1 C + a_1 C \left(\frac{\mu_1 C - \lambda_1}{a_1 C} + \epsilon \right) = a_1 C\epsilon >0
\end{align*}
and,
\begin{align*}
    \varphi(\bar{x}_{e,1}^{(2)} - \epsilon) &= \lambda_1 - \left(\mu_1 - a_1 \left(\frac{\mu_1 C - \lambda_1}{a_1 C} - \epsilon \right)\right)C = \lambda_1 - \mu_1 C + a_1 C \left(\frac{\mu_1 C - \lambda_1}{a_1 C} - \epsilon \right) = -a_1 C\epsilon <0.
\end{align*}
This implies that if $||\Bar{X}_1(0) - \bar{x}_{e,1}^{(1)}|| < \epsilon \rightarrow ||\Phi_1^{\bar{\pi}}(\Bar{X}_1(0), t) - \bar{x}_{e,1}^{(1)}|| > ||\Bar{X}_1(0) - \bar{x}_{e,1}^{(1)}|| $ and $\bar{x}_{e,1}^{(2)} = \frac{\mu_1 C - \lambda_1}{a_1 C}$ is not a stable equilibrium point.

Now, we consider the system with truncation. In addition to the equilibrium points identified in the system without truncation, the truncated system also has the truncation threshold $\kappa_1$. Note that our system has the base assumptions $\lambda_1/(\mu_1 - a_1 C) \leq C$, $\mu_1 \geq a_1 \kappa_1$ and $\kappa_1>C$. Then, by Lemma \ref{lemma:kappa_stab}, $\kappa_1$ is characterized as a stable pseudo equilibrium point. Thus, $\Bar{X}_1$ has one stable equilibrium $\left(\frac{\mu_1 - \sqrt{\mu^2_1 - 4 a_1 \lambda_1}}{2 a_1}\right)$ and one stable pseudo equilibrium $\left(\kappa_1\right)$.

Note that in order system to be in equilibrium both $\Bar{X}_1$ and $\Bar{X}_2$ should be in their respective equilibria. Therefore, next, we examine the behavior of $\Bar{X}_2$ for each $\bar{x}_{e,1} \in \left\{\frac{\mu_1 - \sqrt{\mu^2_1 - 4 a_1 \lambda_1}}{2 a_1}, \kappa_1 \right\}$.

We consider two cases:
\begin{enumerate}
    \item $\bar{x}_{e,1} = \frac{\mu_1 - \sqrt{\mu^2_1 - 4 a_1 \lambda_1}}{2 a_1} \leq C$: Under this case, the stationary points of $\Bar{X}_2$ (denoted by $\bar{x}_{e,2}$) satisfies
    \begin{align}
        \lambda_2 - (\mu_2-a_2 \bar{x}_{e,2})((C-\bar{x}_{e,1})\wedge x_2) = 0.
    \end{align}
    We conduct a similar analysis that we have conducted for $\Bar{X}_1$ to $\Bar{X}_2$ by only replacing $C$ with $C-\bar{x}_{e,1}$. As a result of this analysis, we identify one stable equilibrium and one stable pseudo-equilibrium:

    \begin{align*}
        \bar{x}_{e,2} \in \left\{\frac{\mu_2 \pm \sqrt{\mu_2^2 - 4 a_2 \lambda_2}}{2a_2}, \kappa_2 \right\}.
    \end{align*}

    \item $\bar{x}_{e,1} = \kappa_1$: Under this case, 
    \begin{align}
        \dot{\Bar{X}}_2 = \lambda_2 \mathbbm{1}\{\Bar{X}_2 < \kappa_2\},
    \end{align}
    and therefore, there is no stationary point and only one stable pseudo equilibrium that appears at the point of discontinuity $\bar{x}_{e,2}=\kappa_2$.
\end{enumerate}
Thus, the system defined by the ODE \eqref{equation:two_class_fluid_1}-\eqref{equation:two_class_fluid_2} yields one stable equilibrium and two stable pseudo equilibrium points:
\begin{align*}
    \bar{x}_e \in 
    \Bigg\{ \left(\frac{\mu_1 - \sqrt{\mu_1^2 - 4 a_1 \lambda_1}}{2a_1},\frac{\mu_2 - \sqrt{\mu_2^2 - 4 a_2 \lambda_2}}{2a_2}\right),\ \left(\frac{\mu_1 - \sqrt{\mu_1^2 - 4 a_1 \lambda_1}}{2a_1},\kappa_2\right), \left(\kappa_1,\kappa_2\right) \Bigg\}.    
\end{align*}

\halmos

\section{An Example of Meta-stability} \label{sec:additional_vector_field_example}
\begin{figure}[!h]
    \centering
    \includegraphics[scale=0.4]{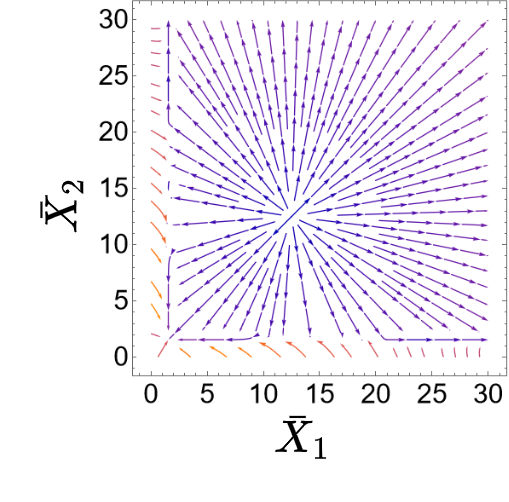}
    \caption{Vector field of the fluid model \eqref{eq:fluid_dynamics} with parameters $\lambda_1=\lambda_2=1.5$, $\mu_1=\mu_2=1$, $a_1=a_2=0.022$, $\kappa_1=\kappa_2=30$, $C=4$ experiencing meta-stability under the policy of assigning equal capacity to each customer class.}
    \label{fig:vector_field_additional}
\end{figure}

We provide an example system that experiences meta-stability with four potential equilibria. The system has $4$ servers with arrival rates $\lambda_1=\lambda_2=1.5$, service rates $\mu_1=\mu_2=1$, slowdown rates $a_1=a_2=0.022$, and blocking thresholds $\kappa_1=\kappa_2=30$. The system operates under a policy that assigns equal capacity to each class of customers. Figure \ref{fig:vector_field_additional} illustrates the vector field for the system.

\section{Fluid-based Policy Iteration}\label{sec:fluid_ADP}
In this section, we introduce a fluid-based policy iteration to find a good initial policy. In particular, we use the fluid model introduced in Section \ref{sec:fluid_approx} to approximate the differences in the relative value function. Note that since the fluid model is a deterministic dynamical system, the approximation does not have any estimation variance, but may incur a large estimation bias. Thus, the policy learned based on this approximation is only used as an initialization.

For a given scheduling policy $\pi$ for the stochastic system, we denote $\bar{\pi}$ as the corresponding fluid policy, which relaxes the integral requirement. Let $\bar X^{\bar{\pi}}(t;\bar{x}_0)$ denote the fluid trajectory under $\bar{\pi}$ starting from $\bar X^{\bar{\pi}}(0)=\bar{x}_0$, i.e.,
\[
\dot{\bar X}_i^{\bar{\pi}}(t)=\lambda_i1\{\bar X_i^{\bar{\pi}}(t)<\kappa_i\} - \Bar{f}_i(\bar X_i^{\bar{\pi}}(t))\bar Z_i^{\bar{\pi}}(t),
\]
where $\bar Z_i^{\bar\pi}(t)$ is the amount of capacity allocated to serve class $i$ fluid at time $t$ under policy $\bar{\pi}$, which maps the fluid trajectory $\bar X^{\bar{\pi}}(t)$ to $\bar Z^{\bar\pi}(t)$. Let $\bar D_i^{\bar{\pi}}(x)$ denote the fluid-based value function difference, which is constructed as follows. If $\bar X^{\bar{\pi}}(t;x)$ and $\bar X^{\bar{\pi}}(t;x-\mathbf{e}_i)$, $i=1,\dots, I$, all converge to the same equilibrium, $\bar D_i^{\bar{\pi}}(x)$ is the cost difference before the two systems reach that equilibrium, i.e.,
\[
\bar D_i^{\bar{\pi}}(x)=\int_{0}^{\infty}c(\bar X^{\bar{\pi}}(t;x))-c(\bar X^{\bar{\pi}}(t;x-\mathbf{e}_i))dt.
\]
Note that once $\bar X^{\bar{\pi}}(t;x)$ and $\bar X^{\bar{\pi}}(t;x-\mathbf{e}_i)$ reach the equilibrium point, $c(\bar X^{\bar{\pi}}(t;x))-c(\bar X^{\bar{\pi}}(t;x-\mathbf{e}_i))=0$.
If they converge to different equilibrium points, let
\[
\mathcal{X}_e^{\bar{\pi}}(x)=\left\{\lim_{t\rightarrow\infty} \bar X^{\bar{\pi}}(t;x-\mathbf{e}_i), i=1,\dots, I\right\},
\]
denote the set of equilibrium points. Let
$\bar{x}_e^{\bar{\pi},*}(x)=\argmin_{x\in \mathcal{X}_e^{\bar{\pi}}(x)} c(x)$, which denotes the equilibrium point with the smallest cost.
We set $\bar D_i^{\bar{\pi}}(x)=1$ if $\lim_{t\rightarrow\infty} \bar X^{\bar{\pi}}(t;x-e_i)=\bar{x}_e^{\bar{\pi},*}(x)$; 
and $\bar D_i^{\bar{\pi}}(x)=0$ otherwise. In this way, we would prefer to move in the direction that reaches the lowest-cost equilibrium point.

We run the policy iteration algorithm with the fluid-based value function approximation scheme introduced above to construct a good initial policy. The algorithm is summarized in Algorithm \ref{algo:initial}.

\begin{algorithm}
  \caption{Fluid-based Policy Iteration} \label{algo:initial}
  Start with an arbitrary policy $\bar{\pi}^{0}(\cdot)$, set $n = 0$ and input $n_{max}$\;
  Uniformly sample $N$ points from the state-space $\mathcal{X}$ and denote the sampled set of points as $\mathcal{X}_0$\;
  \tcc{Policy evaluation:}
  Call Algorithm \ref{algo:fluid_valfunc_estimate} to evaluate $\Bar{D}_i^{\bar{\pi}_{n}}(x)$, $i=1,\dots,I$ for $x \in \mathcal{X}_0$ \;
  \tcc{Policy improvement:}
  Using $\{\bar f_i(x_i) \bar D_i^{\bar{\pi}_n}(x); i=1,\dots, I, x\in \mathcal{X}_0\}$ as the input data, train a classifier that classifies each state into a priority order based on the descending order of $\bar f_i(x_i)\bar D_i^{\bar{\pi}_n}(x)$. Let $\bar{\pi}^{n+1}$ denote the resulting policy\;
  Set $n = n+1$. If $n< \bar n_{max}$, go to Step 3; otherwise, terminate and output $\bar{\pi}^n$.
\end{algorithm}

\begin{algorithm}
    \caption{Fluid-Based Value Function Difference Estimation.}\label{algo:fluid_valfunc_estimate}
    Simulate the fluid trajectories $\bar X^{\bar{\pi}}(t;x)$ and $\bar X^{\bar{\pi}}(t;x-\mathbf{e}_i)$,$i=1,\dots,I$ under policy $\bar{\pi}$\;
    \eIf{$\bar X^{\bar{\pi}_{n}}(t;x)$ and $\bar X^{\bar{\pi}_{n}}(t;x-\mathbf{e}_i),$ $i=1,\dots,I$ all converge to the same equilibrium}{Set 
    \[\bar D_i^{\bar{\pi}_{n}}(x)=\int_{0}^{\infty}c(\bar X^{\bar{\pi}_{n}}(t;x))-c(\bar X^{\bar{\pi}_{n}}(t;x-\mathbf{e}_i))dt, ~i=1,\dots,I;\]}
    {Find the equilibrium point with the smallest cost, i.e., $\bar{x}_e^{\bar{\pi}_{n},*}(x)=\argmin_{x\in \mathcal{X}_e^{\bar{\pi}_{n}}(x)
    } c(x)$, where $\mathcal{X}_e^{\bar{\pi}_{n}}(x)=\left\{\lim\limits_{t\rightarrow\infty} \bar X^{\bar{\pi}_{n}}(t;x-\mathbf{e}_i), i=1,\dots, I\right\}$\;
    \eIf{$\lim_{t\rightarrow\infty} \bar X^{\bar{\pi}_{n}}(t;x-\mathbf{e}_i)=\bar{x}_e^{\bar{\pi}_{n},*}(x)$}
    { $\bar D_i^{\bar{\pi}_{n}}(x)=1$\;}{$\bar D_i^{\bar{\pi}_{n}}(x)=0$\;}}
    Output $\bar D_i^{\bar{\pi}_{n}}(x)$, $i=1,\dots,I$\;
\end{algorithm}

\begin{figure}[!h]
     \centering
     \begin{subfigure}[b]{0.48\textwidth}
         \centering
         \includegraphics[width=.85\textwidth]{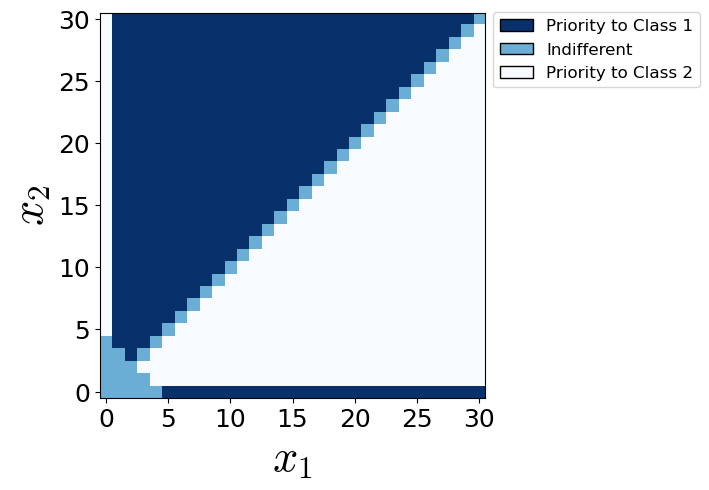}
         \caption{SI: Optimal Policy}
     \end{subfigure}
    \begin{subfigure}[b]{0.48\textwidth}
         \centering
         \includegraphics[width=.85\textwidth]{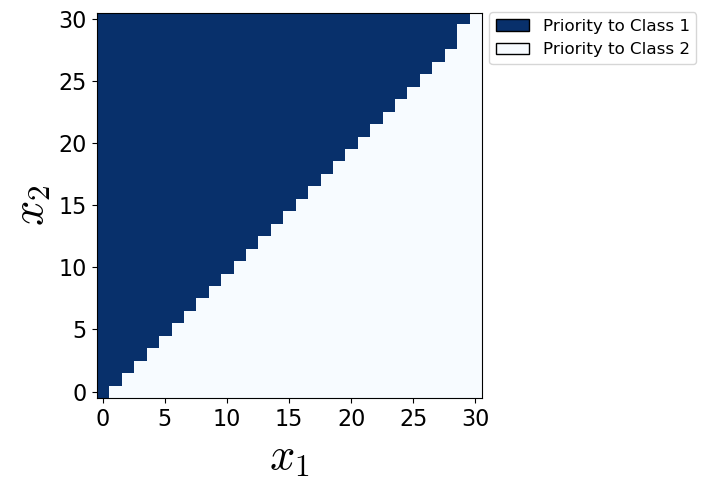}
         \caption{SI: Fluid-Based Policy}
     \end{subfigure}\\
     \begin{subfigure}[b]{0.48\textwidth}
         \centering
         \includegraphics[width=.85\textwidth]{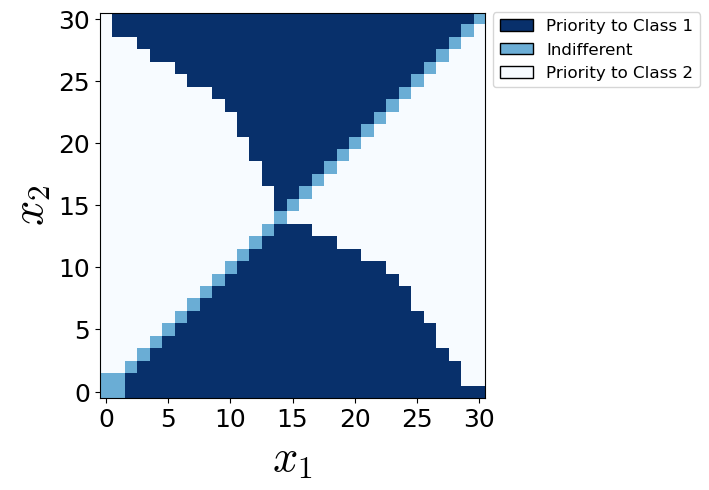}
         \caption{SII: Optimal Policy}
     \end{subfigure}
      \begin{subfigure}[b]{0.48\textwidth}
         \centering
         \includegraphics[width=.85\textwidth]{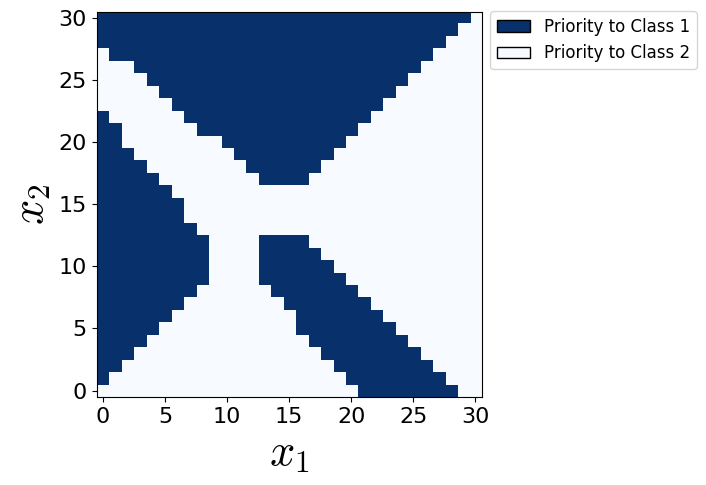}
         \caption{SII: Fluid-Based Policy}
     \end{subfigure}
     \caption{Comparison of the optimal policy versus the policy generated by the fluid-based policy iteration (Fluid-Based Policy). Scenario I (SI):$\lambda = (1.4,1.4), \mu = (1,1), a = (0.01,0.01), \kappa = (30,30), h = (1,1), b = (0,0), C=4$. Scenario II (SII): $\lambda = (0.3,0.3), \mu = (0.9,0.9), a = (0.023,0.023), \kappa = (30,30), h = (1,1), b = (0,0), C=1$.}
     \label{fig:fluid_example}
     \end{figure}

Figure \ref{fig:fluid_example} provides two examples for two-class queues with linear slowdown, i.e., $\bar{f}_i(x)= \mu_i - a_i x_i$. We compare the optimal policy obtained by solving the MDP with the policy obtained by the fluid-based policy iteration (fluid-based policy), i.e., Algorithm \ref{algo:initial}. We observe that in the first scenario, the fluid-based policy is identical to the optimal policy ((a) versus (b)). In the second scenario, the fluid-based policy looks quite different from the optimal policy ((c) versus (d)). However, we observe empirically that it provides a good initial policy that can accelerate the convergence of the ADP algorithm. 

\section{Implementation details}\label{sec:implementation_details}
\begin{figure}[!h]
         \centering
         \includegraphics[scale=0.5]{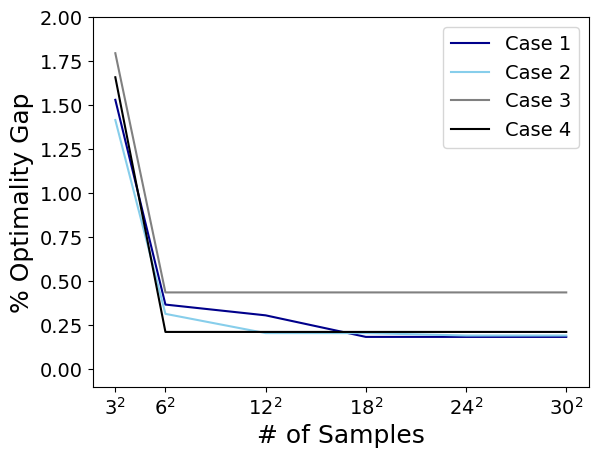}
     \caption{ADP policies generated using different numbers of samples for 4 different parameter regimes.} 

    \raggedright
     \fignote{\textbf{Case 1}: $\lambda=(0.5,2.5),$ $\mu=(1,1),$ $a = (0.0167,0.0167),$ $h=(1.5,1),$ $b=(0,0),$ $\kappa = (30,30),$ $C=4$. \textbf{Case 2}: $\lambda=(1.5,1.5),$ $\mu=(1,1),$ $a = (0.0103,0.0203),$ $h=(3.5,1),$ $b=(0,0),$ $\kappa = (30,30),$ $C=4$. \textbf{Case 3}: $\lambda=(1.5,1.5),$ $\mu=(1,1),$ $a = (0.0103,0.0203),$ $h=(1.5,1),$ $b=(0,0),$ $\kappa = (30,30),$ $C=4$. \textbf{Case 4}: $\lambda=(1.5,1.5),$ $\mu=(1,1),$ $a = (0.0103,0.0203),$ $h=(5,1),$ $b=(0,1000),$ $\kappa = (30,30),$ $C=4$.}
     \label{fig:scaling_experiments}
     \end{figure}

We first examine the effect of $N$ in Algorithm \ref{algo:main}, i.e., the number of uniformly sampled states where the value function differences are evaluated through simulation. We shall refer to $N$ as the sample size since it is the sample size used to train the classifier.
To measure the optimality gap, we focus on a two-class system with linearly decreasing service rate function, i.e., $\bar{f}_i (x_i) =\mu_i - a_i x_i$, $i=1,2$, $C = 4$ servers, and blocking threshold $\kappa = (30,30)$. The size of the state space is $31^2=961$. 
Figure \ref{fig:scaling_experiments} plots the percentage optimality gap as a function of the sample size $N$ under four different parameter settings (see the caption of Figure \ref{fig:scaling_experiments} for details of the model parameters). We observe that even with $N=9$, our method can learn a policy with less than $2\%$ optimality gap. Moreover, a sample size of $36$ is already enough to generate approximately optimal policies, i.e., the optimality gap is less than $0.5\%$. This illustrates the scalability of the proposed ADP approach. In practice, if computationally viable, we recommend setting $N$ at around 5\% of the size of the state space. Otherwise, for larger systems, it is still possible to obtain well-performing policies with a much smaller sample size. For example, in Section \ref{sec:fiveclass}, when we consider a five-class system, the size of the state space is $31^5=2.86\times 10^7$, we only sample $0.004\%$ of the state space, but still achieve good performance. 

Next, we examine the effect of $n_{max}$, i.e., the number of policy iterations to perform. Typically, we observe that if $N$ is selected sufficiently large with respect to the cardinality of the state space (e.g., around 5\% of the size of the state space), our algorithm is able to find a close-to-optimal solution within $5$ iterations. However, if $N$ is very small, less than 0.1\% of the size of the state space, there can be large fluctuations between iterations and we need to run more iterations, 10 to 20 iterations. For example, in Section \ref{sec:fiveclass}, when N is $0.004\%$ of the size of the state space, a well-performing policy can be achieved within $15$ iterations.

For Algorithm \ref{algo:coupling}, we estimate the value function differences by simulating carefully designed coupled systems. Even though we can find $\tau(x)$, i.e., the time the coupled systems coincide, relatively quickly in most cases, there can be some rare cases where the coupled systems do not coincide for a very long time. To avoid excessive simulation time, we introduce a truncation time $T$, when the simulation is terminated if sample paths do not coincide before then. Note that $T$ may introduce some bias in our estimation. Thus, it needs to be chosen to balance the bias and simulation time tradeoff. We suggest simulating $\tau(x)$ and setting $T$ large enough such that it is rarely reached. For example, by looking at the plot in Figure \ref{fig:sim_time_dist}, setting $T\geq 100$ would be reasonable.

For Algorithm \ref{algo:adaptive_sampling}, we need to specify $n_{step}$, $\tilde{n}_{max}$, and $\alpha$. These parameters control the number of simulation replications in value function difference estimation.
First, since our adaptive sampling utilizes Gaussian approximation (based on the Central Limit Theorem) to construct the t-statistics, we suggest setting $n_{step}\geq 30$. Second, setting $\alpha=0.95$ or $0.99$ would be reasonable. Lastly, we would, in general, expect Algorithm \ref{algo:adaptive_sampling} to terminate when $\mathcal{M}^k$ is empty. $\tilde{n}_{max}$ is introduced to avoid excessive long simulations in rare cases. We suggest doing some trial and error and setting $\tilde{n}_{max}$ large enough such that it is rarely reached.

\section{Detailed experiment results}
In this section, we illustrate the detailed experiment results summarized in Sections \ref{sec:worst-case_load} and \ref{sec:service_rate_heterogeneity}. Table \ref{tab:slowdown_effect_experiments} illustrates experiment results for varying worst-case load and holding costs, and Table \ref{tab:service_effect_experiments} illustrates the experiment results for varying service rate differences.

\begin{table}[!h]
\centering
\caption{Long-run average costs under the optimal policy, ADP policy and five benchmark policies and varying service and slowdown rates.}
\label{tab:slowdown_effect_experiments}
\scalebox{0.9}{
\begin{tabular}{lll|ccccccc}
\hline
\textbf{\begin{tabular}[c]{@{}l@{}}Worst-case\\ Load\end{tabular}} & \textbf{$a$}                               & \textbf{$h$}      & \textbf{Optimal}                                           & \textbf{ADP}                                               & \textbf{$h \bar f(0)$}                                        & \textbf{$h \bar f(x)$}                                        & \textbf{\begin{tabular}[c]{@{}c@{}}Max\\ Pressure\end{tabular}} & \textbf{SQF}                                               & \textbf{LQF}                                                \\ \hline
\multirow{11}{*}{\textbf{1}}                                        & \multirow{11}{*}{\textbf{(0.0025, 0.0125)}} & \textbf{(1, 1)}   & \begin{tabular}[c]{@{}c@{}}5.5\\ $\pm$ 0.01\end{tabular}   & \begin{tabular}[c]{@{}c@{}}5.55\\ $\pm$ 0.01\end{tabular}  & \begin{tabular}[c]{@{}c@{}}5.5\\ $\pm$ 0.01\end{tabular}   & \begin{tabular}[c]{@{}c@{}}6.09\\ $\pm$ 0.02\end{tabular}  & \begin{tabular}[c]{@{}c@{}}5.85\\ $\pm$ 0.02\end{tabular}       & \begin{tabular}[c]{@{}c@{}}5.76\\ $\pm$ 0.02\end{tabular}  & \begin{tabular}[c]{@{}c@{}}5.8\\ $\pm$ 0.02\end{tabular}    \\
                                                                   &                                            & \textbf{(1, 1.5)} & \begin{tabular}[c]{@{}c@{}}6.9\\ $\pm$ 0.02\end{tabular}   & \begin{tabular}[c]{@{}c@{}}6.92\\ $\pm$ 0.02\end{tabular}  & \begin{tabular}[c]{@{}c@{}}6.95\\ $\pm$ 0.02\end{tabular}  & \begin{tabular}[c]{@{}c@{}}6.95\\ $\pm$ 0.02\end{tabular}  & \begin{tabular}[c]{@{}c@{}}7.2\\ $\pm$ 0.03\end{tabular}        & \begin{tabular}[c]{@{}c@{}}7.18\\ $\pm$ 0.02\end{tabular}  & \begin{tabular}[c]{@{}c@{}}7.28\\ $\pm$ 0.03\end{tabular}   \\
                                                                   &                                            & \textbf{(1, 2)}   & \begin{tabular}[c]{@{}c@{}}7.82\\ $\pm$ 0.02\end{tabular}  & \begin{tabular}[c]{@{}c@{}}7.82\\ $\pm$ 0.02\end{tabular}  & \begin{tabular}[c]{@{}c@{}}7.82\\ $\pm$ 0.02\end{tabular}  & \begin{tabular}[c]{@{}c@{}}7.82\\ $\pm$ 0.02\end{tabular}  & \begin{tabular}[c]{@{}c@{}}8.3\\ $\pm$ 0.03\end{tabular}        & \begin{tabular}[c]{@{}c@{}}8.6\\ $\pm$ 0.02\end{tabular}   & \begin{tabular}[c]{@{}c@{}}8.76\\ $\pm$ 0.03\end{tabular}   \\
                                                                   &                                            & \textbf{(1, 2.5)} & \begin{tabular}[c]{@{}c@{}}8.68\\ $\pm$ 0.02\end{tabular}  & \begin{tabular}[c]{@{}c@{}}8.68\\ $\pm$ 0.02\end{tabular}  & \begin{tabular}[c]{@{}c@{}}8.68\\ $\pm$ 0.02\end{tabular}  & \begin{tabular}[c]{@{}c@{}}8.68\\ $\pm$ 0.02\end{tabular}  & \begin{tabular}[c]{@{}c@{}}9.28\\ $\pm$ 0.03\end{tabular}       & \begin{tabular}[c]{@{}c@{}}10.03\\ $\pm$ 0.03\end{tabular} & \begin{tabular}[c]{@{}c@{}}10.25\\ $\pm$ 0.04\end{tabular}  \\
                                                                   &                                            & \textbf{(1, 3)}   & \begin{tabular}[c]{@{}c@{}}9.55\\ $\pm$ 0.02\end{tabular}  & \begin{tabular}[c]{@{}c@{}}9.55\\ $\pm$ 0.02\end{tabular}  & \begin{tabular}[c]{@{}c@{}}9.55\\ $\pm$ 0.02\end{tabular}  & \begin{tabular}[c]{@{}c@{}}9.55\\ $\pm$ 0.02\end{tabular}  & \begin{tabular}[c]{@{}c@{}}10.2\\ $\pm$ 0.03\end{tabular}       & \begin{tabular}[c]{@{}c@{}}11.45\\ $\pm$ 0.03\end{tabular} & \begin{tabular}[c]{@{}c@{}}11.73\\ $\pm$ 0.04\end{tabular}  \\
                                                                   &                                            & \textbf{(1, 3.5)} & \begin{tabular}[c]{@{}c@{}}10.41\\ $\pm$ 0.02\end{tabular} & \begin{tabular}[c]{@{}c@{}}10.41\\ $\pm$ 0.02\end{tabular} & \begin{tabular}[c]{@{}c@{}}10.41\\ $\pm$ 0.02\end{tabular} & \begin{tabular}[c]{@{}c@{}}10.41\\ $\pm$ 0.02\end{tabular} & \begin{tabular}[c]{@{}c@{}}11.1\\ $\pm$ 0.03\end{tabular}       & \begin{tabular}[c]{@{}c@{}}12.87\\ $\pm$ 0.03\end{tabular} & \begin{tabular}[c]{@{}c@{}}13.21\\ $\pm$ 0.05\end{tabular}  \\ \hhline{===|=======}
\multicolumn{3}{l}{\textbf{Avg. Opt. Gap (\%)}}                                                                                     & \begin{tabular}[c]{@{}c@{}}8.14\\ $\pm$ 0.04\end{tabular}  & \begin{tabular}[c]{@{}c@{}}0.21\\ $\pm$ 0.03\end{tabular}  & \begin{tabular}[c]{@{}c@{}}0.13\\ $\pm$ 0.03\end{tabular}  & \begin{tabular}[c]{@{}c@{}}1.91\\ $\pm$ 0.11\end{tabular}  & \begin{tabular}[c]{@{}c@{}}6.17\\ $\pm$ 0.09\end{tabular}       & \begin{tabular}[c]{@{}c@{}}13.00\\ $\pm$ 0.20\end{tabular}  & \begin{tabular}[c]{@{}c@{}}15.13\\ $\pm$ 0.24\end{tabular}  \\ \hhline{===|=======}
\multirow{11}{*}{\textbf{1.25}}                                     & \multirow{11}{*}{\textbf{(0.0072,0.0172)}}  & \textbf{(1, 1)}   & \begin{tabular}[c]{@{}c@{}}6.05\\ $\pm$ 0.02\end{tabular}  & \begin{tabular}[c]{@{}c@{}}6.05\\ $\pm$ 0.02\end{tabular}  & \begin{tabular}[c]{@{}c@{}}6.05\\ $\pm$ 0.02\end{tabular}  & \begin{tabular}[c]{@{}c@{}}8.16\\ $\pm$ 0.06\end{tabular}  & \begin{tabular}[c]{@{}c@{}}8.28\\ $\pm$ 0.1\end{tabular}        & \begin{tabular}[c]{@{}c@{}}7.05\\ $\pm$ 0.04\end{tabular}  & \begin{tabular}[c]{@{}c@{}}10.27\\ $\pm$ 0.25\end{tabular}  \\
                                                                   &                                            & \textbf{(1, 1.5)} & \begin{tabular}[c]{@{}c@{}}7.91\\ $\pm$ 0.03\end{tabular}  & \begin{tabular}[c]{@{}c@{}}7.95\\ $\pm$ 0.03\end{tabular}  & \begin{tabular}[c]{@{}c@{}}9.07\\ $\pm$ 0.06\end{tabular}  & \begin{tabular}[c]{@{}c@{}}9.07\\ $\pm$ 0.06\end{tabular}  & \begin{tabular}[c]{@{}c@{}}9.8\\ $\pm$ 0.09\end{tabular}        & \begin{tabular}[c]{@{}c@{}}8.6\\ $\pm$ 0.04\end{tabular}   & \begin{tabular}[c]{@{}c@{}}12.87\\ $\pm$ 0.31\end{tabular}  \\
                                                                   &                                            & \textbf{(1, 2)}   & \begin{tabular}[c]{@{}c@{}}9.39\\ $\pm$ 0.04\end{tabular}  & \begin{tabular}[c]{@{}c@{}}9.41\\ $\pm$ 0.04\end{tabular}  & \begin{tabular}[c]{@{}c@{}}9.94\\ $\pm$ 0.06\end{tabular}  & \begin{tabular}[c]{@{}c@{}}9.94\\ $\pm$ 0.06\end{tabular}  & \begin{tabular}[c]{@{}c@{}}10.96\\ $\pm$ 0.09\end{tabular}      & \begin{tabular}[c]{@{}c@{}}10.15\\ $\pm$ 0.05\end{tabular} & \begin{tabular}[c]{@{}c@{}}15.48\\ $\pm$ 0.38\end{tabular}  \\
                                                                   &                                            & \textbf{(1, 2.5)} & \begin{tabular}[c]{@{}c@{}}10.72\\ $\pm$ 0.05\end{tabular} & \begin{tabular}[c]{@{}c@{}}10.75\\ $\pm$ 0.05\end{tabular} & \begin{tabular}[c]{@{}c@{}}10.82\\ $\pm$ 0.06\end{tabular} & \begin{tabular}[c]{@{}c@{}}10.82\\ $\pm$ 0.06\end{tabular} & \begin{tabular}[c]{@{}c@{}}12.15\\ $\pm$ 0.09\end{tabular}      & \begin{tabular}[c]{@{}c@{}}11.7\\ $\pm$ 0.05\end{tabular}  & \begin{tabular}[c]{@{}c@{}}18.09\\ $\pm$ 0.44\end{tabular}  \\
                                                                   &                                            & \textbf{(1, 3)}   & \begin{tabular}[c]{@{}c@{}}11.7\\ $\pm$ 0.06\end{tabular}  & \begin{tabular}[c]{@{}c@{}}11.72\\ $\pm$ 0.06\end{tabular} & \begin{tabular}[c]{@{}c@{}}11.7\\ $\pm$ 0.06\end{tabular}  & \begin{tabular}[c]{@{}c@{}}11.7\\ $\pm$ 0.06\end{tabular}  & \begin{tabular}[c]{@{}c@{}}13.09\\ $\pm$ 0.09\end{tabular}      & \begin{tabular}[c]{@{}c@{}}13.25\\ $\pm$ 0.05\end{tabular} & \begin{tabular}[c]{@{}c@{}}20.69\\ $\pm$ 0.5\end{tabular}   \\
                                                                   &                                            & \textbf{(1, 3.5)} & \begin{tabular}[c]{@{}c@{}}12.58\\ $\pm$ 0.06\end{tabular} & \begin{tabular}[c]{@{}c@{}}12.61\\ $\pm$ 0.06\end{tabular} & \begin{tabular}[c]{@{}c@{}}12.58\\ $\pm$ 0.06\end{tabular} & \begin{tabular}[c]{@{}c@{}}12.58\\ $\pm$ 0.06\end{tabular} & \begin{tabular}[c]{@{}c@{}}13.93\\ $\pm$ 0.08\end{tabular}      & \begin{tabular}[c]{@{}c@{}}14.8\\ $\pm$ 0.06\end{tabular}  & \begin{tabular}[c]{@{}c@{}}23.3\\ $\pm$ 0.57\end{tabular}   \\ \hhline{===|=======}
\multicolumn{3}{l}{\textbf{Avg. Opt. Gap (\%)}}                                                                                     & \begin{tabular}[c]{@{}c@{}}9.73\\ $\pm$ 0.06\end{tabular}  & \begin{tabular}[c]{@{}c@{}}0.29\\ $\pm$ 0.10\end{tabular}  & \begin{tabular}[c]{@{}c@{}}3.59\\ $\pm$ 0.19\end{tabular}  & \begin{tabular}[c]{@{}c@{}}9.37\\ $\pm$ 0.36\end{tabular}  & \begin{tabular}[c]{@{}c@{}}18.78\\ $\pm$ 0.42\end{tabular}      & \begin{tabular}[c]{@{}c@{}}12.29\\ $\pm$ 0.19\end{tabular} & \begin{tabular}[c]{@{}c@{}}70.08\\ $\pm$ 1.57\end{tabular}  \\ \hhline{===|=======}
\multirow{11}{*}{\textbf{1.5}}                                      & \multirow{11}{*}{\textbf{(0.0103,0.0203)}}  & \textbf{(1, 1)}   & \begin{tabular}[c]{@{}c@{}}6.66\\ $\pm$ 0.02\end{tabular}  & \begin{tabular}[c]{@{}c@{}}6.68\\ $\pm$ 0.02\end{tabular}  & \begin{tabular}[c]{@{}c@{}}6.66\\ $\pm$ 0.02\end{tabular}  & \begin{tabular}[c]{@{}c@{}}12.97\\ $\pm$ 0.12\end{tabular} & \begin{tabular}[c]{@{}c@{}}15.76\\ $\pm$ 0.26\end{tabular}      & \begin{tabular}[c]{@{}c@{}}10.02\\ $\pm$ 0.11\end{tabular} & \begin{tabular}[c]{@{}c@{}}41.41\\ $\pm$ 0.77\end{tabular}  \\
                                                                   &                                            & \textbf{(1, 1.5)} & \begin{tabular}[c]{@{}c@{}}8.94\\ $\pm$ 0.04\end{tabular}  & \begin{tabular}[c]{@{}c@{}}9.07\\ $\pm$ 0.04\end{tabular}  & \begin{tabular}[c]{@{}c@{}}14.11\\ $\pm$ 0.13\end{tabular} & \begin{tabular}[c]{@{}c@{}}14.11\\ $\pm$ 0.13\end{tabular} & \begin{tabular}[c]{@{}c@{}}16.59\\ $\pm$ 0.21\end{tabular}      & \begin{tabular}[c]{@{}c@{}}11.65\\ $\pm$ 0.1\end{tabular}  & \begin{tabular}[c]{@{}c@{}}51.81\\ $\pm$ 0.96\end{tabular}  \\
                                                                   &                                            & \textbf{(1, 2)}   & \begin{tabular}[c]{@{}c@{}}10.9\\ $\pm$ 0.05\end{tabular}  & \begin{tabular}[c]{@{}c@{}}11.07\\ $\pm$ 0.06\end{tabular} & \begin{tabular}[c]{@{}c@{}}15.0\\ $\pm$ 0.13\end{tabular}  & \begin{tabular}[c]{@{}c@{}}15.0\\ $\pm$ 0.13\end{tabular}  & \begin{tabular}[c]{@{}c@{}}17.84\\ $\pm$ 0.2\end{tabular}       & \begin{tabular}[c]{@{}c@{}}13.27\\ $\pm$ 0.1\end{tabular}  & \begin{tabular}[c]{@{}c@{}}62.21\\ $\pm$ 1.16\end{tabular}  \\
                                                                   &                                            & \textbf{(1, 2.5)} & \begin{tabular}[c]{@{}c@{}}12.75\\ $\pm$ 0.07\end{tabular} & \begin{tabular}[c]{@{}c@{}}13.03\\ $\pm$ 0.08\end{tabular} & \begin{tabular}[c]{@{}c@{}}15.89\\ $\pm$ 0.13\end{tabular} & \begin{tabular}[c]{@{}c@{}}15.89\\ $\pm$ 0.13\end{tabular} & \begin{tabular}[c]{@{}c@{}}18.24\\ $\pm$ 0.18\end{tabular}      & \begin{tabular}[c]{@{}c@{}}14.9\\ $\pm$ 0.1\end{tabular}   & \begin{tabular}[c]{@{}c@{}}72.61\\ $\pm$ 1.35\end{tabular}  \\
                                                                   &                                            & \textbf{(1, 3)}   & \begin{tabular}[c]{@{}c@{}}14.53\\ $\pm$ 0.08\end{tabular} & \begin{tabular}[c]{@{}c@{}}14.81\\ $\pm$ 0.09\end{tabular} & \begin{tabular}[c]{@{}c@{}}16.77\\ $\pm$ 0.13\end{tabular} & \begin{tabular}[c]{@{}c@{}}16.77\\ $\pm$ 0.13\end{tabular} & \begin{tabular}[c]{@{}c@{}}19.69\\ $\pm$ 0.18\end{tabular}      & \begin{tabular}[c]{@{}c@{}}16.53\\ $\pm$ 0.1\end{tabular}  & \begin{tabular}[c]{@{}c@{}}83.02\\ $\pm$ 1.54\end{tabular}  \\
                                                                   &                                            & \textbf{(1, 3.5)} & \begin{tabular}[c]{@{}c@{}}16.31\\ $\pm$ 0.1\end{tabular}  & \begin{tabular}[c]{@{}c@{}}16.42\\ $\pm$ 0.11\end{tabular} & \begin{tabular}[c]{@{}c@{}}17.66\\ $\pm$ 0.13\end{tabular} & \begin{tabular}[c]{@{}c@{}}17.66\\ $\pm$ 0.13\end{tabular} & \begin{tabular}[c]{@{}c@{}}19.99\\ $\pm$ 0.16\end{tabular}      & \begin{tabular}[c]{@{}c@{}}18.15\\ $\pm$ 0.1\end{tabular}  & \begin{tabular}[c]{@{}c@{}}93.42\\ $\pm$ 1.73\end{tabular}  \\ \hhline{===|=======}
\multicolumn{3}{l}{\textbf{Avg. Opt. Gap (\%)}}                                                                                     & \begin{tabular}[c]{@{}c@{}}11.68\\ $\pm$ 0.09\end{tabular} & \begin{tabular}[c]{@{}c@{}}1.41\\ $\pm$ 0.13\end{tabular}  & \begin{tabular}[c]{@{}c@{}}24.05\\ $\pm$ 0.59\end{tabular} & \begin{tabular}[c]{@{}c@{}}39.79\\ $\pm$ 0.86\end{tabular} & \begin{tabular}[c]{@{}c@{}}64.23\\ $\pm$ 1.24\end{tabular}      & \begin{tabular}[c]{@{}c@{}}24.17\\ $\pm$ 0.49\end{tabular} & \begin{tabular}[c]{@{}c@{}}479.55\\ $\pm$ 4.29\end{tabular} \\ \hhline{===|=======}
\end{tabular}

}

\raggedright
\tabnote{Other system parameters are fixed at $\lambda = (1.5,1.5)$, $\mu=(1,1)$, $b=(0,0)$, $C=4$, $\kappa=(30,30)$.}
\end{table}

\begin{table}[!h]
\centering
\caption{Long-run average costs under the optimal policy, ADP policy and five benchmark policies and varying service and slowdown rates. }
\label{tab:service_effect_experiments}

\scalebox{0.85}{\begin{tabular}{lll|ccccccc}
\hline
$\mu$                                    & $a$                                        & $h$                & \textbf{Optimal}                                            & \textbf{ADP}                                                & $h \bar{f}(0)$                                              & $h \bar{f}(x)$                                              & \textbf{\begin{tabular}[c]{@{}c@{}}Max\\ Pressure\end{tabular}} & \textbf{SQF}                                                & \textbf{LQF}                                                 \\ \hline
\multirow{11}{*}{\textbf{(0.975, 1.025)}} & \multirow{11}{*}{\textbf{(0.0107, 0.0207)}} & \textbf{(1, 1)}    & \begin{tabular}[c]{@{}c@{}}6.06 \\ $\pm 0.02$\end{tabular}  & \begin{tabular}[c]{@{}c@{}}6.09 \\ $\pm 0.02$\end{tabular}  & \begin{tabular}[c]{@{}c@{}}6.06 \\ $\pm 0.02$\end{tabular}  & \begin{tabular}[c]{@{}c@{}}9.35 \\ $\pm 0.09$\end{tabular}  & \begin{tabular}[c]{@{}c@{}}14.57 \\ $\pm$ 0.26\end{tabular}     & \begin{tabular}[c]{@{}c@{}}8.59 \\ $\pm 0.08$\end{tabular}  & \begin{tabular}[c]{@{}c@{}}38.71 \\ $\pm 0.69$\end{tabular}  \\
                                         &                                            & \textbf{(1.25, 1)} & \begin{tabular}[c]{@{}c@{}}7.16 \\ $\pm 0.03$\end{tabular}  & \begin{tabular}[c]{@{}c@{}}7.17 \\ $\pm 0.03$\end{tabular}  & \begin{tabular}[c]{@{}c@{}}11.87 \\ $\pm 0.11$\end{tabular} & \begin{tabular}[c]{@{}c@{}}11.87 \\ $\pm 0.11$\end{tabular} & \begin{tabular}[c]{@{}c@{}}14.99 \\ $\pm$ 0.22\end{tabular}     & \begin{tabular}[c]{@{}c@{}}9.38 \\ $\pm 0.08$\end{tabular}  & \begin{tabular}[c]{@{}c@{}}43.58 \\ $\pm 0.78$\end{tabular}  \\
                                         &                                            & \textbf{(1.5, 1)}  & \begin{tabular}[c]{@{}c@{}}8.18 \\ $\pm 0.04$\end{tabular}  & \begin{tabular}[c]{@{}c@{}}8.24 \\ $\pm 0.03$\end{tabular}  & \begin{tabular}[c]{@{}c@{}}12.28 \\ $\pm 0.11$\end{tabular} & \begin{tabular}[c]{@{}c@{}}12.28 \\ $\pm 0.11$\end{tabular} & \begin{tabular}[c]{@{}c@{}}15.33 \\ $\pm$ 0.20\end{tabular}     & \begin{tabular}[c]{@{}c@{}}10.16 \\ $\pm 0.08$\end{tabular} & \begin{tabular}[c]{@{}c@{}}48.44 \\ $\pm 0.87$\end{tabular}  \\
                                         &                                            & \textbf{(1.75, 1)} & \begin{tabular}[c]{@{}c@{}}9.12 \\ $\pm 0.04$\end{tabular}  & \begin{tabular}[c]{@{}c@{}}9.32 \\ $\pm 0.05$\end{tabular}  & \begin{tabular}[c]{@{}c@{}}12.69 \\ $\pm 0.11$\end{tabular} & \begin{tabular}[c]{@{}c@{}}12.69 \\ $\pm 0.11$\end{tabular} & \begin{tabular}[c]{@{}c@{}}15.17 \\ $\pm$ 0.18\end{tabular}     & \begin{tabular}[c]{@{}c@{}}10.94 \\ $\pm 0.08$\end{tabular} & \begin{tabular}[c]{@{}c@{}}53.31 \\ $\pm 0.95$\end{tabular}  \\
                                         &                                            & \textbf{(2, 1)}    & \begin{tabular}[c]{@{}c@{}}9.99 \\ $\pm 0.05$\end{tabular}  & \begin{tabular}[c]{@{}c@{}}10.03 \\ $\pm 0.05$\end{tabular} & \begin{tabular}[c]{@{}c@{}}13.11 \\ $\pm 0.11$\end{tabular} & \begin{tabular}[c]{@{}c@{}}13.11 \\ $\pm 0.11$\end{tabular} & \begin{tabular}[c]{@{}c@{}}15.53 \\ $\pm$ 0.17\end{tabular}     & \begin{tabular}[c]{@{}c@{}}11.73 \\ $\pm 0.08$\end{tabular} & \begin{tabular}[c]{@{}c@{}}58.18 \\ $\pm 1.04$\end{tabular}  \\
                                         &                                            & \textbf{(2.25, 1)} & \begin{tabular}[c]{@{}c@{}}10.90 \\ $\pm 0.06$\end{tabular} & \begin{tabular}[c]{@{}c@{}}11.08 \\ $\pm 0.06$\end{tabular} & \begin{tabular}[c]{@{}c@{}}13.52 \\ $\pm 0.11$\end{tabular} & \begin{tabular}[c]{@{}c@{}}13.52 \\ $\pm 0.11$\end{tabular} & \begin{tabular}[c]{@{}c@{}}15.94 \\ $\pm$ 0.16\end{tabular}     & \begin{tabular}[c]{@{}c@{}}12.51 \\ $\pm 0.08$\end{tabular} & \begin{tabular}[c]{@{}c@{}}63.05 \\ $\pm 1.12$\end{tabular}  \\ \hhline{===|=======}
\multicolumn{3}{l|}{\textbf{Avg. Opt. Gap (\%)}}                                                           & -                                                           & \begin{tabular}[c]{@{}c@{}}0.96 \\ $\pm$ 0.09\end{tabular}  & \begin{tabular}[c]{@{}c@{}}35.00 \\ $\pm$ 0.64\end{tabular} & \begin{tabular}[c]{@{}c@{}}44.05 \\ $\pm$ 0.56\end{tabular} & \begin{tabular}[c]{@{}c@{}}83.62 \\ $\pm$ 1.27\end{tabular}     & \begin{tabular}[c]{@{}c@{}}24.90 \\ $\pm$ 0.42\end{tabular} & \begin{tabular}[c]{@{}c@{}}496.53 \\ $\pm$ 4.28\end{tabular} \\ \hhline{===|=======}
\multirow{11}{*}{\textbf{(0.952, 1.052)}} & \multirow{11}{*}{\textbf{(0.0111, 0.0211)}} & \textbf{(1, 1)}    & \begin{tabular}[c]{@{}c@{}}6.13 \\ $\pm 0.02$\end{tabular}  & \begin{tabular}[c]{@{}c@{}}6.13 \\ $\pm 0.02$\end{tabular}  & \begin{tabular}[c]{@{}c@{}}6.13 \\ $\pm 0.02$\end{tabular}  & \begin{tabular}[c]{@{}c@{}}6.78 \\ $\pm 0.05$\end{tabular}  & \begin{tabular}[c]{@{}c@{}}15.60 \\ $\pm$ 0.30\end{tabular}     & \begin{tabular}[c]{@{}c@{}}8.18 \\ $\pm 0.08$\end{tabular}  & \begin{tabular}[c]{@{}c@{}}39.90 \\ $\pm 0.69$\end{tabular}  \\
                                         &                                            & \textbf{(1.25, 1)} & \begin{tabular}[c]{@{}c@{}}7.26 \\ $\pm 0.03$\end{tabular}  & \begin{tabular}[c]{@{}c@{}}7.28 \\ $\pm 0.03$\end{tabular}  & \begin{tabular}[c]{@{}c@{}}11.35 \\ $\pm 0.09$\end{tabular} & \begin{tabular}[c]{@{}c@{}}11.35 \\ $\pm 0.09$\end{tabular} & \begin{tabular}[c]{@{}c@{}}14.89 \\ $\pm$ 0.23\end{tabular}     & \begin{tabular}[c]{@{}c@{}}9.02 \\ $\pm 0.08$\end{tabular}  & \begin{tabular}[c]{@{}c@{}}44.92 \\ $\pm 0.77$\end{tabular}  \\
                                         &                                            & \textbf{(1.5, 1)}  & \begin{tabular}[c]{@{}c@{}}8.33 \\ $\pm 0.04$\end{tabular}  & \begin{tabular}[c]{@{}c@{}}8.36 \\ $\pm 0.04$\end{tabular}  & \begin{tabular}[c]{@{}c@{}}11.77 \\ $\pm 0.09$\end{tabular} & \begin{tabular}[c]{@{}c@{}}11.77 \\ $\pm 0.09$\end{tabular} & \begin{tabular}[c]{@{}c@{}}14.63 \\ $\pm$ 0.18\end{tabular}     & \begin{tabular}[c]{@{}c@{}}9.86 \\ $\pm 0.08$\end{tabular}  & \begin{tabular}[c]{@{}c@{}}49.95 \\ $\pm 0.86$\end{tabular}  \\
                                         &                                            & \textbf{(1.75, 1)} & \begin{tabular}[c]{@{}c@{}}9.31 \\ $\pm 0.05$\end{tabular}  & \begin{tabular}[c]{@{}c@{}}9.35 \\ $\pm 0.05$\end{tabular}  & \begin{tabular}[c]{@{}c@{}}12.20 \\ $\pm 0.09$\end{tabular} & \begin{tabular}[c]{@{}c@{}}12.20 \\ $\pm 0.09$\end{tabular} & \begin{tabular}[c]{@{}c@{}}15.32 \\ $\pm$ 0.18\end{tabular}     & \begin{tabular}[c]{@{}c@{}}10.70 \\ $\pm 0.08$\end{tabular} & \begin{tabular}[c]{@{}c@{}}54.97 \\ $\pm 0.94$\end{tabular}  \\
                                         &                                            & \textbf{(2, 1)}    & \begin{tabular}[c]{@{}c@{}}10.24 \\ $\pm 0.05$\end{tabular} & \begin{tabular}[c]{@{}c@{}}10.42 \\ $\pm 0.06$\end{tabular} & \begin{tabular}[c]{@{}c@{}}12.62 \\ $\pm 0.09$\end{tabular} & \begin{tabular}[c]{@{}c@{}}12.62 \\ $\pm 0.09$\end{tabular} & \begin{tabular}[c]{@{}c@{}}15.81 \\ $\pm$ 0.17\end{tabular}     & \begin{tabular}[c]{@{}c@{}}11.53 \\ $\pm 0.08$\end{tabular} & \begin{tabular}[c]{@{}c@{}}59.99 \\ $\pm 1.03$\end{tabular}  \\
                                         &                                            & \textbf{(2.25, 1)} & \begin{tabular}[c]{@{}c@{}}11.14 \\ $\pm 0.06$\end{tabular} & \begin{tabular}[c]{@{}c@{}}11.26 \\ $\pm 0.06$\end{tabular} & \begin{tabular}[c]{@{}c@{}}13.05 \\ $\pm 0.09$\end{tabular} & \begin{tabular}[c]{@{}c@{}}13.05 \\ $\pm 0.09$\end{tabular} & \begin{tabular}[c]{@{}c@{}}15.95 \\ $\pm$ 0.16\end{tabular}     & \begin{tabular}[c]{@{}c@{}}12.37 \\ $\pm 0.08$\end{tabular} & \begin{tabular}[c]{@{}c@{}}65.01 \\ $\pm 1.11$\end{tabular}  \\ \hhline{===|=======}
\multicolumn{3}{l|}{\textbf{Avg. Opt. Gap (\%)}}                                                           & -                                                           & \begin{tabular}[c]{@{}c@{}}0.71 \\ $\pm$ 0.10\end{tabular}  & \begin{tabular}[c]{@{}c@{}}28.22 \\ $\pm$ 0.56\end{tabular} & \begin{tabular}[c]{@{}c@{}}30.00 \\ $\pm$ 0.52\end{tabular} & \begin{tabular}[c]{@{}c@{}}82.34 \\ $\pm$ 1.38\end{tabular}     & \begin{tabular}[c]{@{}c@{}}19.19 \\ $\pm$ 0.37\end{tabular} & \begin{tabular}[c]{@{}c@{}}504.25 \\ $\pm$ 4.19\end{tabular} \\ \hhline{===|=======}
\multirow{11}{*}{\textbf{(0.930, 1.080)}} & \multirow{11}{*}{\textbf{(0.0115, 0.0215)}} & \textbf{(1, 1)}    & \begin{tabular}[c]{@{}c@{}}6.26 \\ $\pm 0.03$\end{tabular}  & \begin{tabular}[c]{@{}c@{}}6.26 \\ $\pm 0.03$\end{tabular}  & \begin{tabular}[c]{@{}c@{}}6.26 \\ $\pm 0.03$\end{tabular}  & \begin{tabular}[c]{@{}c@{}}6.33 \\ $\pm 0.03$\end{tabular}  & \begin{tabular}[c]{@{}c@{}}18.23 \\ $\pm$ 0.35\end{tabular}     & \begin{tabular}[c]{@{}c@{}}7.88 \\ $\pm 0.07$\end{tabular}  & \begin{tabular}[c]{@{}c@{}}41.10 \\ $\pm 0.67$\end{tabular}  \\
                                         &                                            & \textbf{(1.25, 1)} & \begin{tabular}[c]{@{}c@{}}7.43 \\ $\pm 0.03$\end{tabular}  & \begin{tabular}[c]{@{}c@{}}7.45 \\ $\pm 0.03$\end{tabular}  & \begin{tabular}[c]{@{}c@{}}11.02 \\ $\pm 0.09$\end{tabular} & \begin{tabular}[c]{@{}c@{}}11.02 \\ $\pm 0.09$\end{tabular} & \begin{tabular}[c]{@{}c@{}}15.38 \\ $\pm$ 0.24\end{tabular}     & \begin{tabular}[c]{@{}c@{}}8.78 \\ $\pm 0.07$\end{tabular}  & \begin{tabular}[c]{@{}c@{}}46.28 \\ $\pm 0.75$\end{tabular}  \\
                                         &                                            & \textbf{(1.5, 1)}  & \begin{tabular}[c]{@{}c@{}}8.53 \\ $\pm 0.04$\end{tabular}  & \begin{tabular}[c]{@{}c@{}}8.64 \\ $\pm 0.04$\end{tabular}  & \begin{tabular}[c]{@{}c@{}}11.46 \\ $\pm 0.09$\end{tabular} & \begin{tabular}[c]{@{}c@{}}11.46 \\ $\pm 0.09$\end{tabular} & \begin{tabular}[c]{@{}c@{}}14.80 \\ $\pm$ 0.19\end{tabular}     & \begin{tabular}[c]{@{}c@{}}9.68 \\ $\pm 0.07$\end{tabular}  & \begin{tabular}[c]{@{}c@{}}51.46 \\ $\pm 0.84$\end{tabular}  \\
                                         &                                            & \textbf{(1.75, 1)} & \begin{tabular}[c]{@{}c@{}}9.58 \\ $\pm 0.05$\end{tabular}  & \begin{tabular}[c]{@{}c@{}}9.72 \\ $\pm 0.05$\end{tabular}  & \begin{tabular}[c]{@{}c@{}}11.90 \\ $\pm 0.09$\end{tabular} & \begin{tabular}[c]{@{}c@{}}11.90 \\ $\pm 0.09$\end{tabular} & \begin{tabular}[c]{@{}c@{}}15.69 \\ $\pm$ 0.19\end{tabular}     & \begin{tabular}[c]{@{}c@{}}10.58 \\ $\pm 0.07$\end{tabular} & \begin{tabular}[c]{@{}c@{}}56.64 \\ $\pm 0.92$\end{tabular}  \\
                                         &                                            & \textbf{(2, 1)}    & \begin{tabular}[c]{@{}c@{}}10.59 \\ $\pm 0.06$\end{tabular} & \begin{tabular}[c]{@{}c@{}}10.65 \\ $\pm 0.06$\end{tabular} & \begin{tabular}[c]{@{}c@{}}12.34 \\ $\pm 0.09$\end{tabular} & \begin{tabular}[c]{@{}c@{}}12.34 \\ $\pm 0.09$\end{tabular} & \begin{tabular}[c]{@{}c@{}}15.19 \\ $\pm$ 0.15\end{tabular}     & \begin{tabular}[c]{@{}c@{}}11.48 \\ $\pm 0.07$\end{tabular} & \begin{tabular}[c]{@{}c@{}}61.81 \\ $\pm 1.01$\end{tabular}  \\
                                         &                                            & \textbf{(2.25, 1)} & \begin{tabular}[c]{@{}c@{}}11.54 \\ $\pm 0.07$\end{tabular} & \begin{tabular}[c]{@{}c@{}}11.56 \\ $\pm 0.07$\end{tabular} & \begin{tabular}[c]{@{}c@{}}12.78\\ $\pm 0.09$\end{tabular}  & \begin{tabular}[c]{@{}c@{}}12.78 \\ $\pm 0.09$\end{tabular} & \begin{tabular}[c]{@{}c@{}}15.55 \\ $\pm$ 0.15\end{tabular}     & \begin{tabular}[c]{@{}c@{}}12.38 \\ $\pm 0.07$\end{tabular} & \begin{tabular}[c]{@{}c@{}}66.99 \\ $\pm 1.09$\end{tabular}  \\ \hhline{===|=======}
\multicolumn{3}{l|}{\textbf{Avg. Opt. Gap (\%)}}                                                           & -                                                           & \begin{tabular}[c]{@{}c@{}}0.71 \\ $\pm$ 0.11\end{tabular}  & \begin{tabular}[c]{@{}c@{}}22.45 \\ $\pm$ 0.50\end{tabular} & \begin{tabular}[c]{@{}c@{}}22.62 \\ $\pm$ 0.49\end{tabular} & \begin{tabular}[c]{@{}c@{}}85.03 \\ $\pm$ 1.71\end{tabular}     & \begin{tabular}[c]{@{}c@{}}14.05 \\ $\pm$ 0.31\end{tabular} & \begin{tabular}[c]{@{}c@{}}506.20 \\ $\pm$ 4.04\end{tabular} \\ \hhline{===|=======}
\multirow{11}{*}{\textbf{(0.91, 1.11)}}   & \multirow{11}{*}{\textbf{(0.0119, 0.0219)}} & \textbf{(1, 1)}    & \begin{tabular}[c]{@{}c@{}}6.40 \\ $\pm 0.03$\end{tabular}  & \begin{tabular}[c]{@{}c@{}}6.41 \\ $\pm 0.03$\end{tabular}  & \begin{tabular}[c]{@{}c@{}}6.40 \\ $\pm 0.03$\end{tabular}  & \begin{tabular}[c]{@{}c@{}}6.41 \\ $\pm 0.03$\end{tabular}  & \begin{tabular}[c]{@{}c@{}}21.76 \\ $\pm$ 0.44\end{tabular}     & \begin{tabular}[c]{@{}c@{}}7.63 \\ $\pm 0.06$\end{tabular}  & \begin{tabular}[c]{@{}c@{}}42.10 \\ $\pm 0.64$\end{tabular}  \\
                                         &                                            & \textbf{(1.25, 1)} & \begin{tabular}[c]{@{}c@{}}7.63 \\ $\pm 0.04$\end{tabular}  & \begin{tabular}[c]{@{}c@{}}7.67 \\ $\pm 0.04$\end{tabular}  & \begin{tabular}[c]{@{}c@{}}10.63 \\ $\pm 0.08$\end{tabular} & \begin{tabular}[c]{@{}c@{}}10.44 \\ $\pm 0.08$\end{tabular} & \begin{tabular}[c]{@{}c@{}}16.14 \\ $\pm$ 0.26\end{tabular}     & \begin{tabular}[c]{@{}c@{}}8.60 \\ $\pm 0.06$\end{tabular}  & \begin{tabular}[c]{@{}c@{}}47.41 \\ $\pm 0.72$\end{tabular}  \\
                                         &                                            & \textbf{(1.5, 1)}  & \begin{tabular}[c]{@{}c@{}}8.80 \\ $\pm 0.05$\end{tabular}  & \begin{tabular}[c]{@{}c@{}}8.90 \\ $\pm 0.05$\end{tabular}  & \begin{tabular}[c]{@{}c@{}}11.09 \\ $\pm 0.08$\end{tabular} & \begin{tabular}[c]{@{}c@{}}11.09 \\ $\pm 0.08$\end{tabular} & \begin{tabular}[c]{@{}c@{}}15.14 \\ $\pm$ 0.20\end{tabular}     & \begin{tabular}[c]{@{}c@{}}9.57 \\ $\pm 0.06$\end{tabular}  & \begin{tabular}[c]{@{}c@{}}52.72 \\ $\pm 0.80$\end{tabular}  \\
                                         &                                            & \textbf{(1.75, 1)} & \begin{tabular}[c]{@{}c@{}}9.87 \\ $\pm 0.06$\end{tabular}  & \begin{tabular}[c]{@{}c@{}}9.93 \\ $\pm 0.06$\end{tabular}  & \begin{tabular}[c]{@{}c@{}}11.54 \\ $\pm 0.08$\end{tabular} & \begin{tabular}[c]{@{}c@{}}11.54 \\ $\pm 0.08$\end{tabular} & \begin{tabular}[c]{@{}c@{}}15.15 \\ $\pm$ 0.17\end{tabular}     & \begin{tabular}[c]{@{}c@{}}10.54 \\ $\pm 0.07$\end{tabular} & \begin{tabular}[c]{@{}c@{}}58.02 \\ $\pm 0.88$\end{tabular}  \\
                                         &                                            & \textbf{(2, 1)}    & \begin{tabular}[c]{@{}c@{}}10.94 \\ $\pm 0.06$\end{tabular} & \begin{tabular}[c]{@{}c@{}}11.08 \\ $\pm 0.07$\end{tabular} & \begin{tabular}[c]{@{}c@{}}11.99 \\ $\pm 0.08$\end{tabular} & \begin{tabular}[c]{@{}c@{}}11.99 \\ $\pm 0.08$\end{tabular} & \begin{tabular}[c]{@{}c@{}}15.58 \\ $\pm$ 0.17\end{tabular}     & \begin{tabular}[c]{@{}c@{}}11.51 \\ $\pm 0.07$\end{tabular} & \begin{tabular}[c]{@{}c@{}}63.33 \\ $\pm 0.96$\end{tabular}  \\
                                         &                                            & \textbf{(2.25, 1)} & \begin{tabular}[c]{@{}c@{}}11.74 \\ $\pm 0.07$\end{tabular} & \begin{tabular}[c]{@{}c@{}}11.98 \\ $\pm 0.07$\end{tabular} & \begin{tabular}[c]{@{}c@{}}12.44 \\ $\pm 0.08$\end{tabular} & \begin{tabular}[c]{@{}c@{}}12.44 \\ $\pm 0.08$\end{tabular} & \begin{tabular}[c]{@{}c@{}}15.84 \\ $\pm$ 0.16\end{tabular}     & \begin{tabular}[c]{@{}c@{}}12.48 \\ $\pm 0.07$\end{tabular} & \begin{tabular}[c]{@{}c@{}}68.64 \\ $\pm 1.04$\end{tabular}  \\ \hhline{===|=======}
\multicolumn{3}{l|}{\textbf{Avg. Opt. Gap (\%)}}                                                           & -                                                           & \begin{tabular}[c]{@{}c@{}}1.08 \\ $\pm$ 0.15\end{tabular}  & \begin{tabular}[c]{@{}c@{}}16.46 \\ $\pm$ 0.42\end{tabular} & \begin{tabular}[c]{@{}c@{}}16.07 \\ $\pm$ 0.41\end{tabular} & \begin{tabular}[c]{@{}c@{}}91.38 \\ $\pm$ 2.19\end{tabular}     & \begin{tabular}[c]{@{}c@{}}9.97 \\ $\pm$ 0.27\end{tabular}  & \begin{tabular}[c]{@{}c@{}}505.54 \\ $\pm$ 3.77\end{tabular} \\ \hhline{===|=======}
\end{tabular}}

\raggedright
\tabnote{Other system parameters are fixed at $\lambda = (1.5,1.5), b=(0,0), C=4, \kappa=(30,30)$.}
\end{table}

\section{Examining the effect of blocking costs}\label{sec:experiment_blocking}
As an additional analysis, we examine the effect of blocking costs on the performance of different policies. We consider a system with $\lambda = (1.5,1.5)$, $\mu = (1,1)$, $a=(0.0103,0.0203)$ and $h=(5,1)$. 

Table \ref{tab:blocking_effect_experiments} summarizes the long-run average cost and the average optimality gap for different policies under different blocking costs. We observe that (1) the ADP policy yields a very small optimality gap, $1.43\%$ on average. (2) When $b=(0,0)$, $h \bar{f}(0)$ and $h \bar{f}(x)$, which prioritizes class 1, performs well, however, as the blocking cost increases, $h \bar{f}(0)$ and $h \bar{f}(x)$, which does not take the blocking cost into account can perform quite poorly. The other benchmark policies do not perform well overall.

\begin{table}[!h]
\centering
\caption{Long-run average costs under the optimal policy, ADP policy and five benchmark policies and varying blocking costs.}
\label{tab:blocking_effect_experiments}
\begin{tabular}{l|ccccccc}
\hline
\textbf{$b$}                & \textbf{Optimal}                                            & \textbf{ADP}                                                & \textbf{$h \bar{f}(0)$}                                      & \textbf{$h \bar{f}(x)$}                                      & \textbf{\begin{tabular}[c]{@{}c@{}}Max\\ Pressure\end{tabular}} & \textbf{SQF}                                                & \textbf{LQF}                                                   \\ \hline
\textbf{(0, 0)}             & \begin{tabular}[c]{@{}c@{}}18.26 \\ $\pm$ 0.11\end{tabular} & \begin{tabular}[c]{@{}c@{}}18.29 \\ $\pm$ 0.12\end{tabular} & \begin{tabular}[c]{@{}c@{}}18.29 \\ $\pm$ 0.12\end{tabular}  & \begin{tabular}[c]{@{}c@{}}18.29 \\ $\pm$ 0.12\end{tabular}  & \begin{tabular}[c]{@{}c@{}}19.97 \\ $\pm$ 0.14\end{tabular}     & \begin{tabular}[c]{@{}c@{}}20.73 \\ $\pm$ 0.10\end{tabular} & \begin{tabular}[c]{@{}c@{}}112.16 \\ $\pm$ 2.08\end{tabular}   \\ \hline
\textbf{(0, 10)}            & \begin{tabular}[c]{@{}c@{}}19.05 \\ $\pm$ 0.13\end{tabular} & \begin{tabular}[c]{@{}c@{}}19.29 \\ $\pm$ 0.13\end{tabular} & \begin{tabular}[c]{@{}c@{}}19.47 \\ $\pm$ 0.14\end{tabular}  & \begin{tabular}[c]{@{}c@{}}19.47 \\ $\pm$ 0.14\end{tabular}  & \begin{tabular}[c]{@{}c@{}}21.10 \\ $\pm$ 0.16\end{tabular}     & \begin{tabular}[c]{@{}c@{}}21.33 \\ $\pm$ 0.12\end{tabular} & \begin{tabular}[c]{@{}c@{}}114.70 \\ $\pm$ 2.14\end{tabular}   \\ \hline
\textbf{(0, 100)}           & \begin{tabular}[c]{@{}c@{}}19.43 \\ $\pm$ 0.13\end{tabular} & \begin{tabular}[c]{@{}c@{}}19.65 \\ $\pm$ 0.13\end{tabular} & \begin{tabular}[c]{@{}c@{}}30.04 \\ $\pm$ 0.35\end{tabular}  & \begin{tabular}[c]{@{}c@{}}30.04 \\ $\pm$ 0.35\end{tabular}  & \begin{tabular}[c]{@{}c@{}}31.34 \\ $\pm$ 0.37\end{tabular}     & \begin{tabular}[c]{@{}c@{}}26.75 \\ $\pm$ 0.27\end{tabular} & \begin{tabular}[c]{@{}c@{}}137.55 \\ $\pm$ 2.66\end{tabular}   \\ \hline
\textbf{(0, 1000)}          & \begin{tabular}[c]{@{}c@{}}19.45 \\ $\pm$ 0.13\end{tabular} & \begin{tabular}[c]{@{}c@{}}19.65 \\ $\pm$ 0.14\end{tabular} & \begin{tabular}[c]{@{}c@{}}135.77 \\ $\pm$ 2.48\end{tabular} & \begin{tabular}[c]{@{}c@{}}135.77 \\ $\pm$ 2.48\end{tabular} & \begin{tabular}[c]{@{}c@{}}133.67 \\ $\pm$ 2.45\end{tabular}    & \begin{tabular}[c]{@{}c@{}}80.94 \\ $\pm$ 1.97\end{tabular} & \begin{tabular}[c]{@{}c@{}}366.05 \\ $\pm$ 7.87\end{tabular}   \\ \hline
\textbf{(10, 10)}           & \begin{tabular}[c]{@{}c@{}}19.05 \\ $\pm$ 0.13\end{tabular} & \begin{tabular}[c]{@{}c@{}}19.47 \\ $\pm$ 0.14\end{tabular} & \begin{tabular}[c]{@{}c@{}}19.47 \\ $\pm$ 0.14\end{tabular}  & \begin{tabular}[c]{@{}c@{}}19.47 \\ $\pm$ 0.14\end{tabular}  & \begin{tabular}[c]{@{}c@{}}21.10 \\ $\pm$ 0.16\end{tabular}     & \begin{tabular}[c]{@{}c@{}}21.34 \\ $\pm$ 0.12\end{tabular} & \begin{tabular}[c]{@{}c@{}}117.82 \\ $\pm$ 2.21\end{tabular}   \\ \hline
\textbf{(100, 100)}         & \begin{tabular}[c]{@{}c@{}}20.08 \\ $\pm$ 0.14\end{tabular} & \begin{tabular}[c]{@{}c@{}}20.16 \\ $\pm$ 0.14\end{tabular} & \begin{tabular}[c]{@{}c@{}}30.04 \\ $\pm$ 0.35\end{tabular}  & \begin{tabular}[c]{@{}c@{}}30.04 \\ $\pm$ 0.35\end{tabular}  & \begin{tabular}[c]{@{}c@{}}31.34 \\ $\pm$ 0.37\end{tabular}     & \begin{tabular}[c]{@{}c@{}}26.90 \\ $\pm$ 0.27\end{tabular} & \begin{tabular}[c]{@{}c@{}}168.72 \\ $\pm$ 3.36\end{tabular}   \\ \hline
\textbf{(1000, 1000)}       & \begin{tabular}[c]{@{}c@{}}24.39 \\ $\pm$ 0.26\end{tabular} & \begin{tabular}[c]{@{}c@{}}25.33 \\ $\pm$ 0.31\end{tabular} & \begin{tabular}[c]{@{}c@{}}135.77 \\ $\pm$ 2.48\end{tabular} & \begin{tabular}[c]{@{}c@{}}135.77 \\ $\pm$ 2.48\end{tabular} & \begin{tabular}[c]{@{}c@{}}133.67 \\ $\pm$ 2.45\end{tabular}    & \begin{tabular}[c]{@{}c@{}}82.44 \\ $\pm$ 1.97\end{tabular} & \begin{tabular}[c]{@{}c@{}}677.72 \\ $\pm$ 14.93\end{tabular}  \\ \hhline{=|=======}
\textbf{Avg. Opt. Gap (\%)} & -                                                           & \begin{tabular}[c]{@{}c@{}}1.43 \\ $\pm$ 0.02\end{tabular}  & \begin{tabular}[c]{@{}c@{}}166.21 \\ $\pm$ 4.92\end{tabular} & \begin{tabular}[c]{@{}c@{}}166.21 \\ $\pm$ 4.92\end{tabular} & \begin{tabular}[c]{@{}c@{}}169.08 \\ $\pm$ 4.75\end{tabular}    & \begin{tabular}[c]{@{}c@{}}94.75 \\ $\pm$ 2.49\end{tabular} & \begin{tabular}[c]{@{}c@{}}1049.09 \\ $\pm$ 16.68\end{tabular} \\ \hline
\end{tabular}

\tabnote{Other system parameters are fixed at $\mu = (1,1)$, $\lambda = (1.5,1.5)$, $a = (0.0103,0.0203)$, $h=(5,1)$, $C=4$, $\kappa=(30,30)$.}
\end{table}

\section{Case Study}\label{sec:estimation_details}
\subsection{Estimation Details}\label{ap:case_est}
In this section, we provide the details of the estimation in the case study (Section \ref{sec:case_study}). Let $w_i$ and $L_i$ respectively denote the waiting time and rehabilitation length of stay for patient $i$. We are interested in estimating the causal effect of waiting time on the rehabilitation length of stay. Due to potential endogeneity and identification problems (see \cite{gorgulu2022} for details), we utilize an Instrumental Variable (IV) approach \citep{wooldridge2002econometric} to estimate the causal effect. The estimation procedure is similar to \cite{gorgulu2022}, with the key difference being that we estimate the effect of $w_i$ on $L_i$ directly, rather than on $log(L_i)$

Let $IV_i$ denote a valid instrument and $CV_i$ be the vector of control variables containing various information regarding patients clinical status (e.g., age, sex, commorbidity) and some operational factor. We use following set of equations to describe the system:
\begin{align}
    L_i &= \psi_1 CV_i + \beta w_i + u_i \label{eq:2sls_2}\\
    D_i &= \psi_2 CV_i + \gamma IV_i + v_i \label{eq:2sls_1}
\end{align}
where $\psi$ is the vector of control coefficients, $u_i$ and $v_i$ are random errors, and $\beta$ denote the effect of one additional day of waiting on the rehabilitation length of stay.

We conduct separate analysis for Neuro/MSK and Medicine categories and estimate the value of $\beta$ using two-stage least squares method (2SLS). We utilize rehab occupancy at the time that patients receive Alternative Level of Care (ALC) status as the instrument. The details on the data preparation and analysis on the instrument can be found in \cite{gorgulu2022}. 

Note that in our case study, we model the rehabilitation length of stay of patient $i$ as 
\[
L_i = \alpha + \beta (w_i \wedge 9) + \epsilon
\]
where $\alpha$ is the intercept and $\epsilon$ is the random error with mean 0 and standard deviation $\gamma$. After estimating $\beta$, we fix it to its estimated value and we estimate $\alpha$ and $\gamma$ based on a single regression analysis. The estimated parameters can be found in Table \ref{tab:estimated_parameters} in Section \ref{sec:case_study}.
\subsection{Preemptive vs. Non-Preemptive Algorithms}\label{ap:case_pre}
\begin{table}[]
\centering
\caption{Average waiting times of Neuro/MSK and Medicine patients under ADP policies with various holding costs generated using preemptive and non-preemptive versions of the algorithm.}
\label{tab:case_study_wait_sensitivity_h}
\begin{tabular}{c|cc|cc}
\hline
        & \multicolumn{2}{c|}{\textbf{Preemptive}} & \multicolumn{2}{c}{\textbf{Non-preemptive}} \\ \hline
$h$     & \textbf{Neuro/MSK}       & \textbf{Medicine}       & \textbf{Neuro/MSK}       & \textbf{Medicine}       \\ \hline
$1,1.8$ &15.983                    &  0.246                 &  0.010                   &   0.214                \\
$1,1.4$ & 0.010                    & 0.215                   &   0.010                  &   0.214                 \\
$1,1$   & 0.010                    & 0.213                   &  0.010                   &   0.214                 \\
$1.4,1$ & 0.010                    & 0.211                   &  0.010                  & 0.214                   \\
$1.8,1$ & 0.010                    & 0.213                 &    0.010               &      0.214              \\ 
\hline
\end{tabular}
\end{table}

Recall from our discussion in Section \ref{sec:nonpreemptive} that utilizing the non-preemptive version of our algorithm requires considering a much larger state space, which adds complexity. As a result, preemptive ADP in general leads to faster policy learning. However, if the real system is non-preemptive, using the policy learned by the preemptive ADP may lead to poor performance.
We next compare the performance of the policies generated by the preemptive version of the ADP algorithm and the non-preemptive version. Note that both policies are evaluated based on the real non-preemptive system with wait-dependent service slowdown and log-normally distributed service times. Table \ref{tab:case_study_wait_sensitivity_h} illustrates the average waiting times for different holding costs. We observe that the preemptive algorithm performs well when holding costs are relatively close to each other and it fails when the differences are large (i.e., the $h=(1.8,1)$ case). Intuitively, this is because preemptive ADP, unlike its non-preemptive counterpart, allows for substantial changes in capacity allocations from one decision epoch to the next. Consequently, in the preemptive setting, short-term cost reduction tends to have a greater impact than in the non-preemptive setting.

\end{document}